\documentclass[9pt]{amsart}
\usepackage{amsmath,amsfonts,amsthm}
\usepackage{color}
\usepackage{verbatim}
\usepackage{eepic,epic}
\usepackage{epsfig,subfigure,epstopdf}
\usepackage[]{graphics}
\usepackage[]{graphics}
\usepackage{url}
\graphicspath{{./pics/}}

\newtheorem{lemma}{Lemma}[section]

\numberwithin{equation}{section}

\newcount\icount

\def\DD#1#2{\icount=#1
  \ifnum\icount<1
  \,_{ 0}\kern -.1em D^{#2}_{\kern -.1em x}
  \else
  \,_{x}\kern -.2em D^{#2}_1
  \fi
}

\def\DDRI#1#2{\icount=#1
  \ifnum\icount<1
  \,_{-\infty}^{\kern 1em R}\kern -.2em D^{#2}_{\kern -.1em x}
  \else
  \,_{x}^R \kern -.2em D^{#2}_\infty
  \fi
}

\def\DDR#1#2{\icount=#1
  \ifnum\icount<1
 _{0}^{ \kern -.1em R} \kern -.2em D^{#2}_{\kern -.1em x}
  \else
 _{x}^{ \kern -.1em R} \kern -.2em D^{#2}_{\kern -.1em 1}
  \fi
}

\def\DDCI#1#2{\icount=#1
  \ifnum\icount<1
  \,_{-\infty}^{\kern 1em C}  \kern -.2em D^{#2}_{\kern -.1em x}
  \else
  \,_{x}^C \kern -.2em  D^{#2}_\infty
  \fi
}

\def\DDC#1#2{\icount=#1
  \ifnum\icount<1
  \,_{0}^C \kern -.2em  D^{#2}_{\kern -.1em x}
  \else
  \,_{x}^C \kern -.2em D^{#2}_1
  \fi
}

\def\Hdi#1#2{\icount=#1
  \ifnum\icount<1
  \widetilde H_{L}^{#2}\II
  \else
  \widetilde H_{R}^{#2}\II
  \fi
}

\begin{document}
\title[Inverse Problems for FDEs]
{A Tutorial on Inverse Problems for Anomalous Diffusion Processes}
\author {Bangti Jin \and William Rundell}

\address{Department of Computer Science, University College London, Gower Street,
London WC1E 6BT, UK ({bangti.jin@gmail.com})}
\address{Department of Mathematics,
Texas A\&M University, College Station, TX 77843-3368, USA ({{rundell@math.tamu.edu}})}

\begin{abstract}
Over the last two decades, anomalous diffusion processes in which the mean squares variance grows slower or faster
than that in a Gaussian process have found many applications. At a macroscopic level, these processes are adequately
described by fractional differential equations, which involves fractional derivatives in time or/and space. The
fractional derivatives describe either history mechanism or long range interactions of particle motions at a
microscopic level. The new physics can change dramatically the behavior of the forward problems. For example, the
solution operator of the time fractional diffusion diffusion equation has only limited smoothing property, whereas
the solution for the space fractional diffusion equation may contain weakly singularity. Naturally
one expects that the new physics will impact related inverse problems in terms of uniqueness, stability, and degree
of ill-posedness. The last aspect is especially important from a practical point of view, i.e., stably reconstructing
the quantities of interest.

In this paper, we employ a formal analytic and numerical way, especially the two-parameter Mittag-Leffler function and singular
value decomposition, to examine the degree of ill-posedness of several ``classical'' inverse problems for fractional differential
equations involving a Djrbashian-Caputo fractional derivative in either time or space, which represent the fractional analogues
of that for classical integral order differential equations. We discuss four inverse problems, i.e., backward fractional
diffusion, sideways problem, inverse source problem and inverse potential problem for time fractional diffusion, and inverse
Sturm-Liouville problem, Cauchy problem, backward fractional diffusion and sideways problem for space fractional diffusion. It
is found that contrary to the wide belief, the influence of anomalous diffusion on the degree of ill-posedness is not definitive: it can either significantly
improve or worsen the conditioning of related inverse problems, depending crucially on the specific type of given data and quantity
of interest. Further, the study exhibits distinct new features of ``fractional'' inverse problems, and a partial list of surprising
observations is given below.
\begin{itemize}
  \item[(a)] Classical backward diffusion is exponentially ill-posed, whereas time fractional backward diffusion is only mildly
  ill-posed in the sense of norms on the domain and range spaces. However, this does not imply that the latter always allows a more effective reconstruction.
  \item[(b)] Theoretically, the time fractional sideways problem is severely ill-posed like its classical counterpart, but
  numerically can be nearly well-posed.
  \item[(c)] The classical Sturm-Liouville problem requires two pieces of spectral data to uniquely determine a general potential, but
  in the fractional case, one single Dirichlet spectrum may suffice.
  \item[(d)] The space fractional sideways problem can be far more or far less ill-posed than the classical counterpart,
  depending on the location of the lateral Cauchy data.
\end{itemize}
In many cases, the precise mechanism of these surprising observations is unclear, and awaits further analytical and numerical
exploration, which requires new mathematical tools and ingenuities. Further, our findings indicate fractional diffusion inverse
problems also provide an excellent case study in the differences between theoretical ill-conditioning involving domain and range
norms and the numerical analysis of a finite-dimensional reconstruction procedure.  Throughout we will also describe known analytical and
numerical results in the literature.\\
\textbf{Keywords}: fractional inverse problem; fractional differential equation; anomalous diffusion; Djrbashian-Caputo
fractional derivative; Mittag-Leffler function.
\end{abstract}

\maketitle

\section{Introduction}\label{sec:intro}

Diffusion is one of the most prominent transport mechanisms found in nature. At a microscopic level, it is
related to the random motion of individual particles, and the use of the Laplace operator and the
first-order derivative in the canonical diffusion model rests on a Gaussian process assumption on the
particle motion, after Albert Einstein's groundbreaking work \cite{Einstein:1905}. Over the last two decades a large body of
literature has shown that anomalous diffusion models in which the mean square variance grows faster
(superdiffusion) or slower (subdiffusion) than that in a Gaussian process under certain circumstances
can offer a superior fit to experimental data (see the comprehensive reviews \cite{MetzlerKlafter:2000,
SokolovKlafterBlumen:2002,BerkowitzCortisDentzScher:2006,Mainardi:2010} for physical background and
practical applications). For example, anomalous diffusion is often observed in materials with memory, e.g.,
viscoelastic materials, and heterogeneous media, such as soil, heterogeneous aquifer, and underground fluid flow. At a microscopic
level, the subdiffusion process can be described by a continuous time random walk \cite{MontrollWeiss:1965}, where the waiting time
of particle jumps follows some heavy tailed distribution, whereas the superdiffusion process can be
described by L\'{e}vy flights or L\'{e}vy walk, where the length of particle jumps follows some heavy
tailed distribution, reflecting the long-range interactions among particles. Following the aforementioned micro-macro
correspondence, the macroscopic counterpart of a continuous time random walk is a differential equation
with a fractional derivative in time, and that for a L\'{e}vy flight is a differential equation with a fractional
derivative in space. We will refer to these two cases as time fractional diffusion and space fractional diffusion,
respectively, and it is generically called a fractional derivative equation (FDE) below.
In general the fractional derivative can appear in both time and space variables.

Next we give the mathematical model. Let $\Omega=(0,1)$ be the unit interval. Then a general one-dimensional
FDE is given by
\begin{equation}\label{eqn:fde}
  \partial_t^\alpha u - {\DDC 0 \beta } u + q u = f\quad (x,t) \in \Omega\times(0,T),
\end{equation}
where $T>0$ is a fixed time, and it is equipped with suitable boundary and initial conditions. The fractional
orders $\alpha\in(0,1)$ and $\beta\in (1,2)$ are related to the parameters specifying the large-time behavior
of the waiting-time distribution or long-range behavior of the particle jump distribution. For example, in
hydrological studies, the parameter $\beta$ is used to characterize the heterogeneity of porous medium
\cite{ClarkeMeerschaertWheatcraft:2005}.  In theory, these
parameters can be determined from the underlying stochastic model, but often in practice, they are determined
from experimental data \cite{HatanoHatano:1998,HatanoNakagawaWangYamamoto:2013,LiYamamoto:2014b}. The notation $\partial_t^\alpha
= {\,_{0}^C \kern -.2em  D^{\alpha}_{\kern -.1em t}}$ is the Djrbashian-Caputo derivative operator of order $\alpha\in(0,1)$
in the time variable $t$, and $\DDC 0 \beta$ denotes the Djrbashian-Caputo derivative of order $\beta\in(1,2)$ in the space variable $x$. For a real
number $n-1<\gamma<n$, $n\in\mathbb{N}$, and $f\in H^n(0,1)$, the left-sided Djrbashian-Caputo derivative
$\DDC 0 \gamma f$ of order $\gamma$ is defined by \cite[pp. 91]{KilbasSrivastavaTrujillo:2006}
\begin{equation}\label{eqn:Drjbashian-Caputo}
  \DDC 0 \gamma  f = \frac{1}{\Gamma(n-\gamma)}\int_0^x (x-s)^{n-1-\gamma}f^{(n)}(s)ds,
\end{equation}
where $\Gamma(z)$ denotes Euler's Gamma function defined by
\begin{equation*}
 \Gamma(z) =\int_0^\infty s^{z-1}e^{-s}ds,\quad \ \Re(z)>0.
\end{equation*}
The Djrbashian-Caputo derivative was first introduced by Armenian mathematician Mkhitar M. Djrbashian for studies
on space of analytical functions and integral transforms in 1960s (see \cite{Dzharbashyan:1964,Djrbashian:1989,
Djrbashian:1993} for surveys on related works). Italian geophysicist Michele Caputo independently proposed the use of the
derivative for modeling the dynamics of viscoelastic materials in 1967 \cite{Caputo:1967}. We note that there are
several alternative (and different) definitions of fractional derivatives, notably the Riemann-Liouville fractional
derivative, which formally is obtained from \eqref{eqn:Drjbashian-Caputo} by interchanging the order of integration
and differentiation, i.e., the left-sided Riemann-Liouville fractional derivative
$\DDR 0 \gamma f$ of order $\gamma\in (n-1,n)$, $n\in\mathbb{N}$, is defined by \cite[pp. 70]{KilbasSrivastavaTrujillo:2006}
\begin{equation*}
  \DDR 0 \gamma f = \frac{d^n}{dx^n} \frac{1}{\Gamma(n-\gamma)}\int_0^x (x-s)^{n-1-\gamma}f(s)ds.
\end{equation*}
In this work, we shall focus mostly on the Djrbashian-Caputo derivative since it allows a convenient treatment
of the boundary and initial conditions.

Under certain regularity assumption on the functions, with an integer order $\gamma$, the Djrbashian-Caputo and Riemann-Liouville
derivatives both recover the usual integral order derivative (see for example \cite[pp. 100]{NakagawaSakamotoYamamoto:2010} for
the Djrbashian-Caputo case). For example, with $\alpha=1$ and $\beta=2$, the Djrbashian-Caputo fractional derivatives
$\partial_t^\alpha u$ and $\DDC 0 \alpha u$ coincide with the usual first- and second-order derivatives $\frac{\partial u}{
\partial t}$ and $\frac{\partial^2u}{\partial x^2}$, respectively, for which the model \eqref{eqn:fde} recovers the standard
one-dimensional diffusion equation, and thus generally the model \eqref{eqn:fde} is regarded as a fractional counterpart. The Djrbashian-Caputo derivative (and
many others) is an integro-differential operator, and thus it is nonlocal in nature. As a consequence, many useful rules, e.g.,
product rule and integration by parts, from PDEs are either invalid or require significant modifications. The nonlocality underlies
most analytical and numerical challenges associated with the model \eqref{eqn:fde}. It significantly complicates the mathematical
and numerical analysis of the model, including relevant inverse problems.

In a fractional model, there are a number of parameters, e.g., fractional order(s), diffusion and potential coefficients
(when using a second-order elliptic operator in space), initial condition, source term, boundary conditions and domain geometry,
that cannot be measured/specified directly, and have to be inferred indirectly from measured data. Typically, the data is the forward
solution restricted to either the boundary or the interior of the physical domain. This gives rise to a large variety
of inverse problems for FDEs, which have started to attract much attention in recent years, since the pioneering work
\cite{ChengNakagawaYamamotoYamazaki:2009}. An interesting question is how the nonlocal physics behind anomalous diffusion
processes will influence the behavior of related inverse problems, e.g., uniqueness, stability, and the degree of
ill-posedness. The degree of ill-posedness is especially important for developing practical numerical reconstruction
procedures. There is a now well known example of backward fractional diffusion, i.e., recovering the initial condition
in a time fractional diffusion equation from the final time data, which is only mildly ill-posed, instead of severely
ill-posed for the classical backward diffusion problem. In some sense, this example has led to the belief that ``fractionalizing''
inverse problems can always mitigate the degree of ill-posedness, and thus allows a better chance of an accurate numerical reconstruction.

In this paper, we examine the degree of ill-posedness of ``fractional'' inverse problems from a formal analytic and
numerical point of view, and contrast their numerical stability properties with their classical, that is, the Gaussian diffusion counterparts,
for which there are many deep analytical results \cite{PrilepkoOrlovskyVasin:2000,Isakov:2006,Isakov:1999}. Specifically,
we revisit a number of ``classical'' inverse problems for the FDEs, e.g., the backward diffusion problem,
sideways diffusion problem and inverse source problem, and numerically exhibit their degree of ill-posedness. These examples
indicate that the answer to the aforementioned question is not definitive: it depends crucially on the type (unknown and data)
of the inverse problem we look at, and the nonlocality of the problem (fractional derivative) can either greatly improve or
worsen the degree of ill-posedness.

The mathematical theory of inverse problems for FDEs is still in its infancy, and thus in this work, we only discuss the
topic formally to give a flavor of inverse problems for FDEs -- our goal is to give insight rather than to pursue
an in-depth analysis.  The technical developments that are available we leave to the references cited. In addition,
known theoretical results and computational techniques in the literature will be briefly described, which however
are not meant to be exhaustive. The rest of the paper
is organized as follows. In Section \ref{sec:prelim} we review two special functions, i.e., Mittag-Leffler
function and Wright function, and their basic properties. The Mittag-Leffler function plays an extremely important
role in understanding anomalous diffusion processes. We also recall the basic tool -- singular value decomposition --
for analyzing discrete inverse problems. Then in Section \ref{sec:time} we study several inverse problems
for FDEs with a time fractional derivative, including backward diffusion, inverse source problem, sideways problem
and inverse potential problem. In Section \ref{sec:space} we consider inverse problems for FDEs with a space fractional
derivative, including the inverse Sturm-Liouville problem, Cauchy problem, backward diffusion and
sideways problem. In the appendices, we give the implementation details of the computational methods for solving
the time- and space fractional differential equations. These methods are employed throughout for computing the
forward map (unknown-to-measurement map) so as to gain insight into related inverse problems. Throughout the
notation $c$, with or without a subscript, denote a generic constant, which may differ at different
occurrences, but it is always independent of the unknown of interest.

\section{Preliminaries}\label{sec:prelim}
We recall two important special functions, Mittag-Leffler function and Wright function, and one useful tool for
analyzing discrete ill-posed problems, singular value decomposition.
\subsection{Mittag-Leffler function}
We shall use extensively the two-parameter Mittag-Leffler function $E_{\alpha,\beta}(z)$ (with $\alpha>0$ and $\beta\in\mathbb{R}$)
defined by \cite[equation (1.8.17), pp. 40]{KilbasSrivastavaTrujillo:2006}
\begin{equation}\label{eqn:mitlef}
  E_{\alpha,\beta}(z) = \sum_{k=0}^\infty \frac{z^k}{\Gamma(k\alpha+\beta)} \ \ z\in\mathbb{C}.
\end{equation}
This function with $\beta=1$ was first introduced by G\"{o}sta Mittag-Leffler in 1903 \cite{Mittag-Leffler:1903}
and then generalized by others \cite{Humbert:1953,Agarwal:1953}. It can be verified directly that
\begin{equation*}
  \begin{aligned}
    E_{1,1}(z) = e^z, \quad
    E_{2,1}(z) = \cosh \sqrt{z},\quad
    E_{2,2}(z) = \frac{\sinh \sqrt{z}}{\sqrt{z}}.
  \end{aligned}
\end{equation*}
Hence it represents a generalization of the exponential function in that $E_{1,1}(z)=e^z$. The Mittag-Leffler function $E_{\alpha,
\beta}(z)$ is an entire function of $z$ with order $\alpha^{-1}$ and type 1 \cite[pp. 40]{KilbasSrivastavaTrujillo:2006}. Further,
the function $E_{\alpha,1}(-t)$ is completely monotone on the positive real axis $\mathbb{R}^+$ \cite{Pollard:1948}, and thus
it is positive on $\mathbb{R}^+$; see also \cite{Schneider:1996} for extension to the two-parameter Mittag-Leffler
function $E_{\alpha,\beta}(z)$. It appears in the solution representation for FDEs: The functions $E_{\alpha,1} (-\lambda t^\alpha)$ and $t^{\alpha-1}
E_{\alpha,\alpha}(-\lambda t^\alpha)$ appear in the kernel of the time fractional diffusion problem with initial data and the right
hand side, respectively, and also are eigenfunctions to the fractional Sturm-Liouville problem with a zero potential, cf. Section \ref{ssec:eig}.

In our discussions, the asymptotic behavior of the function $E_{\alpha,\beta}(z)$ will play a crucial role.
It satisfies the following exponential asymptotics  \cite[pp. 43, equations (2.8.17) and (2.8.18)]{KilbasSrivastavaTrujillo:2006},
which was first derived by Djrbashian \cite{Dzharbashyan:1964}, and refined by many researchers \cite{WongZhao:2002,Paris:2002}.
\begin{lemma}\label{lem:mitlef}
Let $\alpha\in(0,2)$, $\beta\in\mathbb{R}$, and $\mu\in(\alpha\pi/2,\min(\pi,\alpha\pi))$. Then for $N\in\mathbb{N}$
\begin{equation*}
  \begin{aligned}
    E_{\alpha,\beta}(z) &= \frac{1}{\alpha}z^{(1-\beta)/\alpha}e^{z^{1/\alpha}}-\sum_{k=1}^N\frac{1}{\Gamma(\beta-\alpha k)}\frac{1}{z^k} + O(\frac{1}{z^{N+1}})\ \mbox{with } |z|\to\infty,\ |\mathrm{arg}(z)|\leq \mu,\\
    E_{\alpha,\beta}(z) & = -\sum_{k=1}^N\frac{1}{\Gamma(\beta-\alpha k)}\frac{1}{z^k} + O(\frac{1}{z^{N+1}})\ \mbox{with } |z|\to\infty,\ \mu \leq |\mathrm{arg}(z)|\leq \pi.
  \end{aligned}
\end{equation*}
\end{lemma}

\begin{figure}[hbt!]
  \centering
  \begin{tabular}{cc}
  \includegraphics[trim = 1cm .1cm 1cm 0.5cm, clip = true, width=7.5cm]{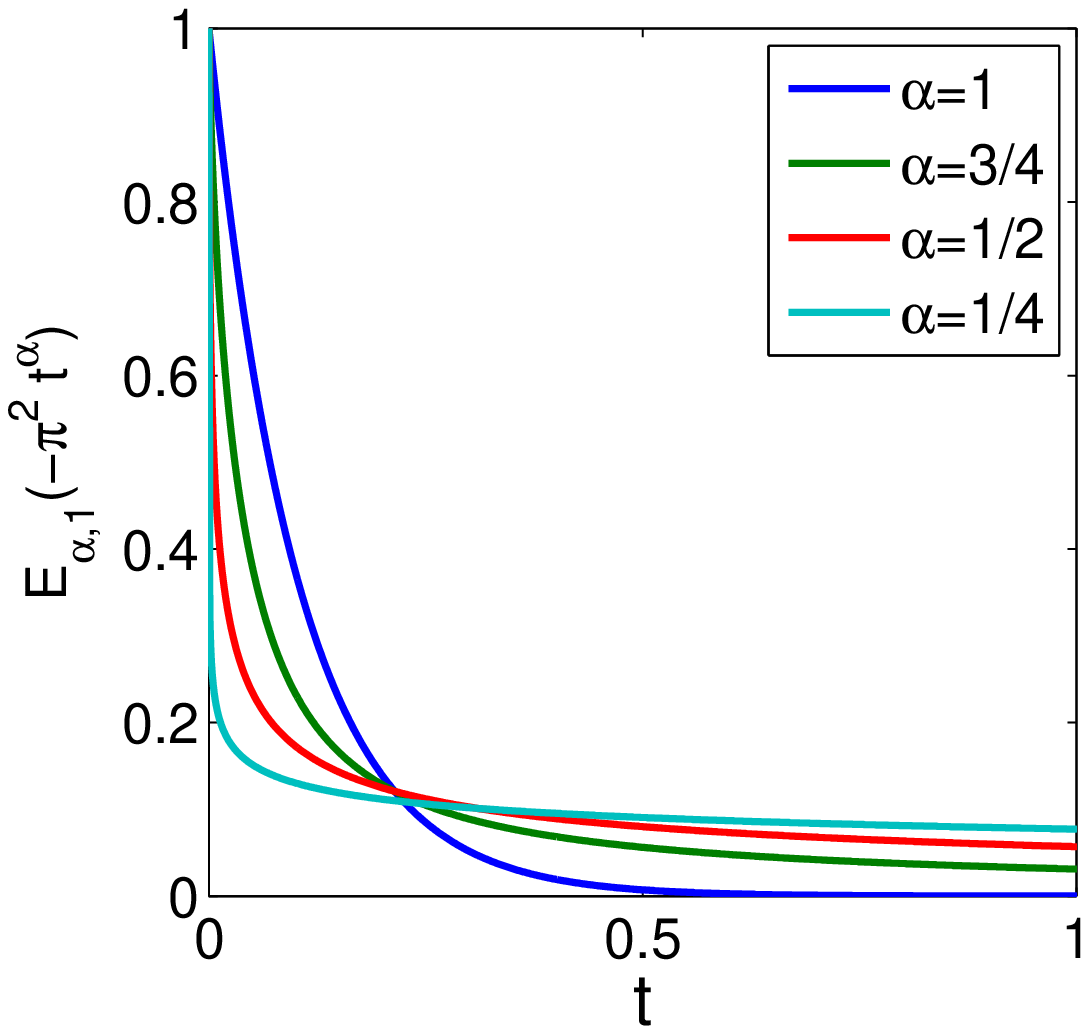}&
  \includegraphics[trim = 1cm .1cm 1cm 0.5cm, clip = true, width=7.5cm]{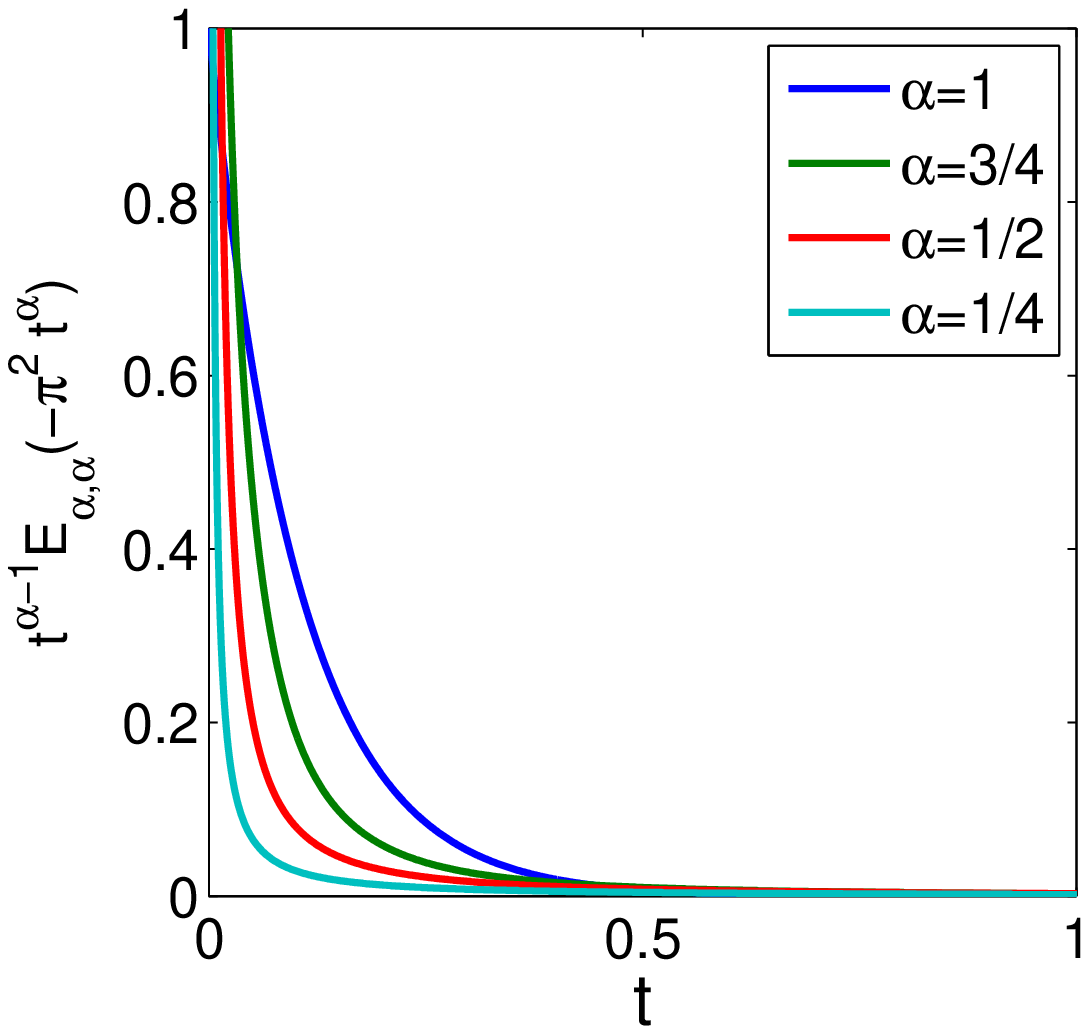}\\
   (a) $E_{\alpha,1}(-\pi^2t^\alpha)$ & (b) $t^{\alpha-1}E_{\alpha,\alpha}(-\pi^2 t^\alpha)$
  \end{tabular}
  \caption{The profiles of Mittag-Leffler functions (a) $E_{\alpha,1}(-\pi^2t^\alpha)$ and (b)
  $t^{\alpha-1}E_{\alpha,\alpha}(-\pi^2 t^\alpha)$. The value $\pi^2$ is the first Dirichlet eigenvalue
  of the negative Laplacian on the interval $(0,1)$.}\label{fig:mitlef}
\end{figure}

From these asymptotics, the Mittag-Leffler function $E_{\alpha,\beta}(z)$ decays only linearly on the negative real axis
$\mathbb{R}^-$, which is much slower than the exponential decay for the exponential function $e^{z}$. However, on the
positive real axis $\mathbb{R}^+$, it grows exponentially, and the growth rate increases with the fractional
order $0<\alpha<2$. To illustrate the distinct feature, we plot the functions $E_{\alpha,1}(-\pi^2 t^\alpha)$ and
$t^{\alpha-1}E_{\alpha,\alpha}(-\pi^2 t^\alpha)$ in Fig. \ref{fig:mitlef} for several different $\alpha$ values, where
$\lambda=\pi^2$ is the first Dirichlet eigenvalue of the negative Laplacian on the unit interval $\Omega=(0,1)$; see
Appendix \ref{app:mitlef} for further details on the computation of the Mittag-Leffler function.
Fig. 1(a) can be viewed as the time evolution of $u(1/2,t)$, where $\partial^\alpha_t u - u_{xx} = 0$ with $u(0,t)=u(1,t)
= 0$, and initial data $u_0(x) = \sin\pi x$ (the lowest Fourier eigenmode). The slow decay behavior at large time is clearly
observed. In particular, at $t=1$, the function $E_{\alpha,1}(-\pi^2t)$ still takes values distinctly away from zero for
any $0<\alpha<1$, whereas the exponential function $e^{-\pi^2 t}$ almost vanishes identically. In contrast, for $t$ close
to zero, the picture is reversed: the Mittag-Leffler function $E_{\alpha,1}(-\pi^2t)$ decays much faster than the exponential
function $e^{-\pi^2t}$. The drastically different behavior of the function $E_{\alpha,1}(-z)$, in comparison with the
exponential function $e^{-z}$, explains many unusual phenomena with inverse problems for FDEs to be described below.
According to the exponential asymptotics, the function $E_{\alpha,\alpha}(z)$ decays faster on the negative real axis
$\mathbb{R}^-$, since $1/\Gamma(0)=0$, i.e., the first term in the expansion vanishes. This is confirmed numerically in Fig.
\ref{fig:mitlef}(b). Even though not shown, it is noted that the function $E_{\alpha,\alpha}(z)$ decays only quadratically on
the negative real axis $\mathbb{R}^-$ for $\alpha\in(0,1)$ or $\alpha\in(1,2)$, which is asymptotically much slower than the exponential decay.

The distribution of zeros of the Mittag-Leffer function $E_{\alpha,\beta}(z)$ is of immense interest, especially in the related
Sturm-Liouville problem; see Section \ref{ssec:eig} below. The case of $\beta=1$ was first studied by Wiman \cite{Wiman:1905}. It was revisited by
Djrbashian \cite{Dzharbashyan:1964}, and many deep results were derived, especially for the case of $\alpha=2$. There are many further refinements
\cite{Sedletskii:1994}; see \cite{PopovSedletskii:2013} for an updated account.

\subsection{Wright function}
The Wright function $W_{\rho,\mu}(z)$ is defined by \cite{Wright:1933,Wright:1940}
\begin{equation*}
W_{\rho,\mu}(z) = \sum_{k=0}^\infty \frac{z^k}{k!\Gamma(\rho k + \mu)},
\qquad \mu,\ \rho \in \mathbb{R},\ \rho>-1,\quad z \in \mathbb{C}.
\end{equation*}
This is an entire function of order $1/(1+\rho)$ \cite[Theorem 2.4.1]{GorenfloLuchkoMainardi:1999}. It was first
introduced in connection with a problem in number theory by  Edward M. Wright, and revived in recent years since it appears
as the fundamental solution for FDEs \cite{Mainardi:1996}. The Wright function $W_{\rho,\mu}(z)$ has the following the asymptotic
expansion in one sector containing the negative real axis $\mathbb{R}^-$ \cite[Theorem 3.2]{Luchko:2008}. Like before, the exponential
asymptotics can be used to deduce the distribution of its zeros \cite{Luchko:2000}.
\begin{lemma}\label{lem:wright}
Let $ -1<\rho<0$, $y=-z$, $\mathrm{arg}(z) \leq \pi$, $-\pi<\mathrm{arg}(y)\leq \pi$,
$|\mathrm{arg}(y)|\leq \min(3\pi(1+\rho)/2,\pi)-\epsilon$, $\epsilon>0$. Then
\begin{equation*}
  W_{\rho,\mu}(z) = Y^{1/2-\mu}e^{-Y}\left\{\sum_{m=0}^{M-1}A_mY^{-m}+O(Y^{-M})\right\},\ \ Y\to\infty,
\end{equation*}
where $Y=(1+\rho)((-\rho)^{-\rho}y)^{1/(1+\rho)}$ and the coefficients $A_m$, $m=0,1,\ldots$
are defined by the asymptotic expansion
\begin{equation*}
  \begin{aligned}
    \frac{\Gamma(1-\mu-\rho t)}{2\pi(-\rho)^{-\rho t}(1+\rho)^{(1+\rho)(t+1)}\Gamma(t+1)}&= \sum_{m=0}^{M-1}\frac{(-1)^mA_m}{\Gamma((1+\rho)t+\mu+\frac{1}{2}+m)}\\
     &\qquad+O\left(\frac{1}{\Gamma((1+\rho)t+\beta+\tfrac{1}{2}+M)}\right),
  \end{aligned}
\end{equation*}
valid for $\mathrm{arg}(t), \mathrm{arg}(-\rho t)$, and $\mathrm{arg}(1-\mu-\rho t)$ all lying
between $-\pi$ and $\pi$ and $t$ tending to infinity.
\end{lemma}

The Wright function $W_{\rho,\mu}(z)$, $-1<\rho<0$, decays exponentially on the negative real axis $\mathbb{R}^-$,
in a manner similar to the exponential function $e^z$, but at a different decay rate. Its special role
in fractional calculus is underscored by the fact that it forms the free-space fundamental solution
$K_\alpha(x,t)$ to the one-dimensional time fractional diffusion equation \cite{Mainardi:1996} by
\begin{equation}\label{eqn:fundsol}
  K_\alpha(x,t) = \frac{1}{2t^\frac{\alpha}{2}}W_{-\frac{\alpha}{2},1-\frac{\alpha}{2}}\left(-|x|/t^{\alpha/2}\right).
\end{equation}
The multidimensional case is more complex and involves further special functions, in particular, the Fox H function
\cite{SchneiderWyss:1989,Kochubei:1990}. For $\alpha=1$, the formula \eqref{eqn:fundsol} recovers the familiar
free-space fundamental solution for the one-dimensional heat equation, i.e.,
\begin{equation*}
   K(x,t) = \frac{1}{2\sqrt{\pi t}}e^{-\frac{x^2}{4t}},
\end{equation*}
which is a Gaussian distribution in $x$ for any $t>0$. In the fractional case, the fundamental solution $K_\alpha(x,t)$
exhibits quite different behavior than the heat kernel. To see this, we show the profile of $K_\alpha(x,t)$ in Fig.
\ref{fig:wright} for several $\alpha$ values; see Appendix \ref{app:mitlef} for a brief description of the
computational details. For any $0<\alpha<1$, $K_\alpha(x,t)$ decays slower at a polynomial
rate as the argument $|x|/t^{\alpha/2}$ tends to infinity, i.e., having a long tail, when compared with the Gaussian density.
The long tail profile is one of distinct features of slow diffusion \cite{BerkowitzCortisDentzScher:2006}. Further,
for any $\alpha<1$, the profile is only continuous but not differentiable at $x=0$. The kink at the origin implies
that the solution operator to time fractional diffusion may only have a limited smoothing property.

\begin{figure}[hbt!]
  \centering
  \begin{tabular}{cc}
  \includegraphics[trim = 2cm .1cm 2cm 0.5cm, clip = true, width=7.5cm]{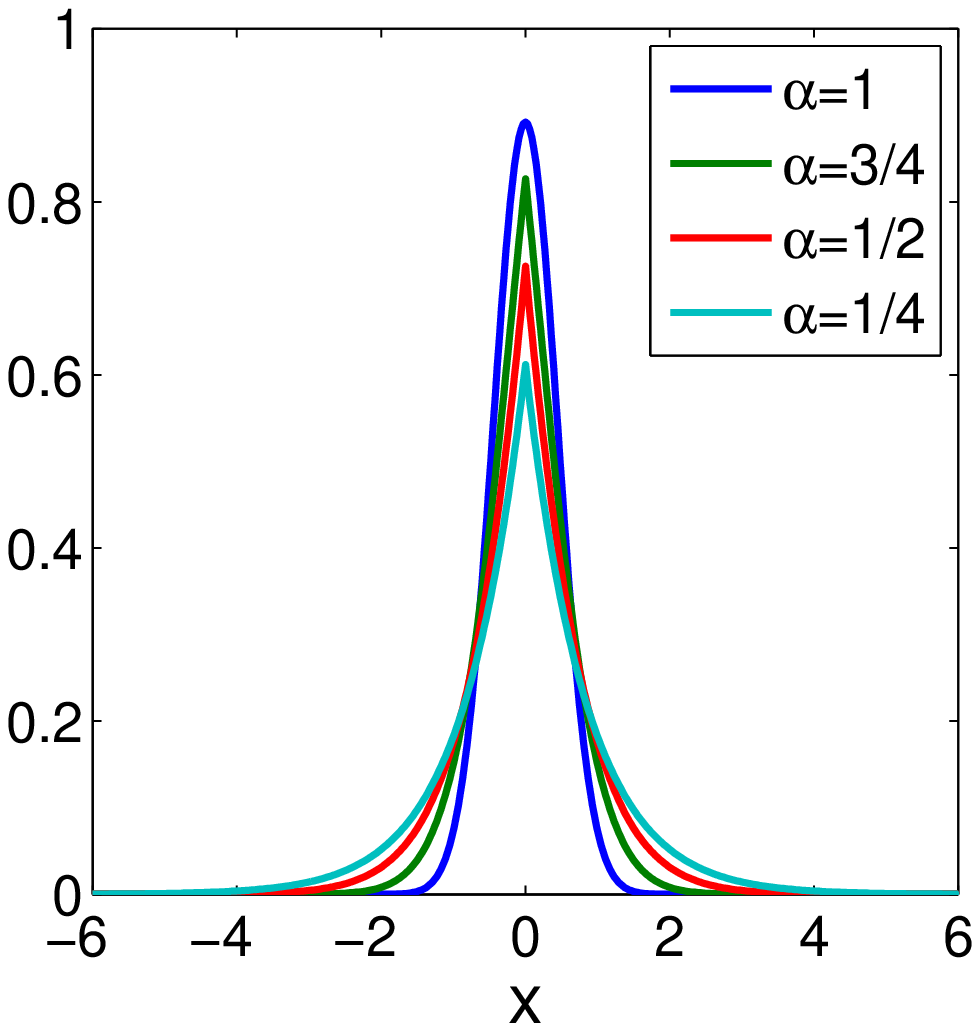}&
  \includegraphics[trim = 2cm .1cm 2cm 0.5cm, clip = true, width=7.5cm]{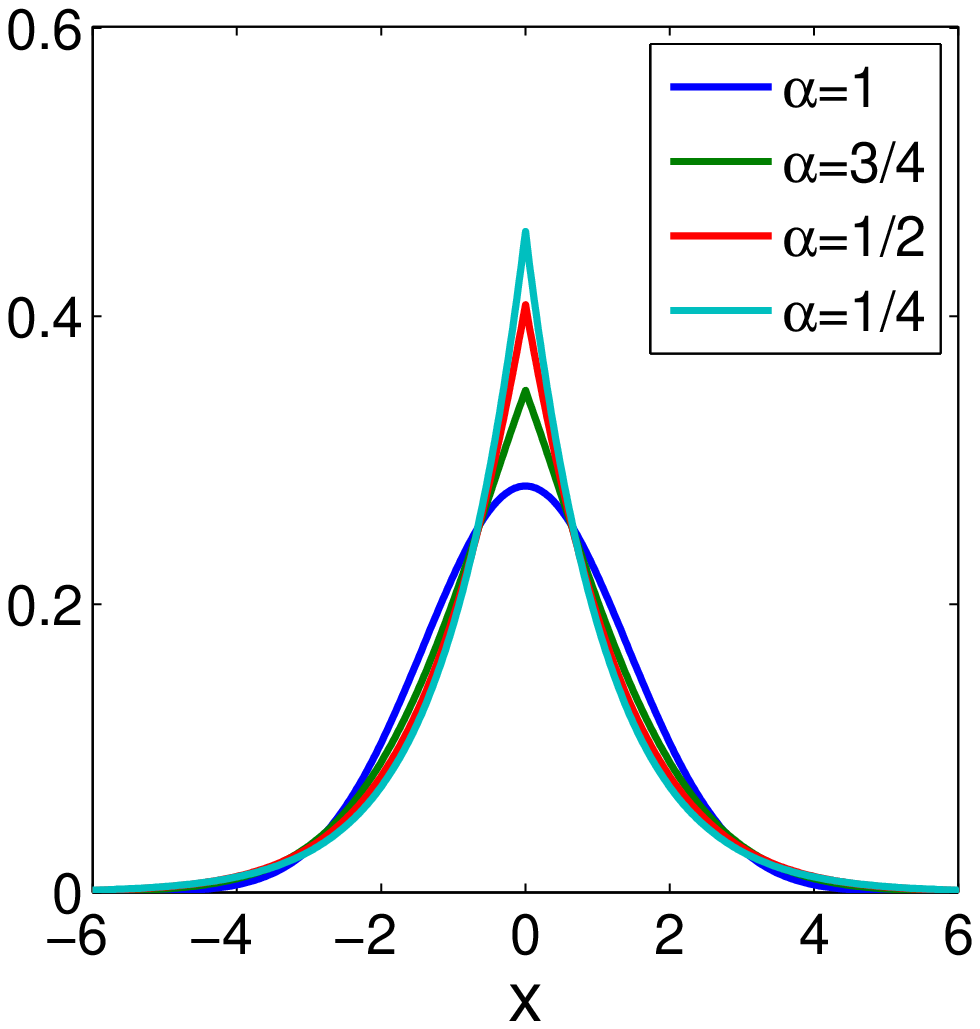}\\
   (a) $t=0.1$ & (b) $t=1$
  \end{tabular}
  \caption{The profile of the fundamental solution $K_\alpha(x,t)$ at (a) $t=0.1$ and (b)
  $t=1$.}\label{fig:wright}
\end{figure}

\subsection{Singular value decomposition} We shall follow the well-established practice in
the inverse problem community, i.e., using the singular value decomposition, as the main tool for numerically analyzing
the problem behavior \cite{Hansen:1998}. Specifically, we shall numerically compute the forward map $F$,
and analyze its behavior to gain insights into the inverse problem. Given a matrix $\mathbf{A}\in
\mathbb{R}^{n\times m}$, its singular value decomposition is given by
\begin{equation*}
  \mathbf{A} = \mathbf{U}\Sigma\mathbf{V}^\mathrm{t},
\end{equation*}
where $\mathbf{U}=[\mathbf{u}_1 \ \ldots \ \mathbf{u}_n]\in\mathbb{R}^{n\times n}$ and $\mathbf{V}=[\mathbf{v}_1\
\ldots \ \mathbf{v}_m]\in\mathbb{R}^{m\times m}$ are column orthonormal matrices, and $\Sigma\in\mathbb{R}^{n\times
m}=\mathrm{diag}(\sigma_1,\ldots,\sigma_r,0,\ldots,0)$ is a diagonal matrix, with the diagonal elements $\{\sigma_i
\}$ being nonnegative and listed in a descending order $\sigma_1>\ldots>\sigma_r>0$, and $r$ being the (numerical)
rank of the matrix $\mathbf{A}$. The diagonal element $\sigma_i$ is known as the $i$th singular value, and the
corresponding columns of $\mathbf{U}$ and $\mathbf{V}$, i.e., $\mathbf{u}_i$ and $\mathbf{v}_i$, are called the left
and right singular vectors, respectively.

One simple measure of the conditioning of a linear inverse problem $\mathbf{Ax}=\mathbf{b}$ is the condition number
$\mathrm{cond}(\mathbf{A})$, which is defined as the ratio of the largest to the smallest nonzero singular value, i.e.,
\begin{equation*}
  \mathrm{cond}(\mathbf{A}) = \sigma_1/\sigma_r.
\end{equation*}
In particular, if the condition number is small, then the data error will not be amplified much. In the
case of a large condition number, the issue is more delicate: it may or may not amplify the data perturbation
greatly. A more complete picture is provided by the singular value spectrum  $(\sigma_1,\sigma_2,\ldots,
\sigma_r)$. Especially, a singular value spectrum gradually decaying to zero without a clear gap is characteristic
of many discrete ill-posed problems, which is reminiscent of the spectral behavior of compact operators. We
shall adopt these simple tools to analyze related inverse problems below.

In addition, using singular value decomposition and regularization techniques, e.g. Tikhonov regularization or
truncated singular value decomposition, one can conveniently obtain numerical reconstructions, even though this
might not be the most efficient way to do so. However, we shall not delve into the extremely important question
of practical reconstructions, since it relies heavily on \textit{a priori} knowledge on the sought for solution
and the statistical nature (Gaussian, Poisson, Laplace ...) of the contaminating noise in the data, which
will depend very much on the specific application. We refer interested readers to the monographs \cite{EnglHankeNeubauer:1996,
SchusterKaltenbacherHofmannKazimierski:2012,ItoJin:2014} and the survey \cite{JinMaass:2012} for updated accounts on
regularization methods for constructing stable reconstructing procedures and efficient computational techniques.
We will also briefly mention below related works on the application of regularization techniques to inverse problems for FDEs.

\section{Inverse problems for time fractional diffusion}\label{sec:time}
In this section, we consider several model inverse problems for the following one-dimensional time
fractional diffusion equation on the unit interval $\Omega=(0,1)$:
\begin{equation}\label{eqn:fde-time}
  \partial_t^\alpha u - u_{xx} + qu = f \quad \mbox{ in } \Omega\times(0,T],
\end{equation}
with the fractional order $\alpha\in(0,1)$, the initial condition $u(0)=v$ and suitable boundary conditions. Although
we consider only the one-dimensional model, the analysis and computation in this part can be extended into the
general multi-dimensional case, upon suitable modifications. Recall that
 $\partial_t^\alpha u$ denotes the Djrbashian-Caputo fractional derivative of order $\alpha$ with respect to
time $t$.  For $\alpha =1$, the fractional derivative $ \partial_t^\alpha u$ coincides with the usual first-order
derivative $u'$, and accordingly, the model \eqref{eqn:fde-time} reduces to the classical diffusion equation. Hence it
is natural to compare inverse problems for the model \eqref{eqn:fde-time} with that for the standard diffusion equation.
We shall discuss the following four inverse problems, i.e., the backward problem, sideways problem, inverse source
problem and inverse potential problem. In the first three cases, we shall assume a zero potential $q=0$.
We will aslo discuss the solution of an inverse coefficient problem using fractional calculus.

\subsection{Backward fractional diffusion}
First we consider the time fractional backward diffusion. By the linearity of the inverse problem, we may assume that
equation \eqref{eqn:fde-time} is prescribed with a homogeneous Dirichlet boundary condition, i.e., $u=0$ at
$x=0,1$, and the initial condition $u(0)=v$. Then the inverse problem reads: given the final time data $g=u(T)$,
find the initial condition $v$. It arises in, for example, the determination of a stationary contaminant source in
underground fluid flow.

To gain insight, we apply the separation of variables. Let $\{(\lambda_j,\phi_j)\}$, with $\lambda_j
=(j\pi)^2$ and $\phi_j=\sqrt{2}\sin j\pi x$, be the Dirichlet eigenpairs of the negative Laplacian on the
interval $\Omega$. The eigenfunctions $\{\phi_j\}$ form an orthonormal basis of the $L^2(\Omega)$ space. Then
using the Mittag-Leffler function $E_{\alpha,\beta}(z)$ defined in \eqref{eqn:mitlef}, the solution $u$ to equation
\eqref{eqn:fde-time} can be expressed as
\begin{equation*}
  u(x,t) = \sum_{j=1}^\infty E_{\alpha,1}(-\lambda_jt^\alpha)(v,\phi_j)\phi_j(x).
\end{equation*}
Therefore, the final time data $g=u(T)$ is given by
\begin{equation*}
  g(x) = \sum_{j=1}^\infty E_{\alpha,1}(-\lambda_j T^\alpha)(v,\phi_j)\phi_j(x).
\end{equation*}
It follows directly that the initial data $v$ is formally given by
\begin{equation*}
  v = \sum_{j=1}^\infty \frac{(g,\phi_j)}{E_{\alpha,1}(-\lambda_j T^\alpha)}\phi_j.
\end{equation*}
Since the function $E_{\alpha,1}(-t)$ is completely monotone on the positive real axis $\mathbb{R}^+$ \cite{Pollard:1948} for any
$\alpha\in (0,1]$, the denominator in the representation does not vanish.
In case of $\alpha=1$, the formula reduces to the familiar expression
\begin{equation*}
  v = \sum_{j=1}^\infty e^{\lambda_j T}(g,\phi_j)\phi_j.
\end{equation*}
This formula shows clearly the well-known severely ill-posed nature of the backward diffusion problem:
the perturbation in the $j$th Fourier mode $(g,\phi_j)$ of the (noisy) data $g$ is amplified by an
exponentially growing factor $e^{\lambda_jT}$, which can be astronomically large, even for a very
small index $j$, if the terminal time $T$ is not very small. Hence it is always severely ill-conditioned
and we must multiply the $j$th Fourier mode of the data $g$ by a factor $e^{\lambda_j T}$ in order to recover
the corresponding mode of the initial data $v$

In the fractional case, by Lemma \ref{lem:mitlef}, the Mittag-Leffler function $E_{\alpha,1}(z)$
decays only linearly on the negative real axis $\mathbb{R}^-$, and thus the multiplier $1/E_{\alpha,1}(-\lambda_jT^\alpha)$
grows only linearly in $\lambda_j$, i.e., $1/E_{\alpha,1}(-\lambda_jT^\alpha)\sim \lambda_j$,
which is very mild compared to the exponential growth $e^{\lambda_jT}$ for the case $\alpha=1$,
and thus the fractional case is only mildly ill-posed. Roughly, the $j$th Fourier mode of the initial
data $v$ now equals the $j$th mode of the data $g$ multiplied by $\lambda_j$. More
precisely, it amounts to the loss of two spatial derivatives \cite[Theorem 4.1]{SakamotoYamamoto:2011}
\begin{equation}\label{eqn:backwardiff}
  \|v\|_{L^2(\Omega)}\leq c\|u(T)\|_{H^2(\Omega)}.
\end{equation}
Intuitively, the history mechanism of the anomalous diffusion process retains the complete dynamics
of the physical process, including the initial data, and thus it is much easier to go backwards to
the initial state $v$. This is in sharp contrast to the classical diffusion, which has only
a short memory and loses track of the preceding states quickly. This result has become quite
well-known in the inverse problems community and has contributed to a belief that ``inverse problems for FDEs are
less ill-conditioned than their classical counterparts'' -- throughout this paper we will see that
this conclusion as a general statement can be quite far from the truth.

Does this mean that for all terminal time $T$ the fractional case is always less ill-posed than the classical one?
The answer is yes, in the sense of the norm on the data space in which the data $g$ lies. Does this mean that from a
computational stability standpoint that one can always solve the backward fractional problem more effectively than
for the classical case? The answer is no, and the difference can be substantial. To illustrate the point, let $J$ be
the highest frequency mode required of the initial data $v$ and assume that we believe we are able to multiply the
first $J$ modes $g_j=(g,\phi_j)$, $j=1,2,\ldots,J$, by a factor no larger than $M$ (which roughly assumes that the
noise levels in both cases are comparable). By the monotonicity of the function $E_{\alpha,1}(-t)$ in $t$, it
suffices to examine the $J$th mode. For the heat equation $v_J:= (v,\phi_J) = e^{\lambda_J T}
g_J$ and provided that $T = T_J < \lambda_J/\log{M}$ this is feasible. For a fixed $J$, if $T_\alpha^\star$ denotes
the point where
\begin{equation*}
  e^{-\lambda_J T_\alpha^\star} = E_{\alpha,1}(-\lambda_J T_\alpha^\star),
\end{equation*}
then in the fractional case for $T<T_\alpha^\star$ the growth factor on
$g_J$ will exceed $M$ for any $T<T_\alpha^\star$. In Table \ref{tab:backwardfrac} below, we present the critical value
$T_\alpha^\star$ for several values of the fractional order $\alpha$ and the maximum number of modes $J$. The numbers
in the table are very telling. For example, for
the case $J=5$, $\alpha=1/4$ and  $T=0.02$ (which is approximately one half the value of $T_\alpha^\star$), the growth
factor is about  $1.6$ for the heat equation but about $113$ for the fractional case. With $J=10$ and $\alpha=1/4$ and
$T = T_\alpha^\star$ the growth factor is around $336$. If $T = T_\alpha^\star/10$ then it has again dropped to less
than $2$ for the heat equation but about $190$ for the fractional case. Of course, for $T> T_\alpha^\star$ the situation
completely reverses. With $J=10$, $\alpha = 1/4$ and $T=10\,T_\alpha^\star$ the growth factor is a possibly workable value
of around $600$; while for the heat equation it is greater than $10^{25}$. We reiterate that the apparent contradiction
between the theoretical ill-conditioning and numerical stability is due to the spectral cutoff present in any practical
reconstruction procedure.

\begin{table}[hbt!]
  \centering
  \caption{The critical values $T_\alpha^*$ for fractional backward diffusion.}\label{tab:backwardfrac}
  \begin{tabular}{cccc}
    \hline
    $\alpha\backslash J$ & 3  & 5 & 10\\
    \hline
      1/4  &  0.0442  &  0.0197  & 0.0059\\
      1/2  &  0.0387  &  0.0163  & 0.0049\\
      3/4  &  0.0351  &  0.0142  & 0.0040\\
    \hline
  \end{tabular}
\end{table}

Next we examine the influence of the fractional order $\alpha$ on the inversion step more closely. To this end, we
expand the initial condition $v$ in the piecewise linear finite element basis functions defined on a uniform partition of
the domain $\Omega=(0,1)$ with 100 grid points. Then we compute the discrete forward map $F$ from the initial condition
to the final time data $g=u(T)$, defined on the same mesh. Numerically, this can be achieved by a fully discrete scheme based on
the L1 approximation in time and the finite difference method in space; see Appendix \ref{app:timefrac} for a description of
the numerical method. The ill-posed behavior of the discrete inverse problem is then analyzed using singular value
decomposition. A similar experimental setup will be adopted for other examples below.

\begin{figure}[hbt!]
  \centering
  \begin{tabular}{cc}
     \includegraphics[trim = 1cm .1cm 1cm 0cm, clip = true, width=7.5cm]{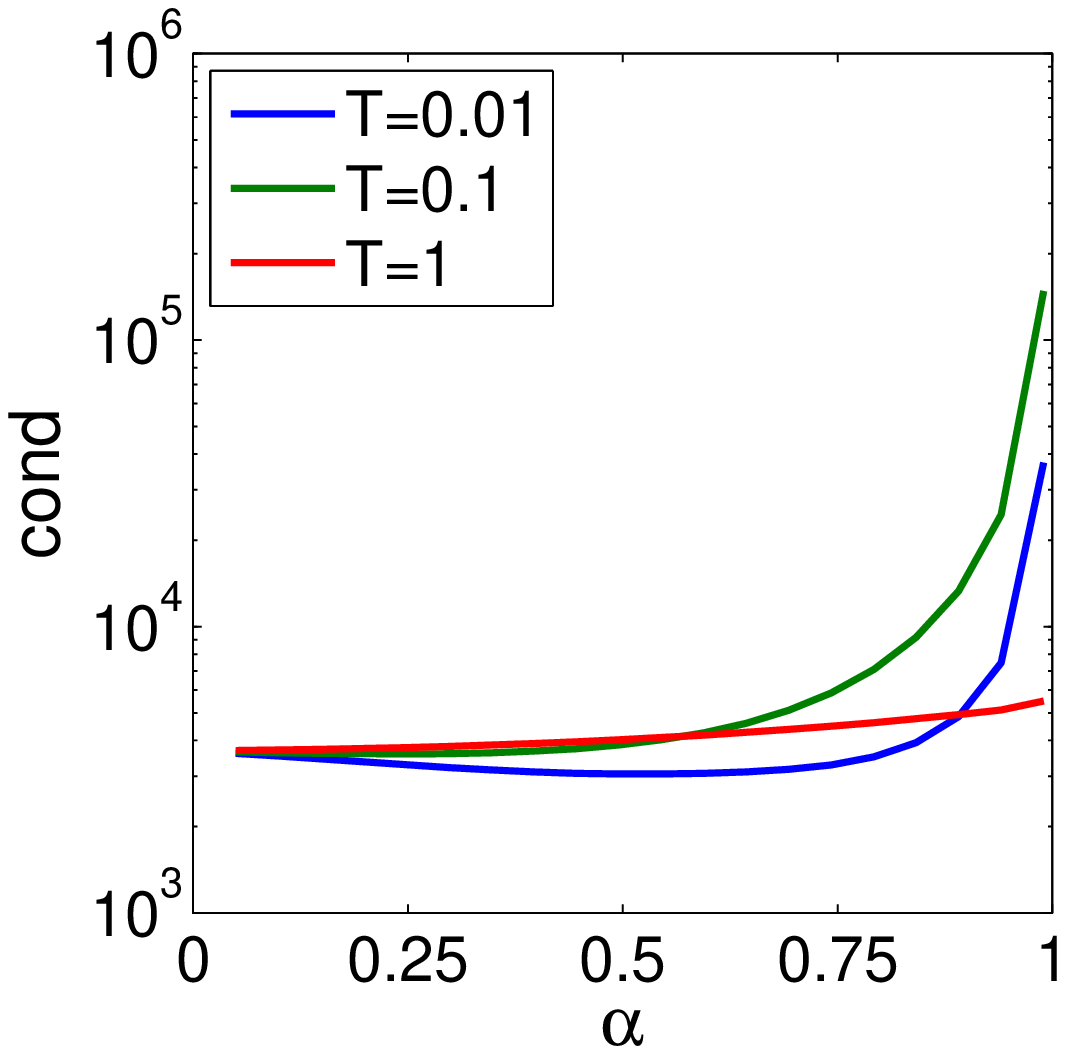} &
     \includegraphics[trim = 1cm .1cm 1cm 0cm, clip = true, width=7.5cm]{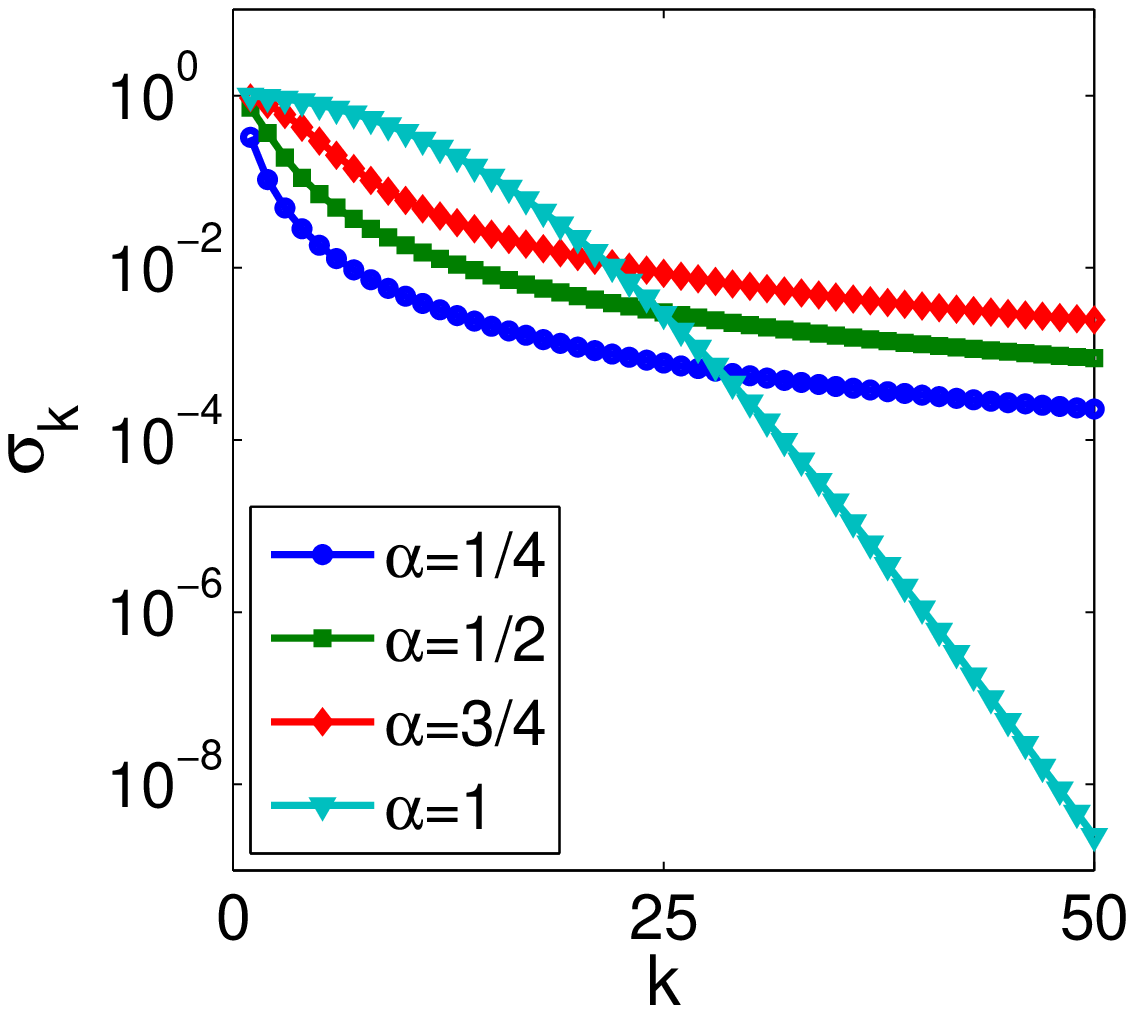}\\
     (a) condition number & (b) singular value spectrum
  \end{tabular}
  \caption{(a) The condition number v.s. the fractional order $\alpha$, and (b) the singular value spectrum at $T=0.01$ for the backward
  fractional diffusion. We only display the first 50 singular values.}\label{fig:backward:cond}
\end{figure}

The numerical results are shown in Fig. \ref{fig:backward:cond}.
The condition number of the (discrete) forward map $F$ stays mostly around $O(10^4)$ for a fairly broad
range of $\alpha$ values, which holds for all three different terminal times $T$. This can be attributed to
the fact that for any $\alpha\in(0,1)$, backward fractional diffusion amounts to a two spacial derivative
loss, cf. \eqref{eqn:backwardiff}. Unsurprisingly, as the fractional order $\alpha$ approaches unity, the
condition number eventually blows up, recovering the severely ill-posed nature of the classical backward
diffusion problem, cf. Fig. \ref{fig:backward:cond}(a). Further, we observe that the smaller is the terminal time $T$,
the quicker is the blowup. The precise mechanism for this observation remains unclear. Interestingly, the condition number
is not monotone with respect to the fractional order $\alpha$, for a fixed $T$. This might imply potential
nonuniqueness in the simultaneous recovery of the fractional order $\alpha$ and the initial data $v$.
The singular value spectra at $T=0.01$ are shown in Fig. \ref{fig:backward:cond}(b). Even though the
condition numbers for $\alpha=1/4$ and $\alpha=1/2$ are quite close, their singular value spectra actually
differ by a multiplicative constant, but their decay rates are almost identical, thereby showing comparable
condition numbers. This shift in singular value spectra can be explained
by the local decay behavior of the Mittag-Leffler function, cf. Fig. \ref{fig:mitlef}(a): the smaller is
the fractional order $\alpha$, the faster is the decay around $t=0$.

\begin{figure}[hbt!]
  \centering
  \begin{tabular}{cc}
     \includegraphics[trim = .5cm .1cm 1cm 0cm, clip = true, width=7.5cm]{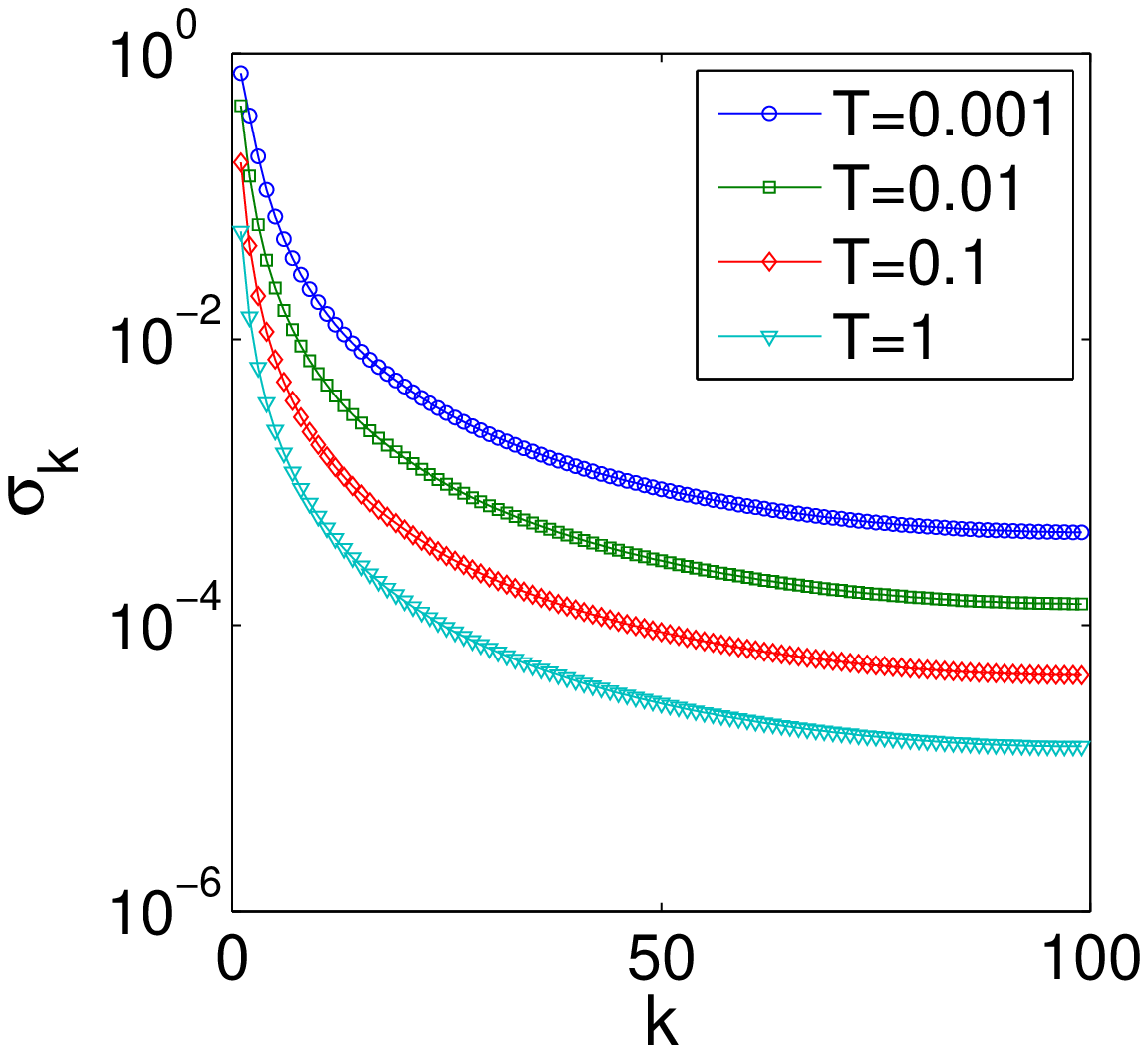} &
     \includegraphics[trim = .5cm .1cm 1cm 0cm, clip = true, width=7.5cm]{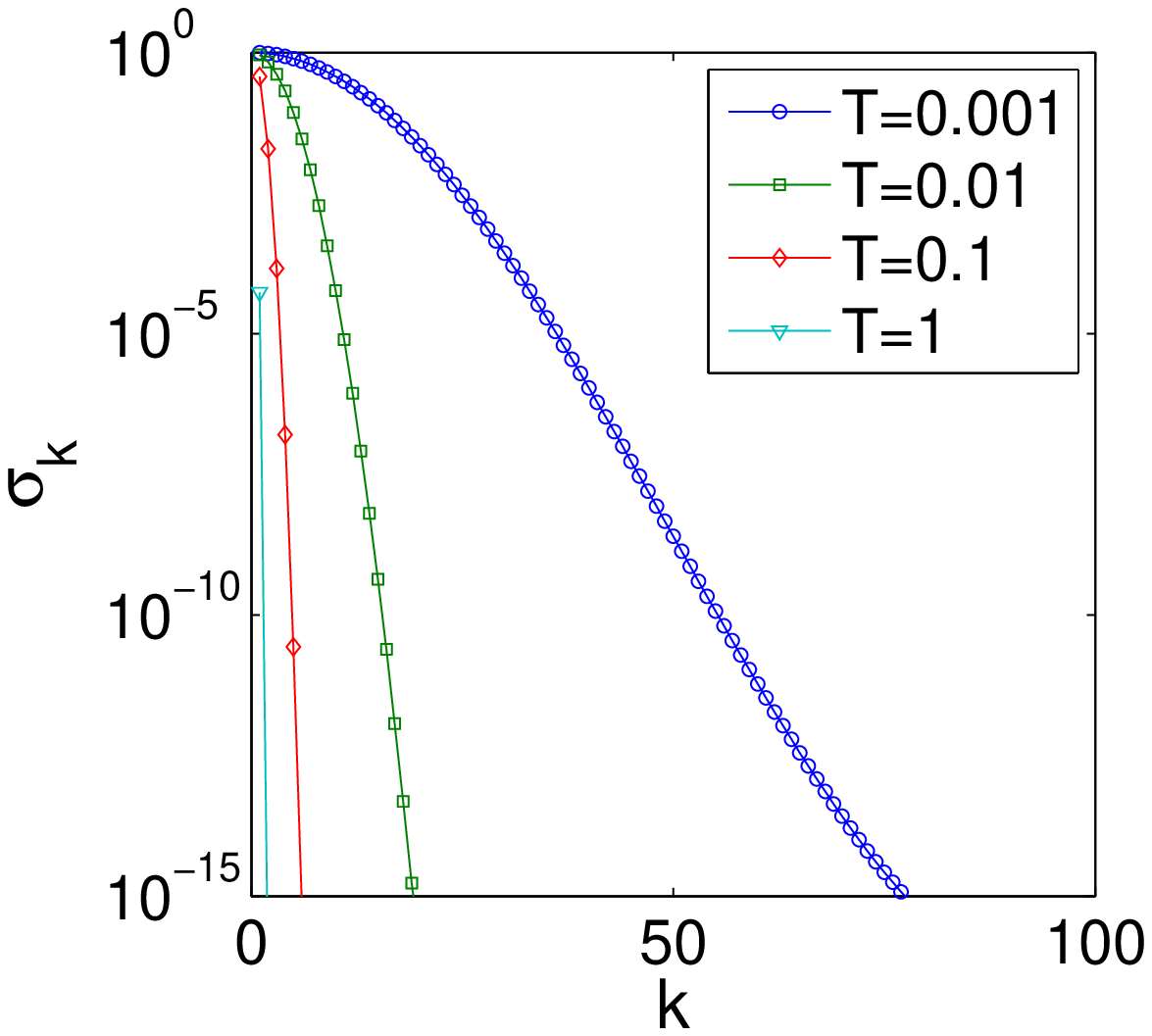}\\
     (a) $\alpha=1/2$ & (b) $\alpha=1$
  \end{tabular}
  \caption{The singular value spectrum of the forward map $F$ from the initial data to the final time data, for (a) $\alpha=1/2$ and (b) $\alpha=1$, at
  four different times for the time fractional backward diffusion.}\label{fig:backward:difftime}
\end{figure}

Even though the condition number is very informative about the (discretized) problem, it does not provide
a full picture, especially when the condition number is large, and the singular value spectrum can be far
more revealing. The spectra for two different $\alpha$ values are given in Fig.
\ref{fig:backward:difftime}. At $\alpha=1/2$, the singular values decay at almost the same algebraic rate,
irrespective of the terminal time $T$. This is expected from the two-derivative loss for any $\alpha<1$.
However, for $\alpha=1$, the singular values decay exponentially, and the decay rate increases dramatically
with the increase of the terminal time $T$. For $T=0.001$, there are a handful of ``significant'' singular
values, say above $10^{-3}$, but when the time $T$ increases to $T=1$, there is only one meaningful
singular value remaining. The distribution of the singular values has important practical consequences.
For a small time $T$, the first few singular values for the classical diffusion
case actually might be much larger than that for the fractional case, which indicates that the classical
case is actually numerically much easier to recover in this regime, concurring with the observations drawn from
 Table \ref{tab:backwardfrac}. For example, at $T=0.001$, the first twenty singular values
are larger than the fractional counterpart,  cf. Fig. \ref{fig:backward:cond}(b), and hence, the
first twenty modes, i.e., left singular vectors, are more stable in the reconstruction procedure.

The mathematical model \eqref{eqn:fde-time} is rescaled with a unit diffusion coefficient. In practice,
there is always a diffusion coefficient $\sigma$ in the elliptic operator, i.e.,
\begin{equation*}
  \partial_t^\alpha u - \nabla\cdot(\sigma\nabla u) + qu = f.
\end{equation*}
For example, the thermal conductivity $\sigma$ of the gun steel at moderate temperature is about $1.8\times10^{-5}
\ \mathrm{m}^2/\mathrm{s}$ \cite{Carasso:1982} and the diffusion coefficient of the oxygen in water at 25
Celsius is	$2.10\times10^{-5}\ \mathrm{cm}^2/\mathrm{s}$ \cite{Cussler:1997}. Mathematically, this does
not change the ill-posed nature of the inverse problem. However, the presence of a diffusion coefficient $\sigma$
has important consequence: it enables the practical feasibility of the classical backward diffusion problem
(and likely for many other inverse problems for the diffusion equation). Physically, a constant conductivity
$\sigma$ amounts to rescaling the final time $T$ by $T'=T/\sigma$. In the fractional case, a similar but nonlinear
scaling law $T'=T/\sigma^{1/\alpha}$ remains valid.

Numerically, time fractional backward diffusion has been extensively studied. Liu and Yamamoto \cite{LiuYamamoto:2010}
proposed a numerical scheme for the one-dimensional fractional backward problem based on the quasi-reversibility method \cite{LattesLions:1969},
and derived error estimates for the approximation, under a priori smoothness assumption on the initial condition.
This represents one of the first works on inverse problems in anomalous diffusion. Later, Wang and Liu
\cite{WangLiu:2013} studied total variation regularization for two-dimensional fractional backward diffusion, and analyzed its well-posedness
of the optimization problem and the convergence of an iterative scheme of Bregman type. Wei and Wang
\cite{WeiWang:2014a} developed a modified quasi-boundary value method for the problem in a general domain, and established error estimates
for both a priori and a posteriori parameter choice rules. In view of better stability results in the fractional case, one naturally
expects better error estimates than the classical diffusion equation, which is confirmed by these studies.

\subsection{Sideways fractional diffusion}
Next we consider the sideways problem for time fractional diffusion. There are several possible formulations,
e.g., the quarter plane and the finite space domain. The quarter plane sideways fractional diffusion problem is as
follows. Let $u(x,t)$ be defined in $(0,\infty)\times(0,\infty)$ by
\begin{equation*}
   \partial^\alpha_t u -  u_{xx} = 0,\qquad x>0,\ t>0,
\end{equation*}
and the boundary and initial conditions
\begin{equation*}
  u(x,0) = 0 \quad \mbox{and}\quad u(0,t) = f(t),
\end{equation*}
where we assume $|u(x,t)| \leq c_1 e^{c_2x^2}$. We do not know the left boundary condition $f$, but are able to
measure $u$ at an intermediate point $x=L>0$, $h(t)=u(L,t)$. The inverse problem is: given the (noisy) data $h$,
find the boundary condition $f$. The solution $u$ of the forward problem is given by a convolution integral
with the kernel being the spatial derivative $K_{\alpha,x}(x,s)$ of the fundamental solution $K_\alpha(x,s)$,
cf. \eqref{eqn:fundsol}, by
\begin{equation*}
  u(x,t) = \int_0^t K_{\alpha,x}(x,t-s)f(s)ds.
\end{equation*}
This representation is well known for the case $\alpha=1$, and it was first derived by Carasso \cite{Carasso:1982};
see also \cite{Cannon:1984} for related discussions. It leads to a convolution integral equation for the unknown
$f$ in terms of the given data $h$
\begin{equation*}
  h(t) = \int_0^t R_\alpha(t-s)f(s) ds,
\end{equation*}
where the convolution kernel $R_\alpha(s)$ is given by a Wright function in the form
\begin{equation*}
  \begin{aligned}
    R_\alpha(s) & = \frac{1}{2s^\alpha}W_{-\frac{\alpha}{2},2-\frac{\alpha}{2}}(-L s^{-\alpha/2})\\
       & = \sum_{k=0}^\infty \frac{(-L)^k}{k!\, \Gamma(-\frac{\alpha}{2}k + 2 -\frac{\alpha}{2})} \,s^{-k\frac{\alpha}{2}-\alpha}.
  \end{aligned}
\end{equation*}
In case of $\alpha=1$, i.e., classical diffusion, the kernel $R(s)$ is given explicitly by
\begin{equation*}
  R(s)= \frac{L}{2\sqrt{\pi}}s^{-\frac{3}{2}} e^{-\frac{L^2}{4s}} \in C^\infty(0,\infty).
\end{equation*}
Since all its derivatives vanish at $s=0$, the classical theory of Volterra integral equations of the first kind
\cite{Lamm:2000} implies the extreme ill-conditioning of the problem. This is not surprising: We are, after all,
mapping a function $f\in C^0(0,\infty)$ to an element in $C^\infty(0,\infty)$. The conditioning of the time fractional
sideways problem again depends on the convolution kernel $R_\alpha$ and its
derivatives at $s=0$ and in this case is the value of the Wright function $W_{-\alpha/2,2-\alpha/2}(-z)$
and its derivatives as $z\to\infty$. These are again zero (in fact the Wright function $W_{-\alpha/2,
2-\alpha/2}(z)$ also decays exponentially to zero for large negative arguments, cf. Lemma \ref{lem:wright}),
and thus the fractional sideways problem is also severely ill-posed. However, this analysis does not
show their difference in the degree of ill-posedness: even though both are severely ill-posed, their
practical computational behavior can still be quite different, as we shall see below.

To see their difference in the degree of ill-posedness, we examine another variant of the sideways problem on a finite
interval $\Omega=(0,1)$, with Cauchy data prescribed on the axis $x=0$, i.e. given zero initial condition $u_0=u(x,0)=0$,
recovering $h=u(1,t)$ from the lateral Cauchy data at $x=0$:
\begin{equation*}
  u(0,t) = f(t), \ \  u_x(0,t) = g(t),\ \ t\geq 0
\end{equation*}
This problem is also known as the lateral Cauchy problem in the literature.
In the case $\alpha=1$, it is known that the inverse problem is severely ill-posed \cite{Cannon:1984,HaoReinhardt:1997}.
To gain insight into the fractional case, we apply the Laplace transform. With $\ \widehat{}\ $ being the Laplace
transform in time, and noting the Laplace transform of the Caputo derivative $\mathcal{L}(\partial_t^\alpha u)
= z^\alpha \widehat{u}(z)-z^{\alpha-1}u_0$ \cite[Lemma 2.24]{KilbasSrivastavaTrujillo:2006}, we deduce
\begin{equation*}
  z^\alpha \widehat u(x,z) - \widehat u_{xx}(x,z) = 0, \ \ \widehat u (0) = \widehat f, \ \ \widehat u_x(0) = \widehat g.
\end{equation*}
The general solution $\widehat u$ is given by $\widehat u(x,z) = \widehat f \cosh z^{\alpha/2}x + \widehat g
\frac{\sinh z^{\alpha/2}x}{z^{\alpha/2}x}$ and thus the solution $\widehat h(z)=\widehat u(1,z)$ at $x=1$ is given by
\begin{equation*}
\widehat h=\widehat f \cosh z^{\alpha/2} + \widehat g \frac{\sinh z^{\alpha/2}}{z^{\alpha/2}}.
\end{equation*}
The solution $h(t)$ can then be recovered by an inverse Laplace transform
\begin{equation}\label{eqn:exp-hhat}
  h(t) = \frac{1}{2\pi \mathrm{i}} \int_{Br} \widehat h e^{zt}d z,
\end{equation}
where $Br=\{z\in \mathbb{C}: \Re z=\sigma, \sigma>0\}$ is the Bromwich path. Upon deforming the contour suitably, this
formula will allow developing
an efficient numerical scheme for the sideways problem via quadrature rules \cite{WeidemanTrefethen:2007},
provided that the lateral Cauchy data is available for all $t>0$. The expression \eqref{eqn:exp-hhat} indicates
that, in the fractional case, the sideways problem still suffers from severe ill-posedness in theory, since
the high frequency modes of the data perturbation are amplified by an exponentially growing multiplier $e^{z^{\alpha/2}}$.
However, numerically, the degree of ill-posedness decreases dramatically as the fractional order $\alpha$ decreases
from unity to zero, since as $\alpha\to0^+$, the multipliers are growing at a much slower rate, and thus we have a better
chance of recovering many more modes of the boundary data. In other words, both the classical and fractional sideways problems are
severely ill-posed in the sense of error estimates between the norms in the data and unknowns; but with a fixed
frequency range, the behavior of the time fractional sideways problem can be much less ill-posed. Hence, anomalous diffusion
mechanism does help substantially since much more effective reconstructions are possible  in the fractional case.

\begin{figure}[hbt!]
  \centering
  \begin{tabular}{cc}
     \includegraphics[trim = .5cm .1cm 1cm 0cm, clip = true, width=7.5cm]{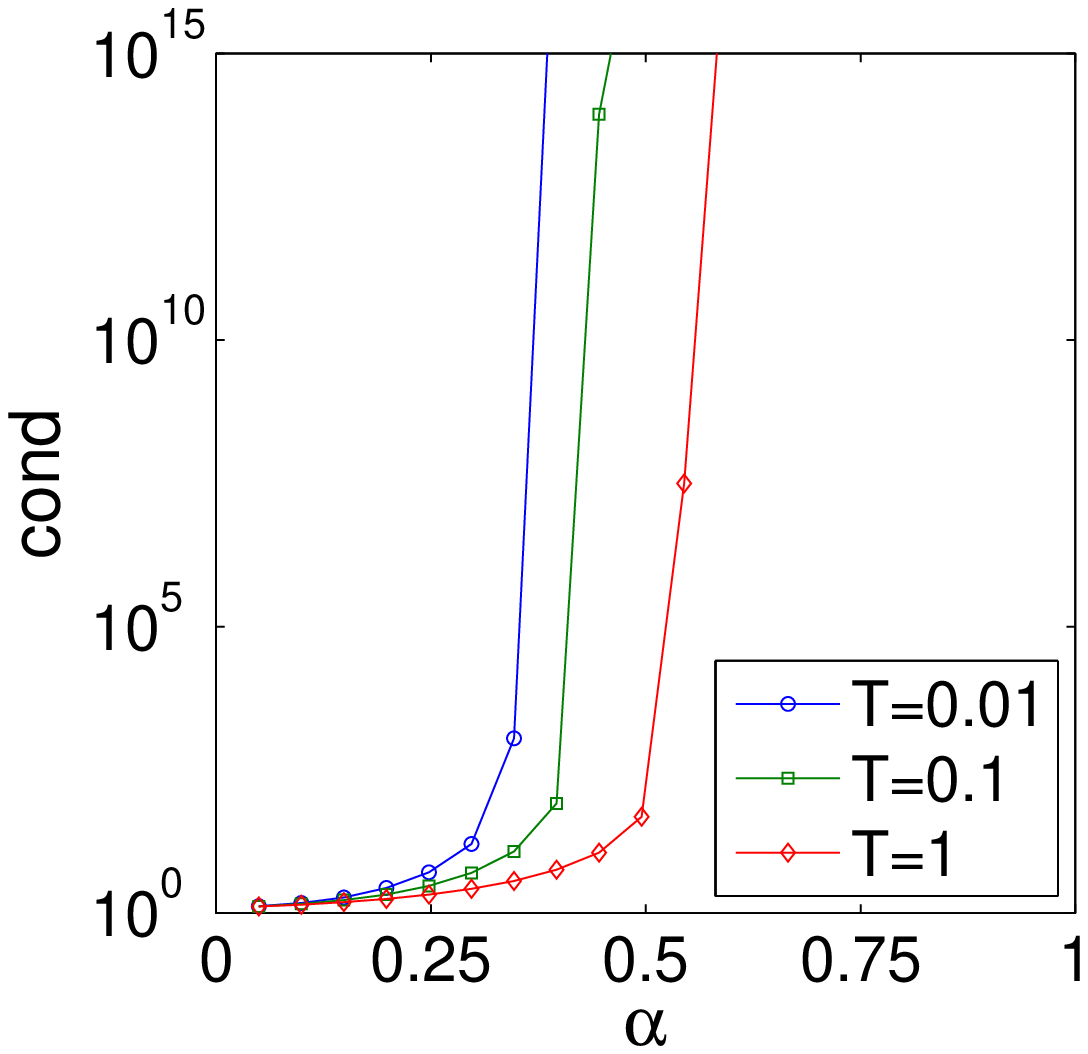} &
     \includegraphics[trim = .5cm .1cm 1cm 0cm, clip = true, width=7.5cm]{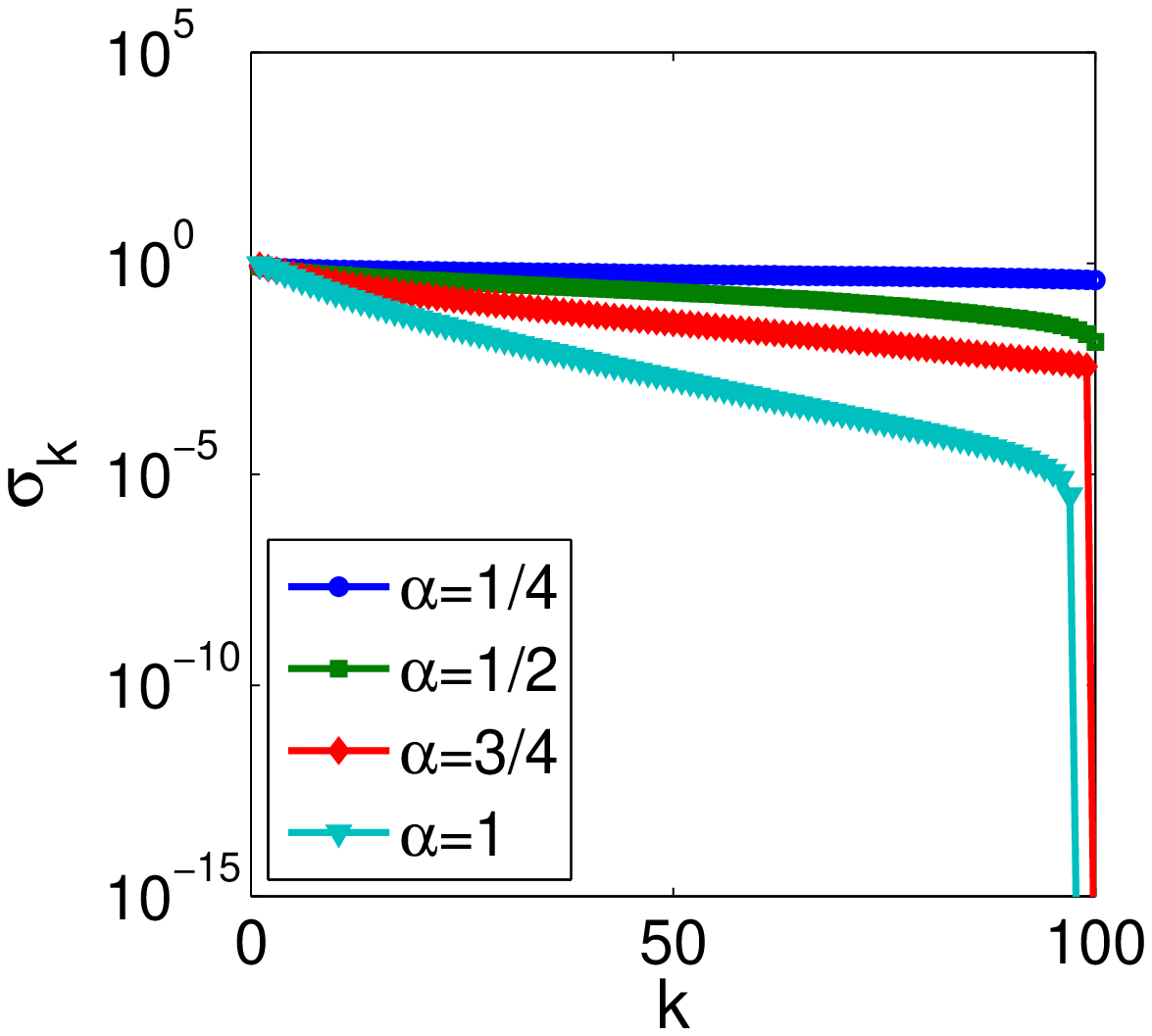}\\
     (a) condition number & (b) singular value spectrum
  \end{tabular}
  \caption{(a) The condition number and (b) singular value spectrum at $T=1$ for the time fractional
  sideways problem.}\label{fig:sideways:cond}
\end{figure}

Next we illustrate the point numerically. The numerical results for the sideways problem are given in Fig.
\ref{fig:sideways:cond}. It is observed that the degree of ill-posedness of the finite-dimensional discretized
version of the inverse problem indeed decreases dramatically with the decrease of the fractional order $\alpha$,
cf. Fig. \ref{fig:sideways:cond}(a), which agrees well with the preceding discussions. Surprisingly, for $T=1$
there is a sudden transition around $\alpha=1/2$, below which the sideways problem behaves as if nearly well-posed,
but above which the conditioning deteriorates dramatically with the increase of the fractional order $\alpha$ and
eventually it recovers the properties of the classical sideways problem. Similar transitions are observed for other
terminal times. This might be related to the discrete setting, for which there is an inherent frequency cutoff.
Further, as the fractional order $\alpha$ approaches zero, the problem reaches a quasi-steady state much quicker
and thus the forward map $F$ can have only fairly localized elements along the main diagonal. To give a more complete
picture, we examine the singular value spectrum in Fig. \ref{fig:sideways:cond}(b). Unlike the backward diffusion
problem discussed earlier, the singular values are actually decaying only algebraically, even for $\alpha=1$, and then
there might be a few tiny singular values contributing to the large condition number. The larger is the fractional
order $\alpha$, the more tiny singular values are in the spectrum. Hence, in the discrete setting, even for
$\alpha=3/4$, the problem is still nearly well-posed, despite the large apparent condition number, since a few tiny
singular values with a distinct gap from the rest of the spectrum are harmless in most regularization techniques.

\begin{figure}[hbt!]
  \centering
  \begin{tabular}{cc}
     \includegraphics[trim = 0cm 0cm 0cm 0cm, clip = true, width=7cm]{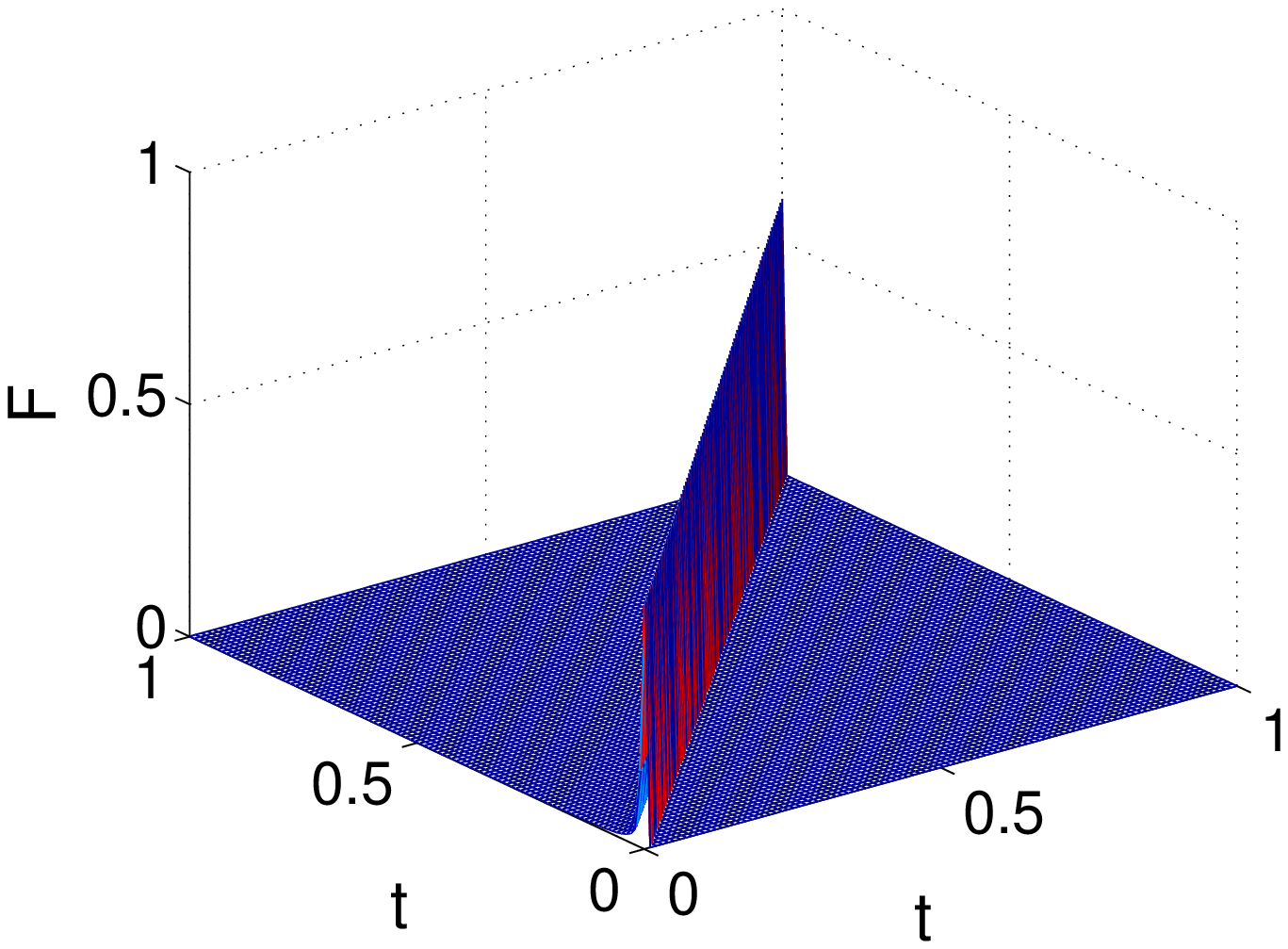} &
     \includegraphics[trim = 0cm 0cm 0cm 0cm, clip = true, width=7cm]{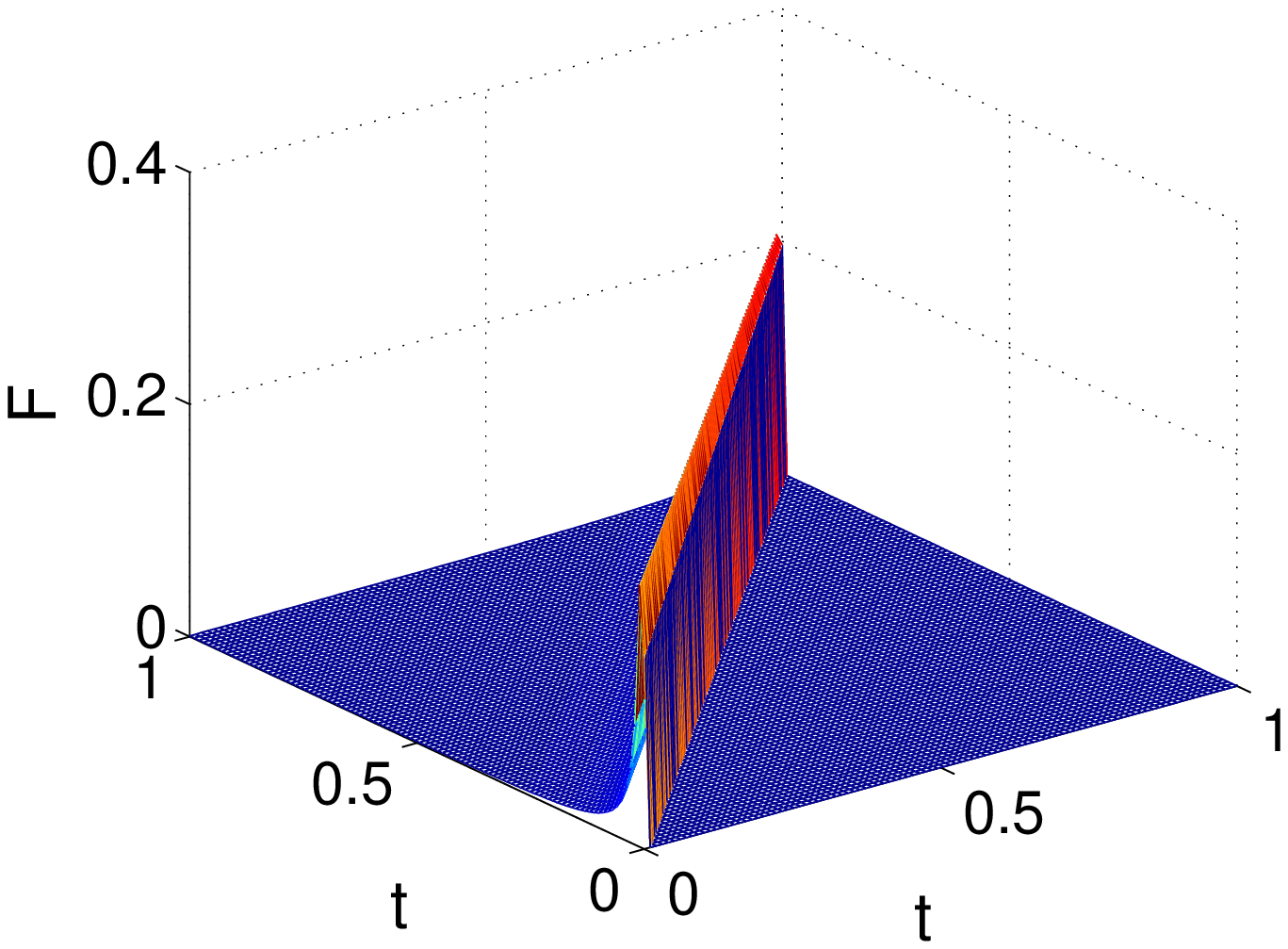}\\
     (a) $\alpha=1/4$ & (b) $\alpha=1/2$\\
     \includegraphics[trim = 0cm 0cm 0cm 0cm, clip = true, width=7cm]{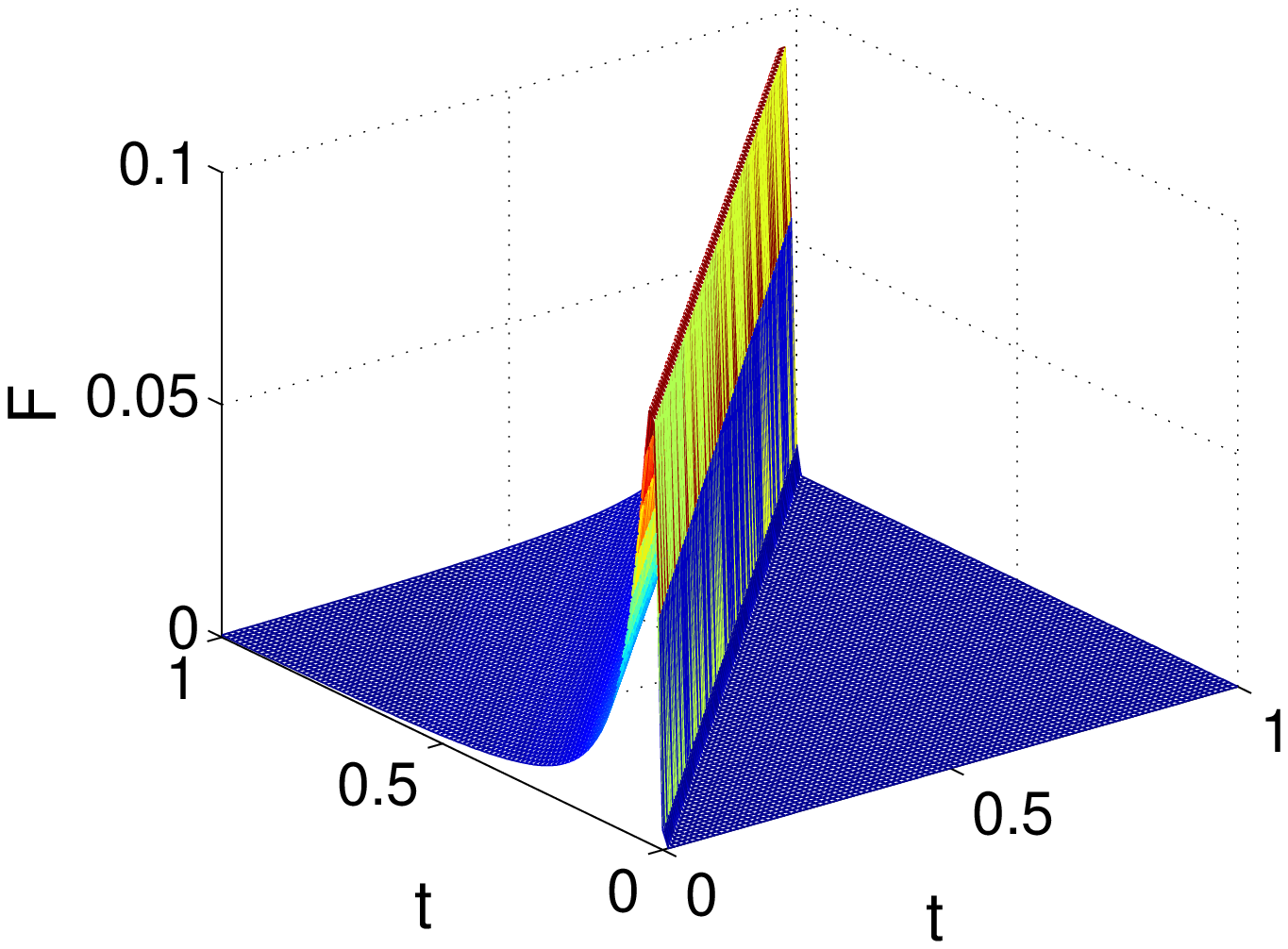} &
     \includegraphics[trim = 0cm 0cm 0cm 0cm, clip = true, width=7cm]{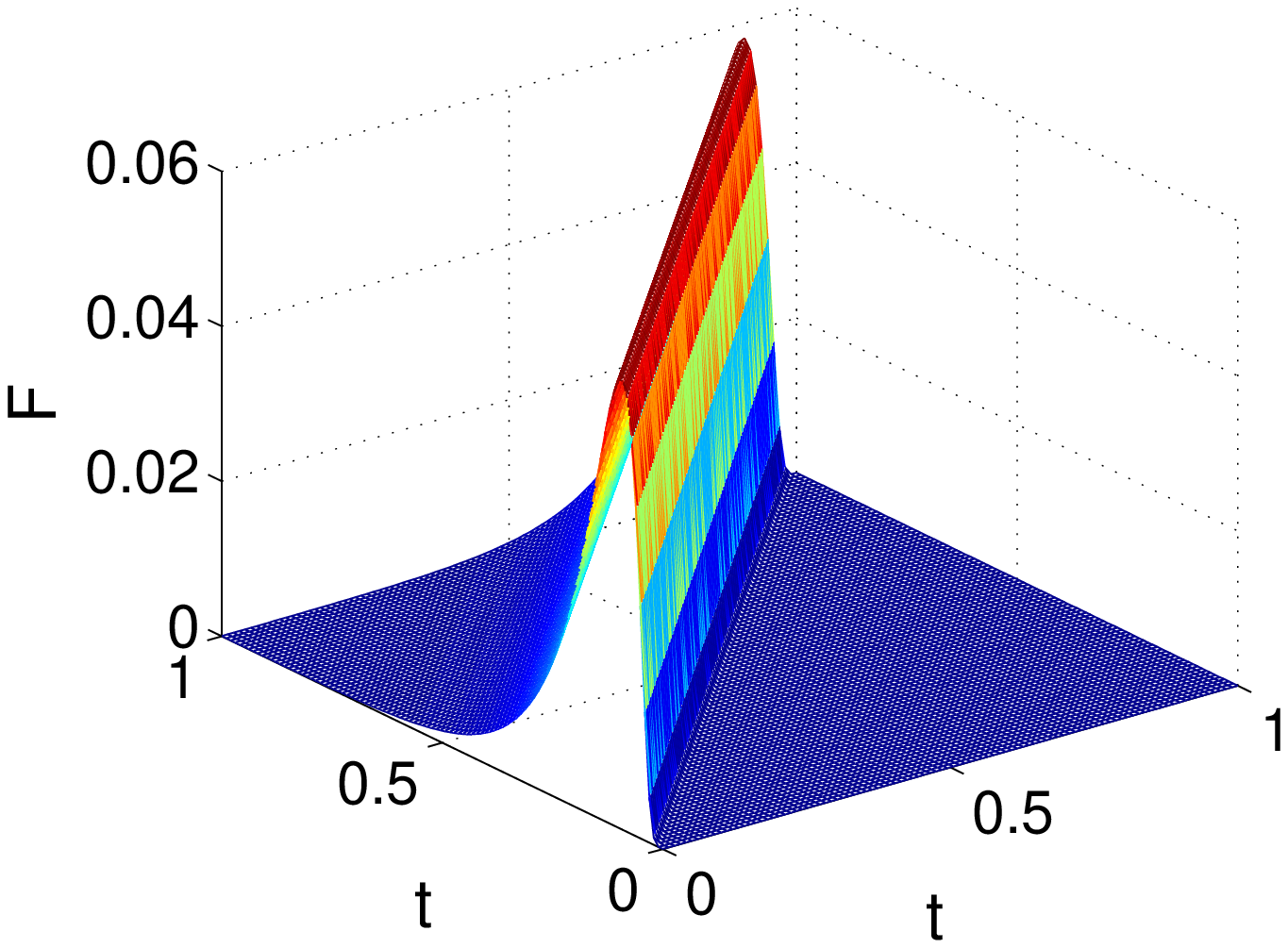}\\
     (c) $\alpha=3/4$ & (d) $\alpha=1$
  \end{tabular}
  \caption{The Jacobian map $F$ for $\alpha=1/4$, $1/2$, $3/4$ and $1$, from the interval $(0,1)$ itself.}\label{fig:sideways:jac}
\end{figure}

Physically this can also be observed in Fig. \ref{fig:sideways:jac}, where the forward map $F$ is from the Dirichlet boundary
condition $x=1$ to the flux boundary condition at $x=0$, in a piecewise linear finite element basis. Pictorially, the forward map $F$ is only
located in the upper left corner and has a triangular structure, which reflects the casual or Volterra nature of the sideways
problem for the fractional diffusion equation. We note that the casual structure should be utilized in developing reconstruction
techniques, via, e.g., Lavrentiev regularization \cite{Lamm:2000}. For small $\alpha$ values, e.g., $\alpha=1/4$, the finite
element basis at the right end point $x=1$ is almost instantly transported to the left end point $x=0$, whose magnitude is slightly decreased,
but with little diffusive effect, resulting a diagonally dominant forward map. However, as the fractional order $\alpha$ increases towards unity,
the diffusive effect eventually kicks in, and the information spreads over the whole interval. Further, for large $\alpha$ values,
it takes much longer time to reach the other side and there is a lag of information arrival, which explains the presence of tiny
singular values. The larger is the fractional order $\alpha$, the smaller is the magnitude, i.e., the less is the amount
of the information reached the other side. Hence, one feasible
approach is to recover only the boundary condition over a smaller subinterval of the measurement time interval. This idea underlies
one popular engineering approach, the sequential function specification method \cite{BeckBlackwellSt:1985,Liu:1996}.

The sideways problem for the classical diffusion has been extensively studied, and many efficient numerical methods have been
developed and analyzed \cite{Carasso:1982,Cannon:1984,Elden:1995,EldenBerntssonReginska:2000}. In the fractional case, however,
there are only a few works on numerical schemes, mostly for one-dimensional problems, and there seems no theoretical study on
stability etc. Murio \cite{Murio:2007,Murio:2008}
developed several numerical schemes, e.g., based on space marching and finite difference method, for the sideways problem, but
without any analysis. Qian \cite{Qian:2010} discussed about the ill-posedness of the quarter plane formulation of the sideways problem
using the Fourier analysis, based on which a mollifier method was proposed, with error estimates provided. In \cite{RundellXuZuo:2013},
the recovery of a nonlinear boundary condition from the lateral Cauchy data was studied using an integral equation approach, and a
convergent fixed point iteration method was suggested. The influence of the imprecise specification of the fractional order
$\alpha$ on the reconstruction was examined. Zheng and Wei \cite{ZhengWei:2011b} proposed a mollification method for the
quarter plane formulation of the sideways problem, by convoluting the fractional derivative with a smooth kernel, and derived error
estimates for the approximation, under a prior bounds on the solution. The Cauchy problem of
the time fractional diffusion has been numerically studied in \cite{ZhengWei:2012}. In particular, with the separation
of variables, a Volterra integral equation reformulation of the problem was derived, from which the ill-posedness of the Cauchy
problem follows directly. All these works are concerned with the one-dimensional case, and the high dimensional case has not been studied.

\subsection{Inverse source problem}
A third classical linear inverse problem for the diffusion equation is the inverse source problem, i.e.,
the recovery of the source term $f$ from lateral boundary data or final time data. Clearly, one piece of boundary
data or final time data alone is insufficient to uniquely determine a general source term, due to dimensionality
disparity. To restore the possible uniqueness, as usual, we look for only a space- or time-dependent component of
the source term $f$. With different combinations of
the data and source term, we get several different (and not equivalent) formulations of the inverse
source problems. Below we examine several of them briefly. By the linearity of the forward problem, we without
loss of generality, assume a zero initial data $v=0$ and a zero potential $q=0$ throughout this part.

\begin{figure}[hbt!]
  \centering
  \begin{tabular}{cc}
     \includegraphics[trim = .5cm .1cm 1cm 0cm, clip = true, width=7.5cm]{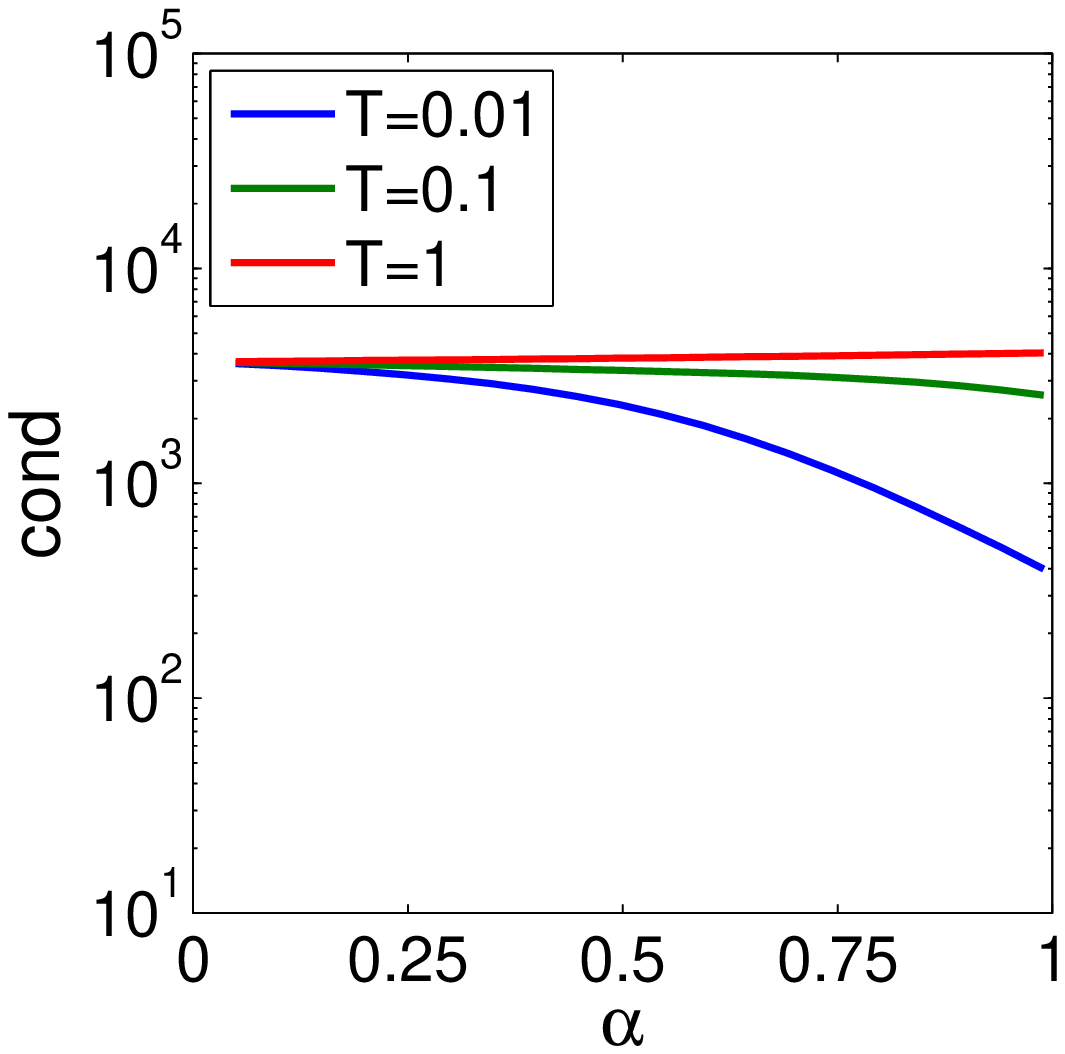} &
     \includegraphics[trim = .5cm .1cm 1cm 0cm, clip = true, width=7.5cm]{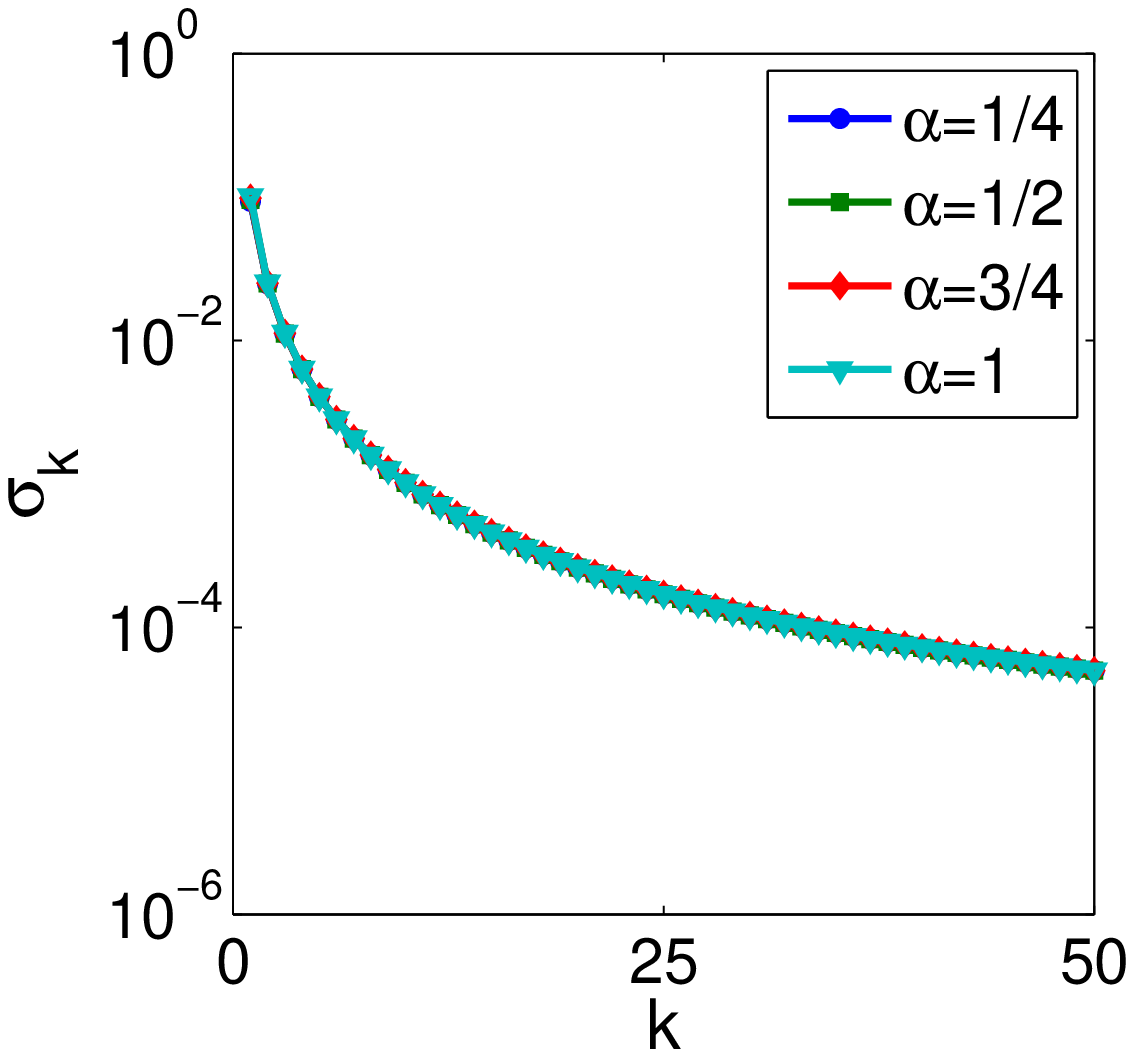}\\
     (a) condition number & (b) singular value spectrum
  \end{tabular}
  \caption{Numerical results for the inverse source problem with final
  time data and a space dependent source term.
  (a) The condition number v.s. the fractional order $\alpha$, and (b) singular value spectrum at $T=1$.}\label{fig:isp_spacespace}
\end{figure}

First, suppose we can measure the solution $u$ at the final time $t=T$, and aim at recovering either a
space dependent or time dependent component of the source term $f$. Like before, we resort to the
separation of variables. For the case of a space dependent only source term $f(x)$, the solution $u$ to
the forward problem is given by
\begin{equation*}
  \begin{aligned}
    u(t) &= \sum_{j=1}^\infty\int_0^t (t-\tau)^{\alpha-1}E_{\alpha,\alpha}(-\lambda_j(t-\tau)^\alpha)(f,\phi_j)\phi_j d\tau\\
        & = \sum_{j=1}^\infty \lambda_j^{-1}(1-E_{\alpha,1}(-\lambda_jt^\alpha))(f,\phi_j)\phi_j.
  \end{aligned}
\end{equation*}
Hence the measured data $g=u(T)$ is given by
\begin{equation*}
  g = \sum_{j=1}^\infty \frac{1}{\lambda_j}(1-E_{\alpha,1}(-\lambda_jT^\alpha))(f,\phi_j)\phi_j.
\end{equation*}
By taking inner product with $\phi_j$ on both sides, we arrive at the following representation of the
source term $f$ in terms of the measured data $g$
\begin{equation}\label{eqn:invsource-p1}
  f = \sum_{j=1}^\infty\lambda_j\frac{(g,\phi_j)}{1-E_{\alpha,1}(-\lambda_jT^\alpha)}\phi_j.
\end{equation}
By the complete monotonicity of the Mittag-Leffler function $E_{\alpha,1}(-t)$ on the positive real axis
$\mathbb{R}^+$ \cite{Pollard:1948}, we deduce
\begin{equation*}
 1 > E_{\alpha,1}(-\lambda_1T^\alpha)>E_{\alpha,1}(-\lambda_2T^\alpha),
\end{equation*}
and thus the formula \eqref{eqn:invsource-p1} is well defined for any $T>0$, and gives the precise
condition for the existence of a source term. Even with a modest value of the terminal time $T$,
the factor $1-E_{\alpha,1}(-\lambda_1T^\alpha)$ is close to unity for all small $\alpha$ values,
especially for those close to zero. Each frequency component $(f,\phi_j)$
differs from $(g,\phi_j)$ essentially by a factor $\lambda_j$, which amounts to two derivative loss
in space.  Actually one can show
\begin{equation*}
  \|f\|_{L^2(\Omega)}\leq c \|g\|_{H^2(\Omega)}.
\end{equation*}
This behavior is identical with that for the backward fractional diffusion. The statement
holds also for the inverse source problem for the classical diffusion case. This is not surprising, since
with a space dependent source term $f$, the solution $u$ to the forward problem can be split into the steady
solution $u_s$ and the decaying solution $u_d$, i.e., $ u = u_s + u_d$, where $u_s$ and $u_d$ solve
\begin{equation*}
  - u_{s}'' =  f, \ \ u_s(0) = u_s(1) = 0,
\end{equation*}
and
\begin{equation*}
  \partial_t^\alpha u_d - u_{d,xx} = 0,\ \ u_d(0,x) = f(x), \  \ u_d(0,t) = u_d(1,t) = 0,
\end{equation*}
respectively. By the decay behavior of the solution $u_d$, the steady state component $u_s$ is dominating, which
amounts to a two spatial derivative loss. This is fully confirmed by the numerical experiments, cf. Fig.
\ref{fig:isp_spacespace}. It is observed that the condition number is almost independent of the fractional
order $\alpha$, and it is of order $O(10^3)$, reflecting the mildly ill-posed nature of the
inverse problem. In particular, for large terminal time $T$, the singular value spectra are almost
identical for all fractional orders, decaying to zero at an algebraic rate, cf. Fig. \ref{fig:isp_spacespace}(b).

\begin{figure}[hbt!]
  \centering
  \begin{tabular}{cc}
     \includegraphics[trim = .5cm .1cm 1cm 0cm, clip = true, width=7.5cm]{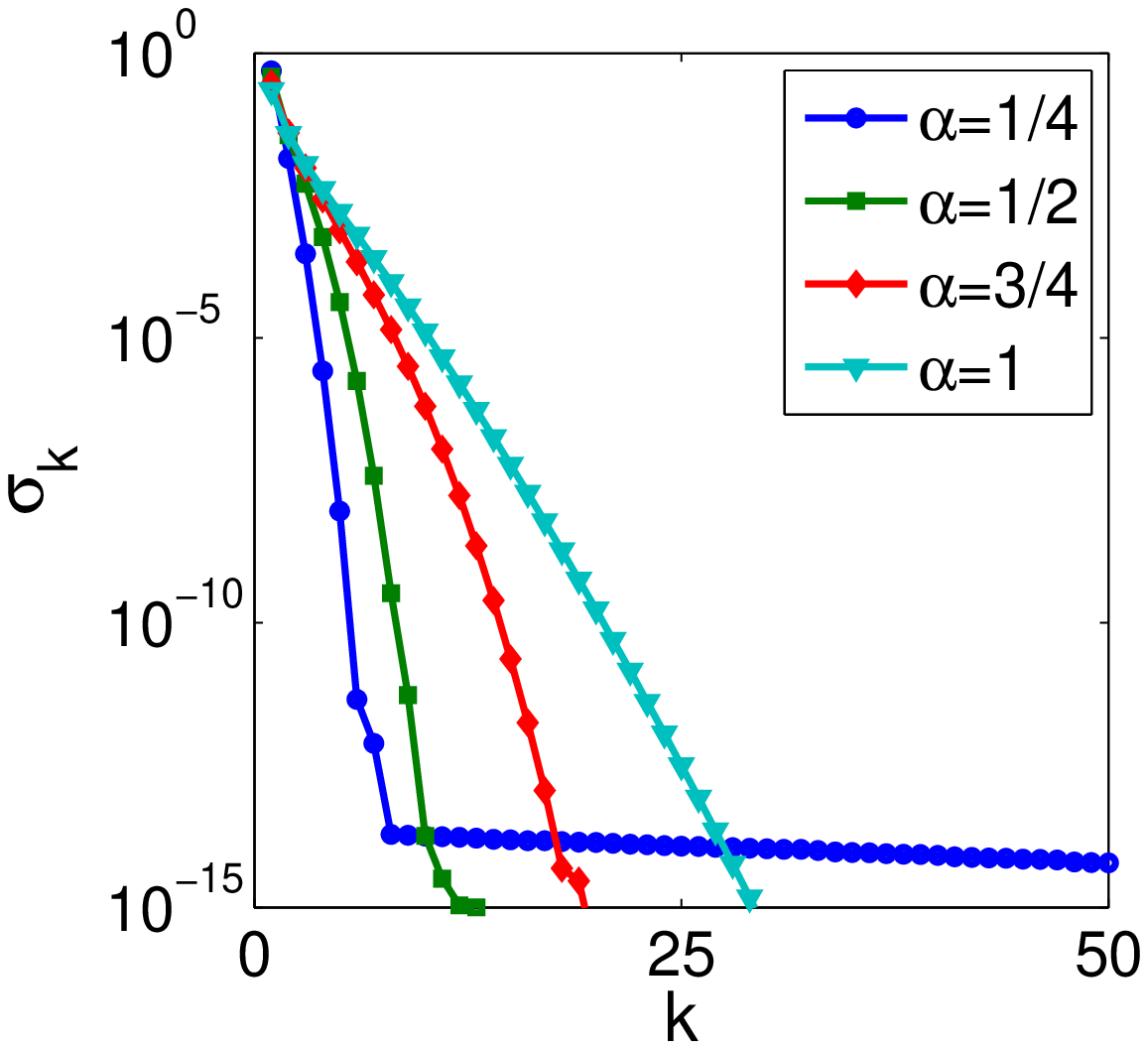} &
     \includegraphics[trim = .5cm .1cm 1cm 0cm, clip = true, width=7.5cm]{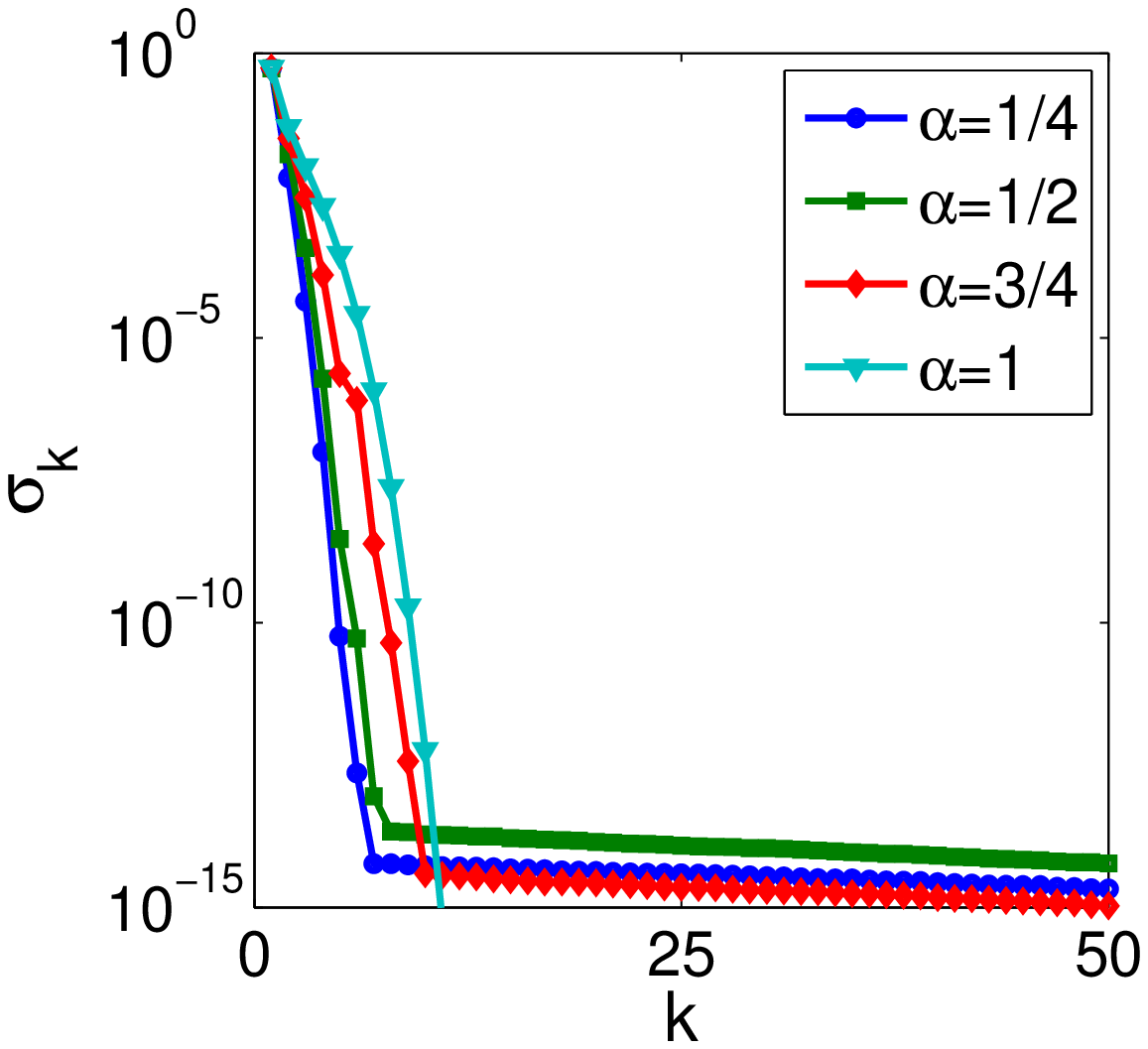}\\
     (a) $T=0.01$ & (b) $T=1$
  \end{tabular}
  \caption{The singular value spectrum at two different terminal times for the inverse source problem with
  a final time data at terminal time $T$, and $f(x,t)=xp(t)$, an unknown time dependent component $p(t)$.}\label{fig:isp_spacetime}
\end{figure}

Next we turn to the time dependent case, i.e., seeking a source term $f$ of the form $f(x,t)=p(t)q(x)$, with a known spacial
component $q(x)$, from the final time data $g=u(T)$. Mathematically, the inverse problem even for the classical diffusion equation
has not been completely analyzed. The inclusion of a nontrivial term $q(x)$ is important since without this there is nonuniqueness.
To see this, we take $u$ to satisfy  $u_t - u_{xx} = f(t)$ on $(0,1)\times(0,T)$ with initial data $u(x,0)=1$ and a homogeneous
Neumann boundary condition $-u_x(0,t)=u_x(1,t)=0$. Then one solution satisfying $u(x,T)=g(x)=1$ is given by $u(x,t)=1$ and $f
\equiv0$, but another is $u(x,t)=\cos(2\pi t/T)$ and $f=(-2\pi/T)\sin(2\pi t/T)$. Likewise, in the fractional case, we can take
$u=\cos(2\pi t/T)$ for the second solution and define $f$ to be its $\alpha$th order Djrbashian-Caputo fractional derivative in time.

Like previously, the solution $u$ to \eqref{eqn:fde-time} is given by
\begin{equation*}
  \begin{aligned}
    u(t) & = \sum_{j=1}^\infty \int_0^t (t-\tau)^{\alpha-1}E_{\alpha,\alpha}(-\lambda_j(t-\tau)^\alpha)p(\tau)d\tau
        (q,\phi_j)\phi_j.
  \end{aligned}
\end{equation*}
Hence the measured data $g(x)=u(x,T)$ is given by
\begin{equation*}
  g(x) = \sum_{j=1}^\infty \int_0^T (T-\tau)^{\alpha-1}E_{\alpha,\alpha}(-\lambda_j(T-\tau)^\alpha)p(\tau)d\tau
        (q,\phi_j)\phi_j(x).
\end{equation*}
By taking inner product with $\phi_j$ on both sides, we deduce
\begin{equation*}
  (g,\phi_j) = (q,\phi_j) \int_0^T (T-\tau)^{\alpha-1}E_{\alpha,\alpha}(-\lambda_j(T-\tau)^\alpha)p(\tau)d\tau.
\end{equation*}
In the case of $\alpha=1$, the formula recovers the relation
\begin{equation*}
  (g,\phi_j) = (q,\phi_j) \int_0^T e^{-\lambda_j(T-\tau)}p(\tau)d\tau,
\end{equation*}
which resembles a finite-time Laplace transform or moment problem, and thus severely smoothing, which renders
the inverse source problem severely ill-posed. Intuitively, the term $e^{-\lambda_j(T-t)}$ can only pick up the
information for $t$ close to the terminal time $T$, and for $t$ away from $T$, the information is severely
damped, especially for high frequency modes, which
leads to the severely ill-posed nature of the inverse problem. In the fractional case, the forward map $F$ from the
unknown to the data is clearly compact, and thus the problem is still ill-posed. However, the kernel $t^{\alpha-1}
E_{\alpha,\alpha}(-\lambda_jt^\alpha)$ is less smooth and decays much slower, and one might expect that the problem is less
ill-posed than the canonical diffusion counterpart. To examine the point, we present the numerical
results for the inverse problem in Fig. \ref{fig:isp_spacetime}. It is severely ill-posed
irrespective of the fractional order $\alpha$: the singular values decay exponentially to zero without
a distinct gap in the spectrum. In particular, for the terminal time $T=1$, the spectrum is almost identical for
all fractional orders $\alpha$. For small $T$, the singular values still decay exponentially, but the rate
is different: the smaller is the fractional order $\alpha$, the faster is the decay, cf. Fig. \ref{fig:isp_spacetime}(a).
Consequently, a few more modes of the source term $p(\tau)$ might be recovered. In other words, due to a slower local decay of the
exponential function $e^{-\lambda t}$, compared with the Mittag-Leffler function $t^{\alpha-1} E_{\alpha,\alpha}
(-\lambda t^\alpha)$, cf. Fig. \ref{fig:mitlef}(a), actually more frequency modes can be picked up by normal
diffusion than the fractional counterpart, cf. Fig. \ref{fig:isp_spacetime}(a). This indicates that with sufficiently
accurate data, at a small time instance, the sideways problem for normal diffusion may allow recovering more
modes, i.e., anomalous diffusion does not help solve the inverse problem.

In practice, the accessible data can also be the flux data at the end point, e.g., $x=0$ or $x=1$. We briefly
discuss the case of recovering a time dependent component $p(t)$ in the source term $f=q(x)p(t)$ from the
flux data at $x=0$. By repeating the preceding argument, the data $g:=-u_x(0,t)$ is related to the unknown $p(t)$ by
\begin{equation*}
  g(t) = -\sum_{j=1}^\infty \int_0^t (t-\tau)^{\alpha-1}E_{\alpha,\alpha}(-\lambda_j(t-\tau)^\alpha) p(\tau)d\tau (q(x),\phi_j)\phi_j'(0).
\end{equation*}

\begin{figure}[hbt!]
  \centering
  \begin{tabular}{cc}
     \includegraphics[trim = .5cm .1cm 1cm 0cm, clip = true, width=7.5cm]{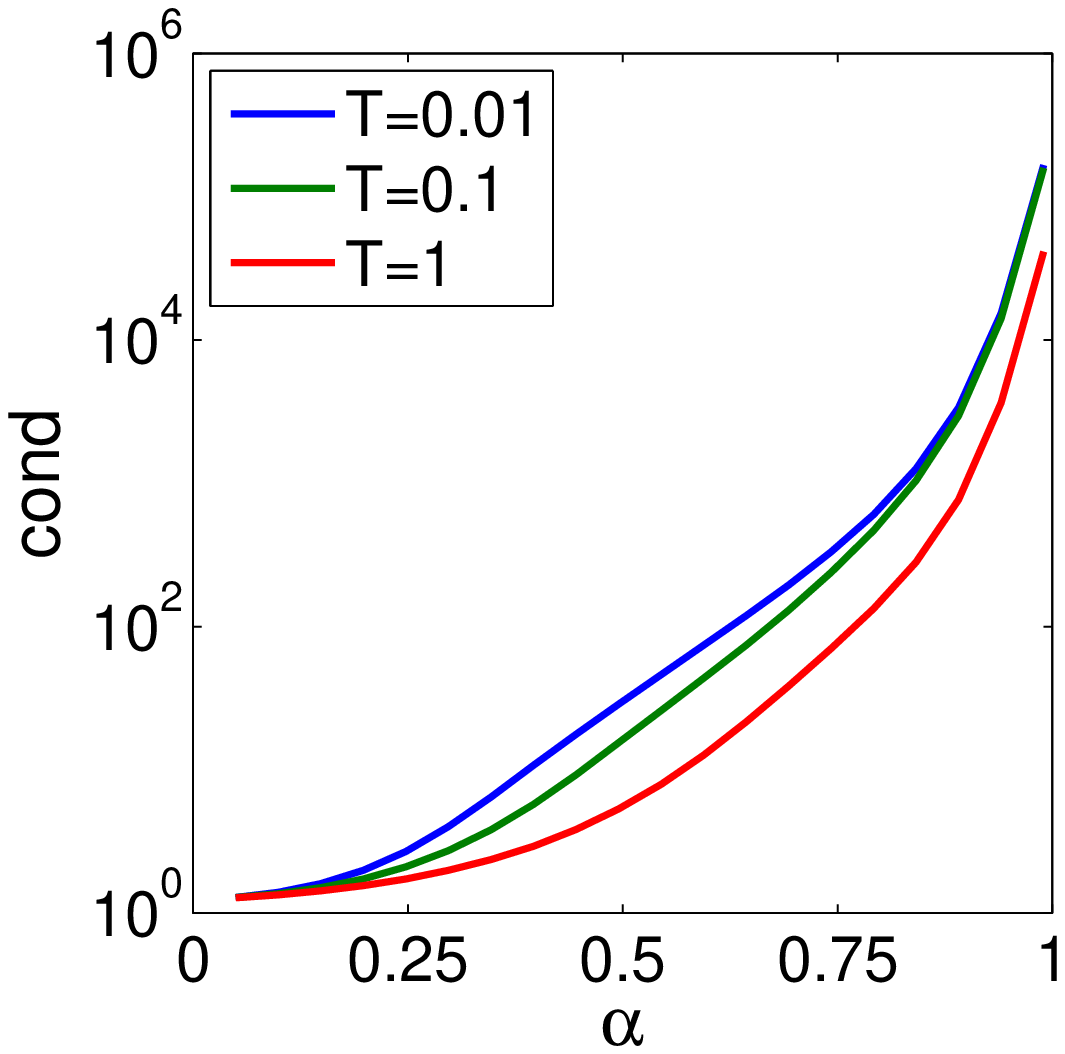} &
     \includegraphics[trim = .5cm .1cm 1cm 0cm, clip = true, width=7.5cm]{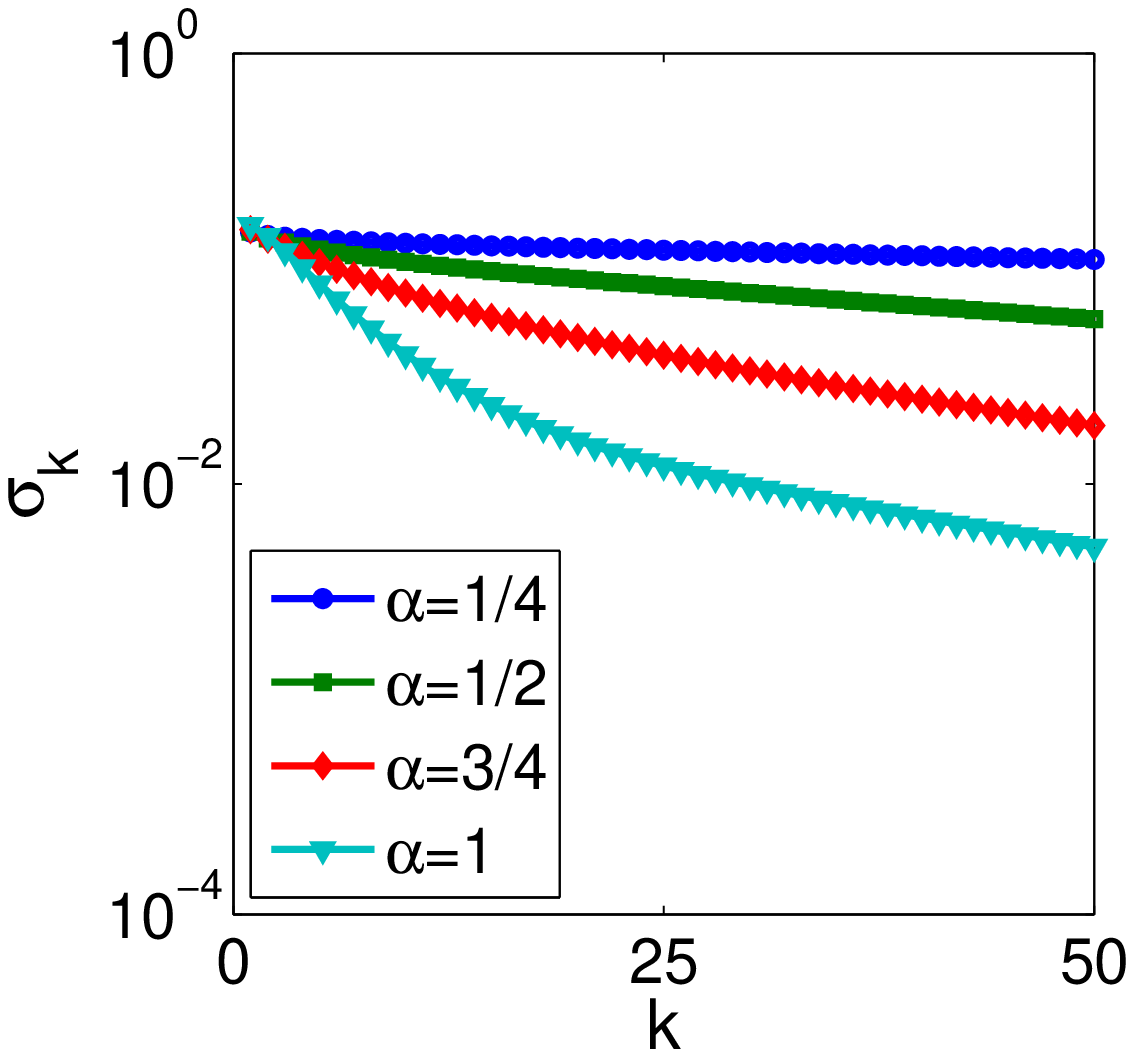}\\
     (a) condition number & (b) singular values
  \end{tabular}
  \caption{Numerical results for the inverse source problem with flux data at $x=0$ and $f(x,t)=xp(t)$, an unknown
  time dependent component $p(t)$. (a) The condition number of the discrete forward map and (b)
  singular value spectrum at $T=1$.}\label{fig:isp_timetime}
\end{figure}

In \cite[Theorem 4.4]{SakamotoYamamoto:2011}, a stability result was established for the recovery of the time
dependent component $p(t)$. Along the same line of thought, under reasonable assumptions, one can deduce that
\begin{equation*}
  \|p\|_{C[0,T]}\leq c\|\partial_t^\alpha g\|_{C[0,T]}.
\end{equation*}
The inverse problem roughly amounts to taking the $\alpha$th order Djrbashian-Caputo fractional derivative in time.
Hence as the fractional order $\alpha$ decreases from unity to zero, it becomes less and less ill-posed. For $\alpha$
close to zero, it is nearly well-posed, at least numerically. In other words, anomalous diffusion can mitigate the
degree of ill-posedness for the inverse problem. To illustrate the discussion, we present in Fig. \ref{fig:isp_timetime}
some numerical results, where the forward map $F$ is from the time dependent component $p(t)$ to the flux data $g(t)$
at $x=0$, both defined over the interval $[0,T]$, discretized using a continuous piecewise linear finite element basis.
The condition number of the discrete forward map $F$ decreases monotonically as the fractional order $\alpha$ decreases
from unity to zero, confirming the preceding discussions. Further, the terminal time $T$ does not affect the condition
number to a large extent.

It is widely accepted in inverse heat conduction that an inverse problem will be severely ill-posed
when the data and unknown are not aligned in the same space/time direction, and only mildly ill-posed
when they do align with each other. Our discussions with the inverse source problems indicate that the
observation remains valid in the time fractional diffusion case. In particular, although not presented, we
note that the inverse source problem of recovering a space dependent component from the lateral Cauchy
data is severely ill-posed for both fractional and normal diffusion. In the simplest case of a space
dependent only source term, it is mathematically equivalent to unique continuation, a well known
example of severely ill-posed inverse problems.

The inverse source problems for the classical diffusion equation have been extensively studied; see e.g.,
\cite{Cannon:1968,ImanuvilovYamamoto:1998,CannonDuChateau:1998}. Inverse source problems for FDEs have also
been numerically studied. Zhang and Xu \cite{ZhangXu:2011} established the unique recovery of a space
dependent source term in \eqref{eqn:fde-time} with pure Neumann boundary data and overspecified Dirichlet
data at $x=0$. This is achieved by an eigenexpansion and Laplace transform, and the uniqueness follows from
a unique continuation principle of analytic functions. Sakamoto and Yamamoto \cite{SakamotoYamamoto:2011a}
discussed the inverse problem of determining a spatially varying function of the source term by final
overdetermined data in multi-dimension, and established its well-posedness in the Hadamard sense except for
a discrete set of values of the diffusion constant, using an analytic Fredholm theory. Very recently, Luchko
et al \cite{LuchkoRundellYamamotoZuo:2013} showed the uniqueness of recovering a nonlinear source term from
the boundary measurement, and developed a numerical scheme of fixed point iteration type. Aleroev et al \cite{AleroevKiraneMalik:2013}
showed the uniqueness of recovering a space dependent source term from integral type observational data. Recently, there are
many numerical studies on this class of inverse problems. In \cite{WeiWang:2014b}, the numerical recovery of
a spatially varying function of the source term from the final time data in a general domain was studied using
a quasi-boundary value problem method; see also \cite{ZhangWei:2013,WangZhouWei:2013} related studies. Wang et
al \cite{WangYamamotoHan:2013} proposed to determine the space-dependent source term from the final time data
in multi-dimension using a reproducing kernel Hilbert space method.

\subsection{Inverse potential problem} Now we consider a nonlinear inverse coefficient
problem for the time fractional diffusion equation: given the final time data $g=u(T)$, find the
potential $q$ in the model
\begin{equation}\label{eqn:invpot}
  \partial_t^\alpha u - u_{xx} + qu = 0 \quad \mbox{ in } \Omega,
\end{equation}
with a homogeneous Neumann boundary condition and initial data $v$. The parabolic counterpart
has been extensively studied \cite{Isakov:1991,ChoulliYamamoto:1996,ChoulliYamamoto:1997},
where it was shown that the problem is nearly well-posed in the Hardamard sense in suitable H\"{o}lder space, under certain
conditions, using the strong maximum principle. In \cite{Isakov:1991}, an elegant fixed point
method was developed, and the monotone convergence of the method was established. It can be adapted
straightforwardly to the fractional case: given an initial guess $q^0$, compute the update $q^k$
recursively by
\begin{equation*}
  q^{k+1} = \frac{g''-\partial_t^\alpha u(x,T;q^k)}{g},
\end{equation*}
where the notation $u(x,T;q^k)$ denote the solution to problem \eqref{eqn:invpot} with the potential
$q^k$ at $t=T$. Since the strong maximum principle is still valid for the time fractional diffusion
equation \cite{Zacher:2013}, the scheme is monotonically convergent, under suitable conditions.

As the terminal time $T\to\infty$, the problem recovers a steady-state problem, and the scheme amounts
to twice numerical differentiation in space and converges within one iteration, provided that the data $g$ is accurate
enough. Hence, it is natural to expect that the convergence of the scheme will depends crucially on the
time $T$: the larger is the time $T$, the closer is the solution $u$ to the steady state solution; and thus
the faster is the convergence of the fixed point scheme. By Lemma \ref{lem:mitlef}, as the fractional order $\alpha$
approaches zero, the solution $u$ decays much faster around $t=0$ than the classical one, i.e.,
$\alpha=1$. In other words, the fractional diffusion problem can reach a ``quasi-steady state'' much faster
than the classical one, especially for $\alpha$ close to zero, and the scheme will then converge much faster.

\begin{figure}[hbt!]
  \centering
  \begin{tabular}{cc}
     \includegraphics[trim = .5cm .1cm 1cm 0cm, clip = true, width=7.5cm]{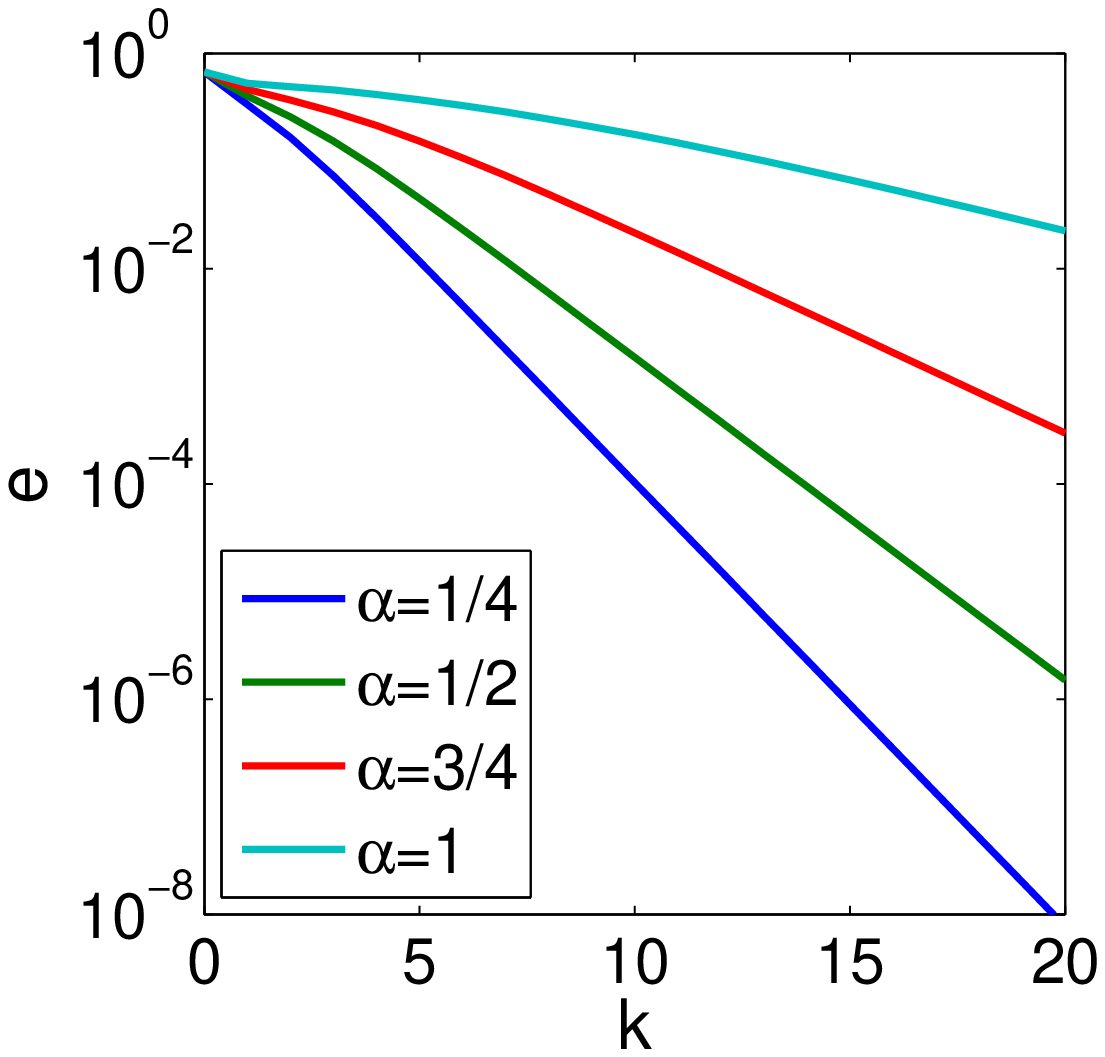} &
     \includegraphics[trim = .5cm .1cm 1cm 0cm, clip = true, width=7.5cm]{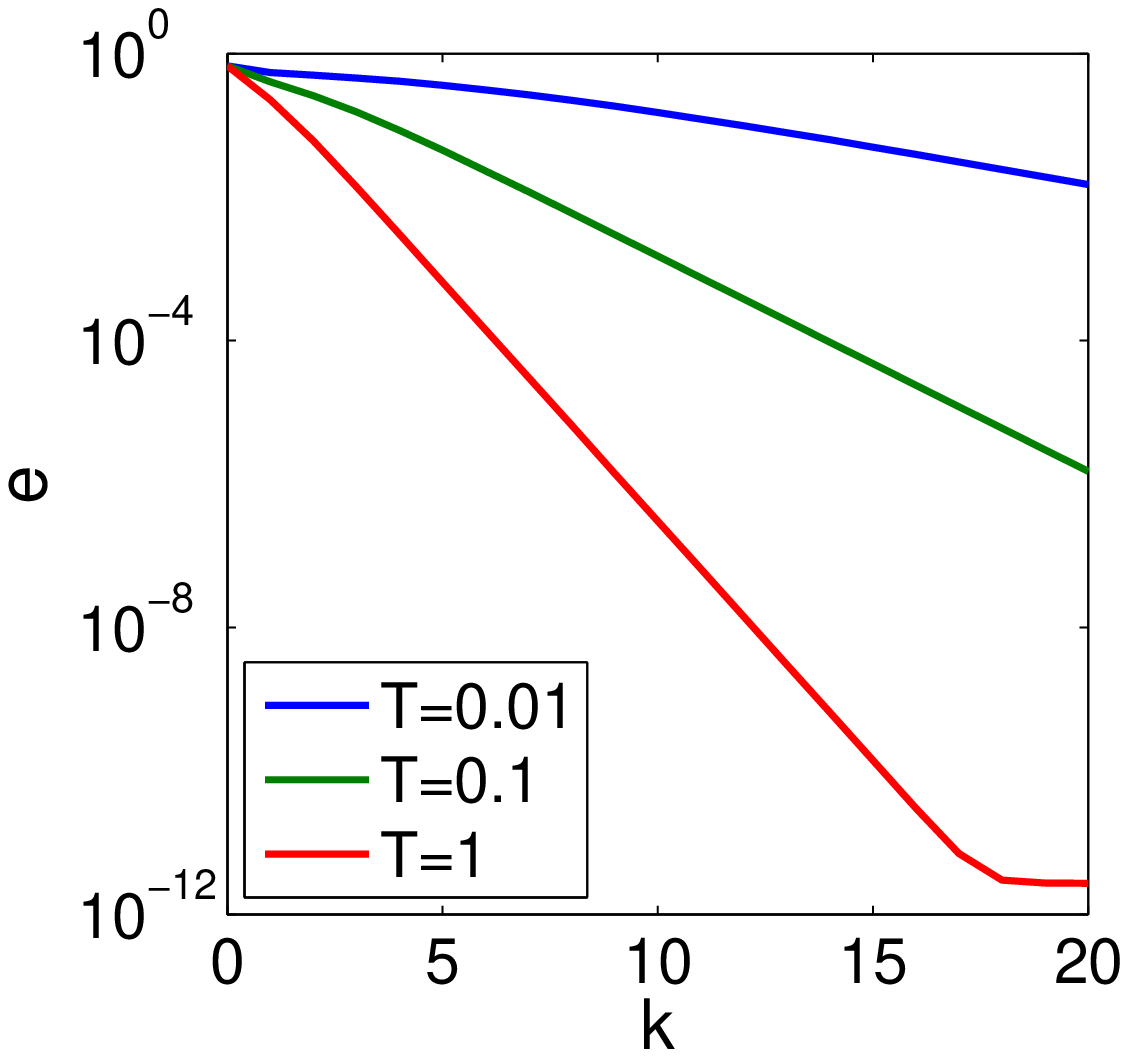}\\
     (a) error at $T=0.1$ & (b) error at $\alpha=1/2$
  \end{tabular}
  \caption{Numerical results, i.e., the relative $L^2(\Omega)$ error $e$, for the inverse potential problem
  from exact final time data at (a) $T=0.1$ and (b) $\alpha=1/2$.}\label{fig:invpot}
\end{figure}

To illustrate the point, we present in Fig. \ref{fig:invpot} some numerical results of reconstructing
a discontinuous potential $q=1 + 2x\chi_{[0,0.5]} + 2(1-\chi_{(0.5,1]}x)$ (with $\chi_S$ being the
characteristic function of the set $S$) from exact data, in order to illustrate the convergence behavior of the
fixed point scheme. In the figure, $e$ denotes the relative $L^2(\Omega)$ error. The numerical results fully
confirm the preceding discussions: at a fixed time $T$, the smaller is the fractional order $\alpha$, the faster
is the convergence; and at fixed $\alpha$, the larger is the time $T$, the faster is the convergence.
Numerically, one also observes the monotone convergence of the scheme.

Generally, the recovery of a coefficient in FDEs has not been extensively studied. Cheng et al \cite{ChengNakagawaYamamotoYamazaki:2009}
established the unique recovery of the fractional order $\alpha\in(0,1)$ and the diffusion coefficient from the
lateral boundary measurements. It represents one of the first mathematical works on invere problems for FDEs, and
has inspired many follow-up works. Yamamoto and Zhang \cite{YamamotoZhang:2012} established conditional stability
in determining a zeroth-order coefficient in a one-dimensional FDE with one half order Caputo derivative by a Carleman estimate.
Carleman estimates for time fractional diffusion were discussed in \cite{XuChengYamamoto:2011,ChengLinNakamura:2013,LinNakamura:2014}.
In \cite{MillerYamamoto:2013}, the unique determination of the spatial coefficient and/or the fractional order from
the data on a subdomain was shown for a positive initial condition. Wang and Wu \cite{WangWu:2014} studied the
simultaneous recovery of two time varying coefficients,
i.e., a kernel function and a source function, from the additional integral observation in multi-dimension, using a fixed point theorem.
All these works are concerned with the theoretical analysis, ant there are even fewer works on the numerical analysis
of related inverse problems. Li et al \cite{LiZhangJiaYamamoto:2013} suggested an optimal perturbation algorithm for
the simultaneous numerical recovery of the diffusion coefficient and fractional order in a one-dimensional time fractional
FDE. In \cite{JinRundell:2012b}, the
authors considered the identification of a potential term from the lateral flux data at one fixed time instance
corresponding to a complete set of source terms, and established the unique determination for ``small'' potentials.
Further, a Newton type method was proposed in \cite{JinRundell:2012b}, and its convergence was shown.

Even though our discussions have focused on time fractional diffusion, which involves one single fractional
derivative in time, it is also possible to consider equations where the time derivative involves multiple factional orders, i.e., $\sum_{k=1
}^m c_k \partial_t^{\alpha_k}$ for a sequence $\alpha_1>\alpha_2>\,\ldots>\alpha_m$ \cite{LiYamamoto:2013,JinLazarovLiuZhou:2014};
see \cite{LiImanuvilovYamamoto:2014} for some first uniqueness results for inverse coefficient problems in the multi-dimensional case.
Further extensions include the distributed-order, spatially and/or temporally variable-order and tempered fractional diffusion, to
better capture certain physical processes, for which, however, related
inverse problems have not been discussed at all.
\subsection{Fractional derivative as an inverse solution}
One of the very first undetermined coefficient problems for PDEs was discussed in the paper by B. Frank Jones Jr.
\cite{Jones:1962} (see also \cite[Chapter 13]{Cannon:1984}). This is to determine the coefficient $a(t)$ from
\begin{equation*}
  \begin{aligned}
    u_t &= a(t) u_{xx},\quad 0<x<\infty,\quad t>0\\
    u(x,0) &= 0,\quad -a(t)u_x(0,t) = g(t),\quad 0<t<T
  \end{aligned}
\end{equation*}
under the over-posed condition of measuring the ``temperature'' at $x=0$
\begin{equation*}
   u(0,t) = \psi(t)
\end{equation*}
In \cite{Jones:1962}, B. Frank Jones Jr. provided a complete analysis of the problem, by giving necessary and sufficient
conditions for a unique solution as well as determining the exact level of ill-conditioning. The key step
in the analysis is a change of variables and conversion of the problem to an equivalent integral equation
formulation. Perhaps surprisingly, this approach involves the use of a fractional derivative as we now show.

The assumptions are that $g$ is continuous and positive and $\psi$ is continuously differentiable with
$\psi(0)=0$ and $\psi'>0$ on $(0,T)$. In addition, the function $h(t)$ defined by
\begin{equation*}
h(t) = \frac{\sqrt{\pi} g(t)}{\int_0^t(t-\tau)^{-1/2}\psi'(\tau)d\tau}
\end{equation*}
satisfies $\lim_{t\to 0} h(t) = h_0 >0$. Note that $h$ is the ratio of the two data functions; the flux $g$
and the Djrbashian-Caputo derivative of order $1/2$ of $\psi$. If we define $h_i = \inf h$ and $h_s = \sup
h$ on $[0,T]$ and look at the space $\mathcal{G} := \{a \in C[0,T): h_i^2 \leq a(t) \leq h_s^2\}$, then it
was shown that any $a \in \mathcal{ G}$ satisfies the inverse problem must also solve the integral equation
\begin{equation*}
a(t) = \frac{\sqrt{\pi} g(t)}{\int_0^t\psi'(\tau)[\int_\tau^t a(s)\,ds ]^{-1/2}d\tau} =: \mathcal{ T}a,
\end{equation*}
and vice-versa. The main result in \cite{Jones:1962} is that the operator $\mathcal{T} $ has a unique fixed point on $\mathcal{G}$
and indeed $\mathcal{T} $ is monotone in the sense of preserving the partial order on $\mathcal{G}$, i.e., if $a_1\geq a_2$
then $\mathcal{T}a_1 \leq \mathcal{T}a_2$.

Given these developments, it might seem that a parallel construction for the time fractional diffusion counterpart, $\partial_t^\alpha u= a(t) u_{xx}$,
would be relatively straightforward but this seems not to be the case. The basic steps for the parabolic version require items
that just are not true in the fractional case, such as the product rule, and without these the above structure cannot be replicated
or at least not without some further ingenuity.

\section{Inverse problems for space fractional diffusion}\label{sec:space}

Now we turn to differential equations involving a fractional derivative in space. There are several possible
choices of a fractional derivative in space, e.g., Djrbashian-Caputo fractional derivative, Riemann-Liouville
fractional derivative, Riesz derivative, and fractional Laplacian \cite{Balakrishnan:1960}, which all have
received considerable attention. In recent years, the use of the fractional Laplacian is especially popular
in high-dimensional spaces, and admits a well-developed analytical theory. We shall focus on the left-sided
Djrbashian-Caputo fractional derivative $\DDC 0 \beta$, $\beta\in(1,2)$, and the one-dimensional case, and consider
the following four inverse problems: inverse Sturm-Liouville problem, Cauchy problem for a fractional elliptic
equation, backwards diffusion, and sideways problem.

\subsection{Inverse Sturm-Liouville problem}\label{ssec:eig}
First we consider the following Sturm-Liouville problem on the unit interval $\Omega=(0,1)$: find $u\in H_0^1(\Omega)\cap H^\beta(\Omega)$ and
$\lambda\in\mathbb{C}$ such that
\begin{equation}\label{eqn:fslp}
  - \DDC 0 \beta u + qu = \lambda u\quad \mbox{in } \Omega,
\end{equation}
with a homogeneous Dirichlet boundary condition $u(0)=u(1)=0$. A Sturm-Liouville problem of this form was
considered by Mkhitar M Djrbashian \cite{Dzrbasjan:1970,Djrbashian:1993} in 1960s to construct certain biorthogonal
basis for spaces of analytic functions; see also \cite{Nahusev:1977}. Like before, with $\beta=2$, it
recovers the classical Sturm-Liouville problem. In the case of a general potential $q$, in the fractional case,
little is known about the analytical properties of the eigenvalues and eigenfunctions. For the case of a
zero potential $q=0$, there are countably many eigenvalues $\{\lambda_j\}$ to \eqref{eqn:fslp}, which are zeros of the
Mittag-Leffler function $E_{\beta,2}(-\lambda)$. The corresponding eigenfunctions are given by $xE_{\beta,2}
(-\lambda_jx^\beta)$. Using the exponential asymptotics on the Mittag-Leffler function in Lemma \ref{lem:mitlef},
one \cite{Sedletskii:1994,JinRundell:2012b} can show that asymptotically, the eigenvalues $\lambda_j$ are distributed as
\begin{equation*}
  |\lambda_j| \sim (2\pi j)^\beta \quad \mbox{and}\quad
  \mathrm{arg}(\lambda_j) \sim \frac{(2-\beta)\pi}{2}.
\end{equation*}
Hence, for any $\beta\in(1,2)$, there are only a finite number of real eigenvalues
to \eqref{eqn:fslp}, and the rest appears as complex conjugate pairs.

\begin{figure}[hbt!]
  \centering
  \begin{tabular}{cc}
  \includegraphics[trim = 1cm .1cm 1cm 0cm, clip = true,width=7.5cm]{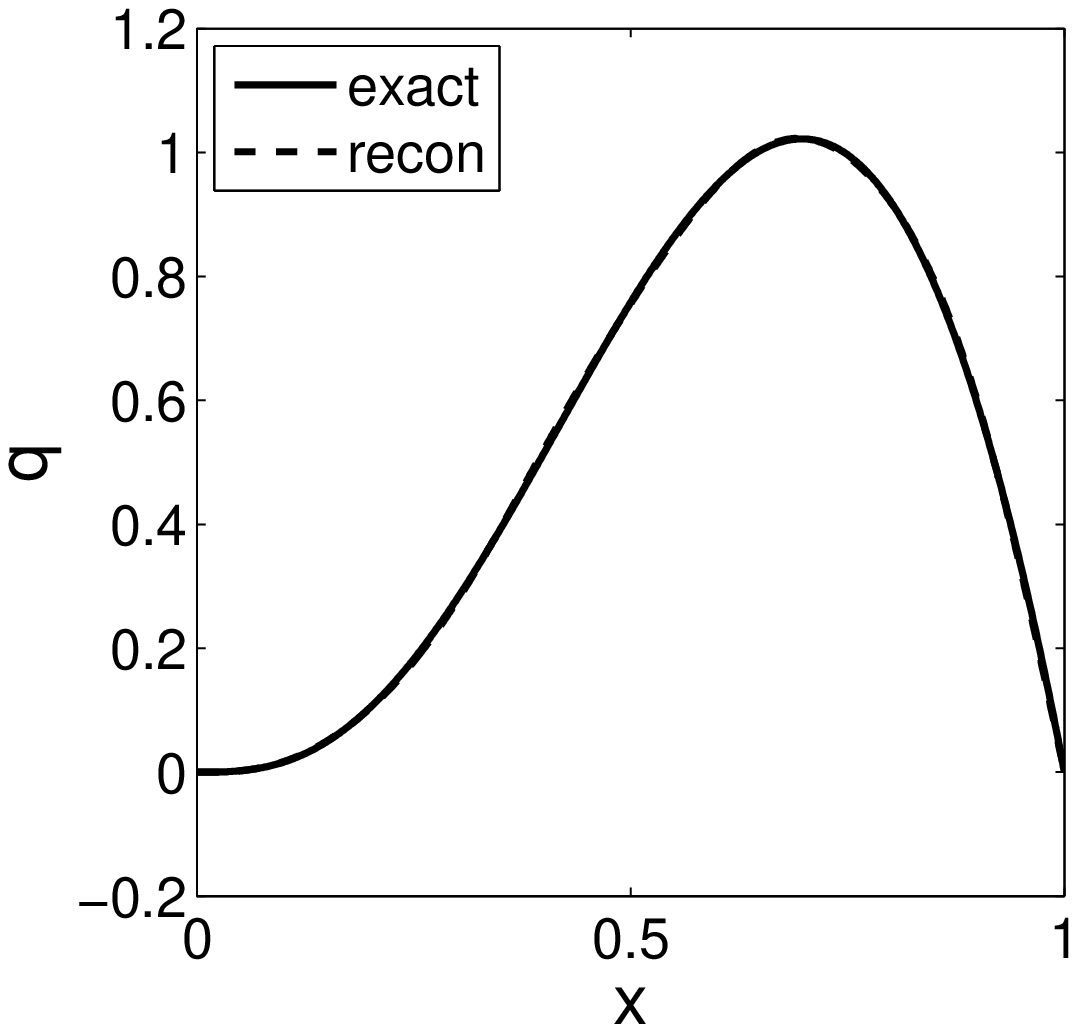} &
  \includegraphics[trim = 1cm .1cm 1cm 0cm, clip = true,width=7.5cm]{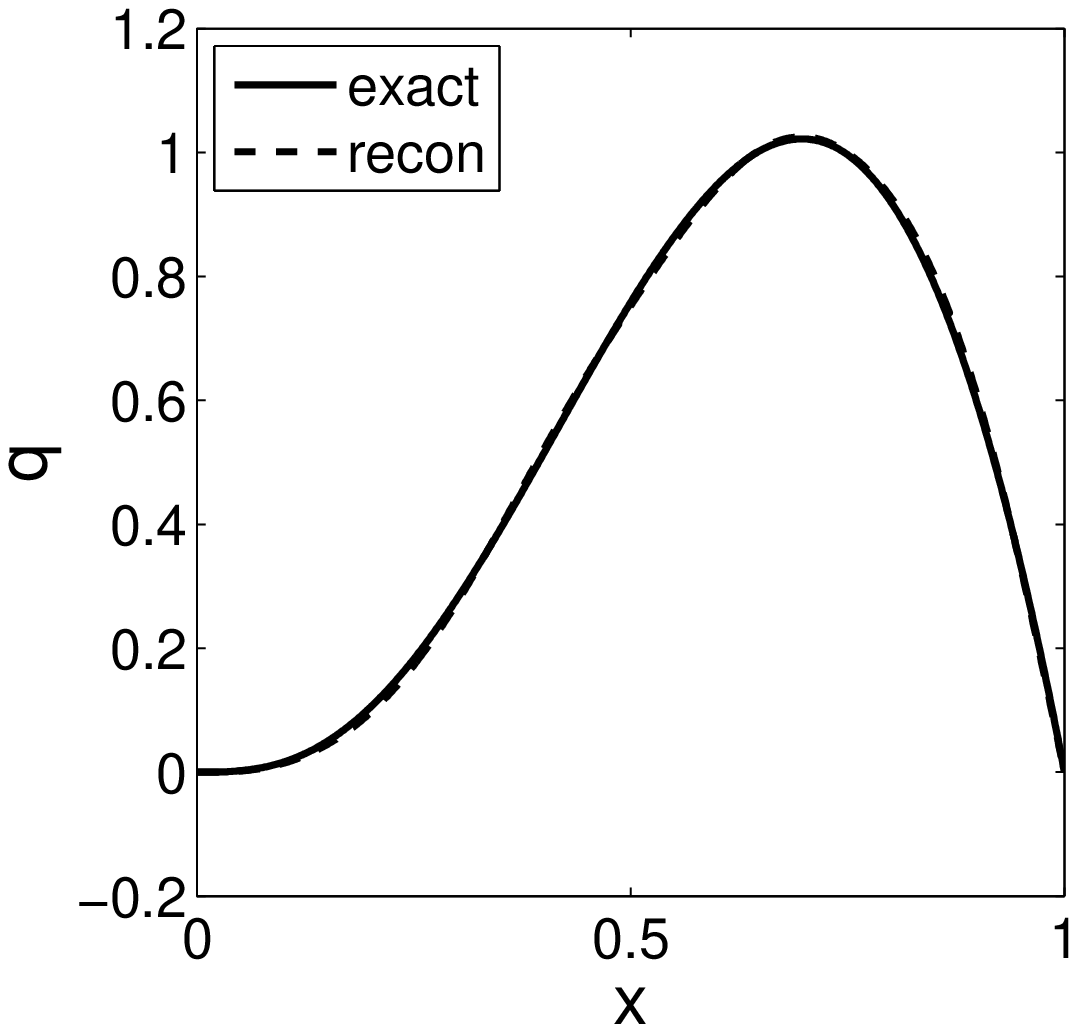}\\
  (a) $\beta=5/3$ & (b) $\beta=7/4$
  \end{tabular}
  \caption{Numerical results for the inverse Sturm-Liouville problem with a Djrbashian-Caputo derivative
  for (a) $\beta=5/3$ and (b) $\beta=7/4$. The
  reconstructions are computed from the first eight eigenvalues (in absolute value) using a frozen
  Newton method \cite{JinRundell:2012b}.}\label{fig:eig:cap}
\end{figure}

It is well known that eigenvalues contain valuable information about the boundary value problem. For example
it is known that the sequence of Dirichlet eigenvalues can uniquely determine a potential $q$ symmetric with
respect to the point $x=1/2$, and together with additional spectral information, one can uniquely determine
a general potential $q$; see \cite{ChadanColtonPaivrintaRundell:1997,RundellSacks:1992} for an overview of results on the classical
inverse Sturm-Liouville problem. In the fractional case, the eigenvalues are generally genuinely complex, and
a complex number may carry more information than a real one. Thus one naturally wonders whether these complex eigenvalues
do contain more information about the potential. Numerically the answer is affirmative. To illustrate this,
we show some numerical reconstructions in Fig. \ref{fig:eig:cap}, obtained by using a frozen Newton method
and representing the sought-for potential $q$ in Fourier series \cite{JinRundell:2012b}. The Dirichlet eigenvalues can be
computed efficiently using a Galerkin finite element method \cite{JinLazarovPasciak:2013a}. One observes that
one single Dirichlet spectrum can uniquely determine a general potential $q$. Unsurprisingly, as the fractional
order $\beta$ tends two, the reconstruction becomes less and less accurate, since in the limit  $\beta=2$,
the Dirichlet spectrum cannot uniquely determine a general potential $q$. Theoretically, the surprising uniqueness
in the fractional case remains to be established.

Naturally, one can also consider the Riemann-Liouville case:
\begin{equation}\label{eqn:fslp:riem}
  -{\DDR 0 \beta u} + qu = \lambda u \quad \mbox{ in } \Omega,
\end{equation}
with $u(0)=u(1)=0$. Like before, little is known about the analytical properties of the eigenvalues and eigenfunctions.
In the case of a zero potential $q=0$, there are countably many eigenvalues to \eqref{eqn:fslp:riem}, which
are zeros of the Mittag-Leffler function $E_{\beta,\beta}(-\lambda)$, and the corresponding eigenfunctions are
given by $x^{\beta-1}E_{\beta,\beta}(-\lambda_jx^\beta)$. Further, the asymptotics of the eigenvalues are still
valid. Hence, for any $\beta\in(1,2)$, there are only a finite number of real eigenvalues to \eqref{eqn:fslp:riem},
and the rest appears as complex conjugate pairs.

\begin{figure}[hbt!]
  \centering
  \begin{tabular}{cc}
  \includegraphics[trim = 1cm .1cm 1cm 0cm, clip = true,width=7.5cm]{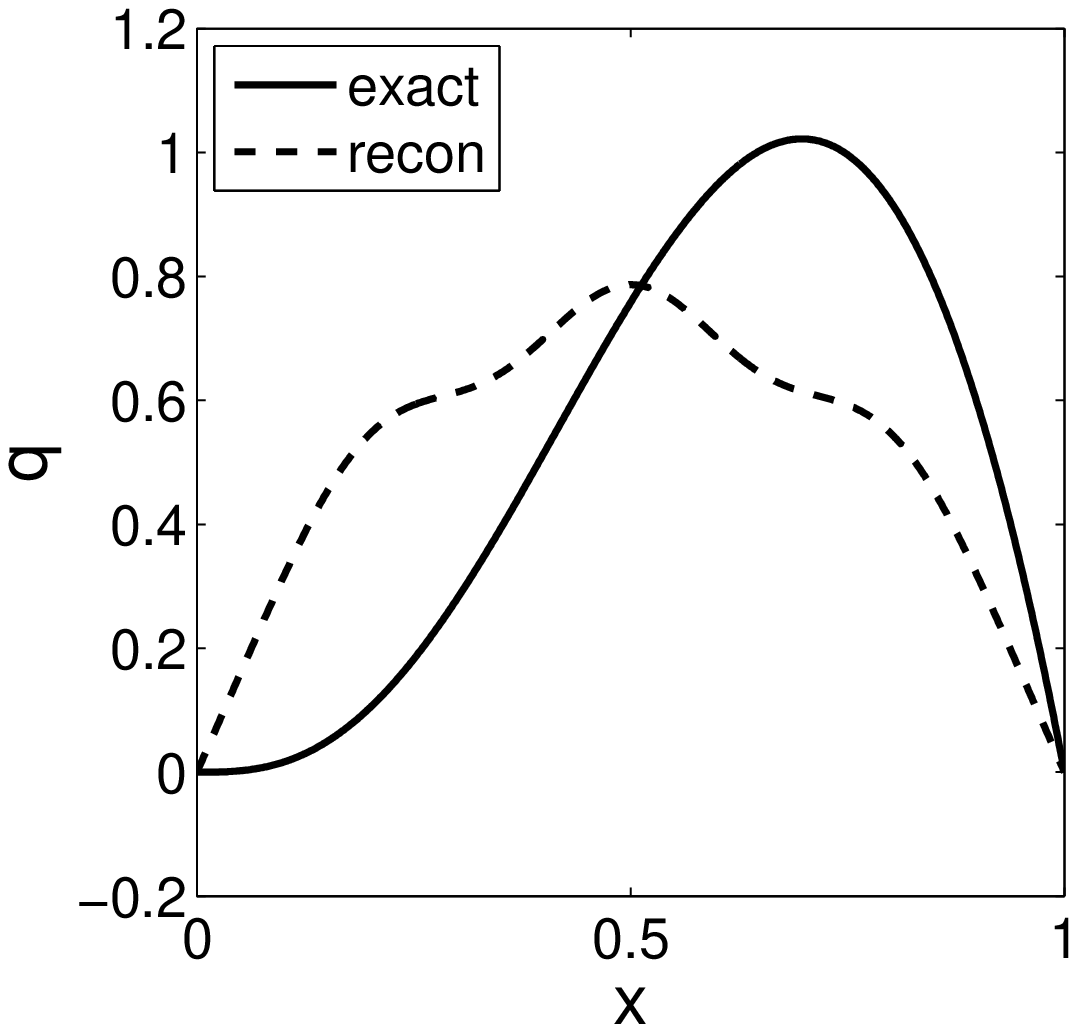}&
  \includegraphics[trim = 1cm .1cm 1cm 0cm, clip = true,width=7.5cm]{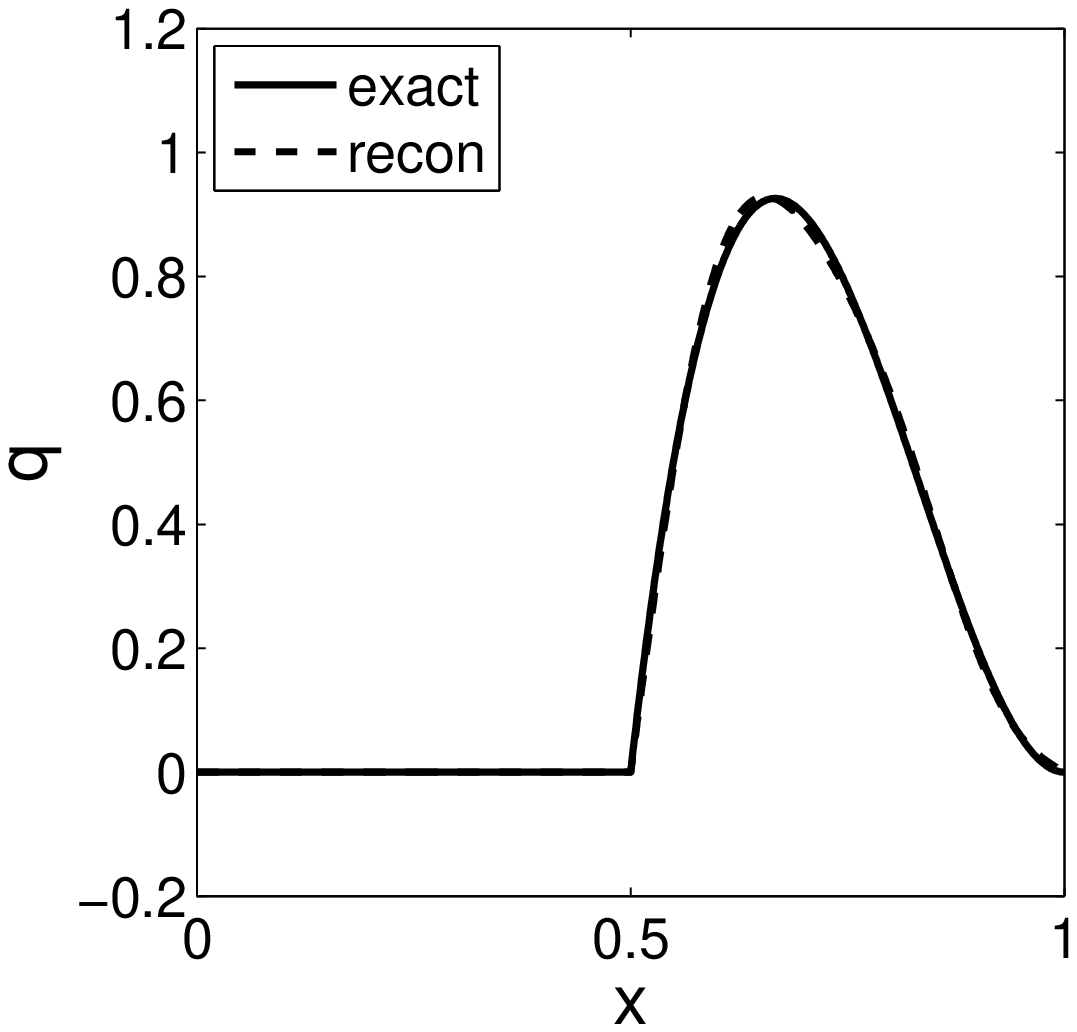}\\
  (a) whole interval & (b) half interval
  \end{tabular}
  \caption{Numerical results for the inverse Sturm-Liouville problem with a Riemann-Liouville fractional
  derivative of order $\beta=4/3$. The reconstructions are computed from the first eight eigenvalues (in absolute value) using a frozen
  Newton method \cite{JinRundell:2012b}.}\label{fig:eig:riem}
\end{figure}

The numerical results from the Dirichlet spectrum in the Riemann-Liouville case are shown in Fig. \ref{fig:eig:riem}.
For a general potential $q$, the reconstruction represents only the symmetric part, which is drastically different
from the Djrbashian-Caputo case, but identical with that for the classical Sturm-Liouville problem. Further,
if we assume that the potential $q$ is known on the left half interval, then the Dirichlet spectrum allows uniquely
reconstructing the potential $q$ on the remaining half interval, cf., Fig. \ref{fig:eig:riem}(b). These results indicate
that in the Riemann-Liouville case the complex spectrum is not more informative than the classical Sturm-Liouville problem. The precise
mechanism underlying the fundamental differences between the Djrbashian-Caputo and Riemann-Liouville cases awaits further study.
However, as in the classical case $\beta=2$, one can show that the linearized derivative of the map $q \to u(1;\lambda,
q)$ around $q=0$ cannot span more than the subspace of even functions in $L^2(\Omega)$.

In general, the Sturm-Liouville problem with a fractional derivative remains completely elusive, and numerical
methods such as finite element method \cite{JinLazarovPasciakRundell:2014} provide a valuable (and often the
only) tool for studying its analytical properties. For a variant of the fractional Sturm-Liouville problem, which
contains a fractional derivative in the lower-order term, Malamud \cite{Malamud:1994} established the existence
of a similarity transformation, analogous to the well-known Gel'fand-Levitan-Marchenko transformation, and also the unique
recovery of the potential from multiple spectra. In the classical case, the Gel'fand-Levitan-Marchenko transformation
lends itself to a constructive algorithm \cite{RundellSacks:1992}; however, it is unclear whether this is true
in the fractional case. In \cite{JinRundell:2012b}, the authors proposed a Newton type
method for reconstructing the potential, which numerically exhibits very good convergence behavior. However, a
rigorous convergence analysis of the scheme is still missing. Further, the uniqueness and nonuniqueness issues of related inverse
Sturm-Liouville problems are outstanding. Last, as noted above, there are other possible choices of the space fractional
derivative, e.g., fractional Laplacian and Riesz derivative. It is unknown whether the preceding observations
are valid for these alternative derivatives.

\subsection{Cauchy problem for fractional elliptic equation}
One classical elliptic inverse problem is the Cauchy problem for the Laplace equation, which plays a fundamental role
in the study of many elliptic inverse problems \cite{Isakov:2006}. A first example was given by Jacques Hadamard
\cite{Hadamard:1923} to illustrate the severe ill-posedness of the Cauchy problem, which motivated him to introduce
the concept of well-posedness and ill-posedness for problems in mathematical physics. So a natural question is whether
the Cauchy problem for the fractional elliptic equation is also as ill-posed? To illustrate this, we consider the
following fractional elliptic problem on the rectangular domain $\Omega=\{(x,y)\in\mathbb{R}^2: 0<x<1,0<y<1\}$
\begin{equation}\label{eqn:cauchy}
    {\DDC 0 \beta} u + {\,_{0}^C \kern -.2em  D^{\beta}_{\kern -.1em y}} u  = 0 \mbox{ in } \Omega,
\end{equation}
with the fractional order $\beta\in(1,2)$ and the Cauchy data
\begin{equation*}
    u(x,0)  = g(x) \quad \mbox{and}\quad
    \frac{\partial u}{\partial \nu}(x,0) = h(x), \  \ 0<x<1
\end{equation*}
where $\nu$ is the unit outward normal direction. With $\beta=2$, it recovers the Cauchy problem for the Laplace equation.
By applying the separation of variables, we assume that $u(x,y) = \phi(x)\psi(y)$,
which directly gives for some scalar $\lambda\in \mathbb{C}$ that
\begin{equation*}
  \begin{aligned}
   {\DDC 0 \beta} \phi(x) &=- \lambda \phi(x),\\
   {\,_{0}^C \kern -.2em  D^{\beta}_{\kern -.1em y}} \psi(y) &= \lambda \psi(y),
  \end{aligned}
\end{equation*}
Let $(\lambda_j,\phi_j)$ be a Dirichlet eigenpair of the Caputo derivative operator $-{\DDC 0 \beta}$ on the unit
interval $D=(0,1)$, i.e., $\phi_j(x) = xE_{\beta,2}(-\lambda_jx^\beta)$, and $|\lambda_j|\to\infty$
as $j\to\infty$ \cite{JinRundell:2012}; see Section \ref{ssec:eig} for further details. With the choice $\phi=\phi_j$ and the
Cauchy data pair $(g,h_j)=(0,-xE_{\beta,2}(\lambda_jx^\beta)/\lambda_j)$, the component $\psi_j$
satisfies
\begin{equation*}
   {\,_{0}^C \kern -.2em  D^{\beta}_{\kern -.1em y}} \psi_j(y) = \lambda_j\psi_j(y) \mbox{ in } y\in(0,\infty),
\end{equation*}
with the initial condition $\psi_j(0)=0$ and $\frac{d}{d y}\psi_j(0)=1/\lambda_j$. Using the relation $ \frac{d}{dx}
x^{\gamma-1}E_{\beta,\gamma}(\lambda x^\beta) = \lambda x^{\gamma-2}E_{\beta,\gamma-1}(\lambda x^\beta)$ \cite[pp.
46]{KilbasSrivastavaTrujillo:2006}, we deduce that the solution $\psi_j$ to the fractional ordinary differential
equation is given by
\begin{equation*}
  \psi_j(y) = yE_{\beta,2}(\lambda_jy^\beta)/\lambda_j.
\end{equation*}
Hence, $u_j(x,y) = xE_{\beta,2}(-
\lambda_jx^\beta)yE_{\beta,2}(\lambda_jy^\beta)/\lambda_j^2$ is a solution to the Cauchy problem
with $g=0$ and $h_j(x)=-xE_{\beta,2}(-\lambda_jx^\beta)/\lambda_j$. By the exponential asymptotics of the Mittag-Leffler
function, cf. Lemma \ref{lem:mitlef}, we deduce that $h_j(x)\to 0$ as $j\to\infty$, whereas for any $y>0$, the
solution $u_j(x,y)\to\infty$ as $j\to\infty$, in view of the exponential growth of the Mittag-Leffler function $E_{\beta,
2}(z)$, cf. Lemma \ref{lem:mitlef}. This indicates that the Cauchy problem for the fractional elliptic equation is also
exponentially ill-posed. However, the interesting question of the degree of ill-posedness, in comparison with the classical
case, is unclear and certainly worthy
of further study. Further, we note that the numerical solution of the fractional elliptic equation \eqref{eqn:cauchy} is
highly nontrivial, and there seems no efficient yet rigorous solver available in the literature, due to a lack
of solution theory for such problems.

\subsection{Backward problem}
\begin{figure}[hbt!]
  \centering
  \begin{tabular}{cc}
     \includegraphics[trim = .5cm .1cm 1cm 0cm, clip = true, width=7.5cm]{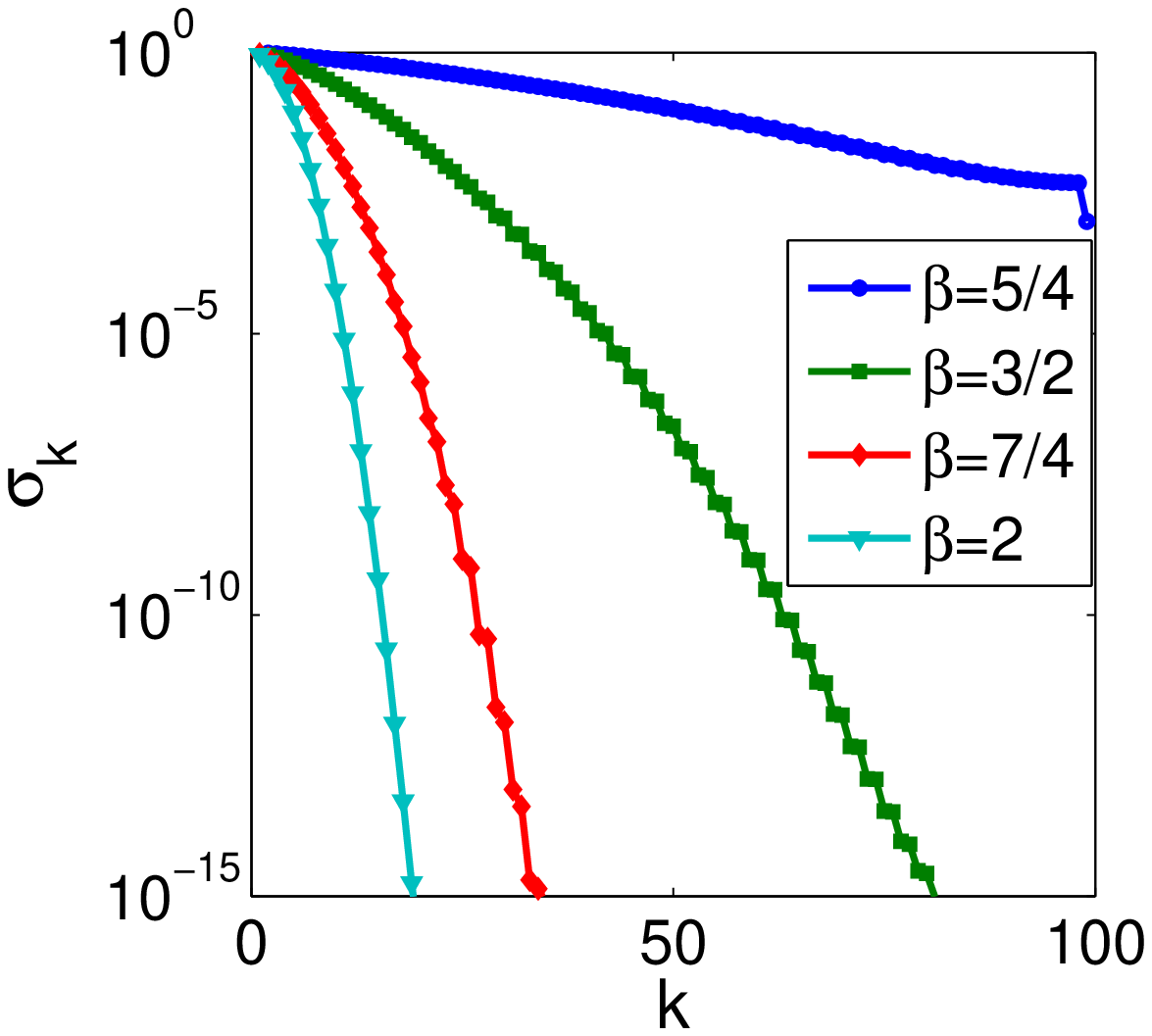} &
     \includegraphics[trim = .5cm .1cm 1cm 0cm, clip = true, width=7.5cm]{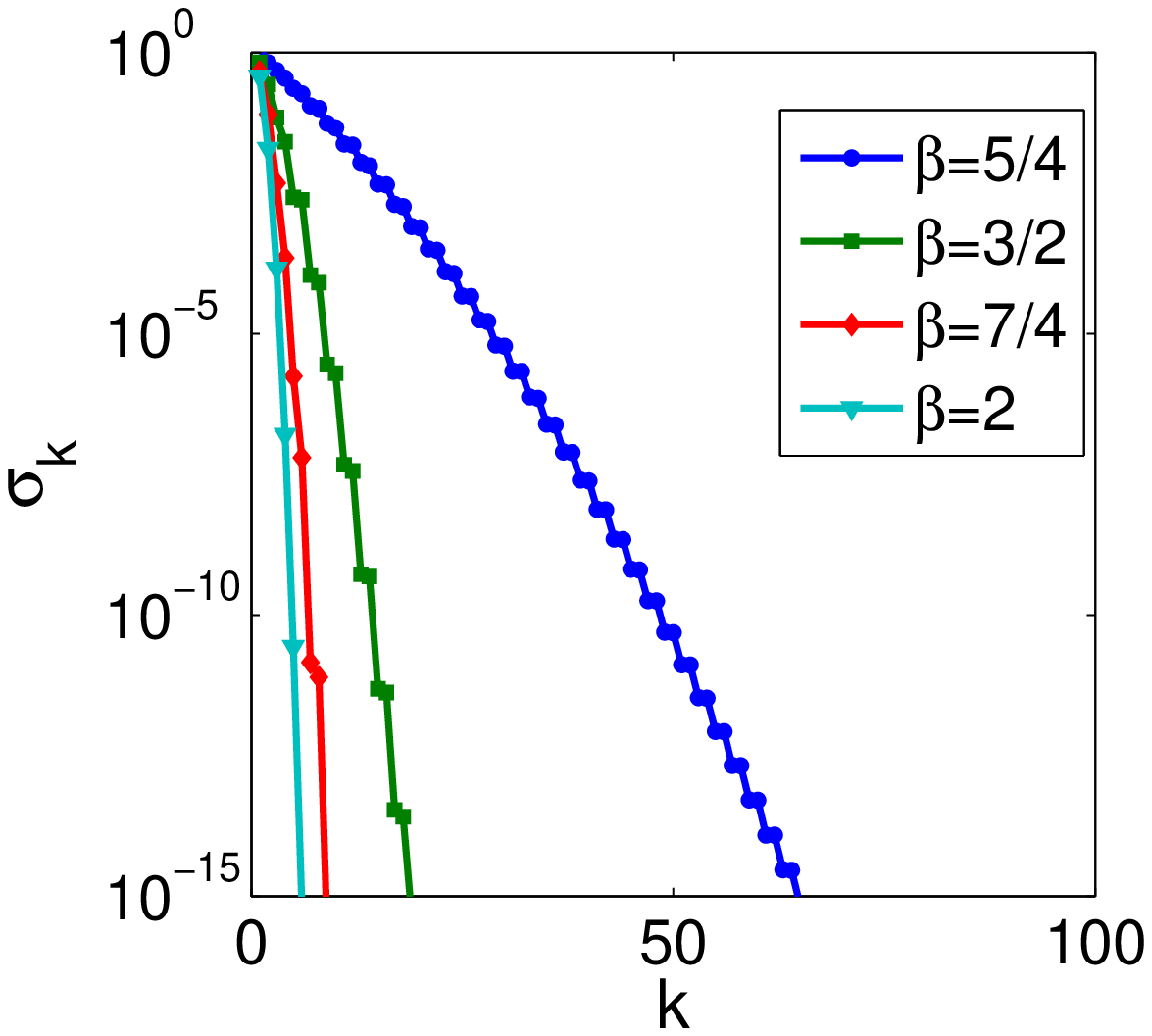}\\
     (a) $T=0.01$ & (b) $T=0.1$
  \end{tabular}
  \caption{Numerical results for the space fractional backward diffusion problem, the singular value spectrum at
  two different time instances, (a) $T=0.01$ and (b) $T=0.1$.}\label{fig:backward_spacefrac}
\end{figure}

Now we return to the backward diffusion problem with fractional derivatives in the space variable(s).
Let $\Omega=(0,1)$ be the unit interval. Then the one-dimensional space fractional diffusion equation is
given by
\begin{equation*}
  u_t - {\DDC 0 \beta } u = 0,\quad (x,t)\in \Omega\times (0,\infty),
\end{equation*}
where the fractional order $\beta\in(1,2)$. The equation is equipped with the following initial condition
$u(x,0)=v$ and zero boundary condition $u(0,t)=u(1,t)=0$. The backward problem is: given the final time
data $g(x)=u(x,T)$, find the initial data $v$. Since the Djrbashian-Caputo derivative operator
$\DDC 0 \beta $ with the zero Dirichlet boundary is sectorial on suitable spaces \cite{ItoJinTakeuchi:2014},
the existence of a solution $u$ follows from the analytic semigroup theory \cite{Pazy:1992,ItoKappel:2002},
and formally it can be represented by
\begin{equation*}
  u(t) = e^{-At}v,
\end{equation*}
where $A$ is the representation of the Djrbashian-Caputo derivative operator $- \DDC 0 \beta$ on its domain. Formally, the solution
$v$ to the space fractional backward problem is given by
\begin{equation*}
  v = e^{AT}g.
\end{equation*}
In case of $\beta=2$, using the eigenpairs $\{(\lambda_j,\phi_j)\}$ and the $L^2(\Omega)$ orthogonality of the
eigenfunctions $\{\phi_j\}$, it recovers the well known formula
\begin{equation*}
  v = \sum_{j=1}^\infty e^{\lambda_jT}(g,\phi_j) \phi_j.
\end{equation*}
The growth factor $e^{\lambda_jT}$ explains the severely ill-posed nature of the inverse problem.
In the fractional case, such an explicit representation is no longer available since the corresponding eigenfunctions $\{\phi_j\}$
are not orthogonal in $L^2(\Omega)$ (actually they can be almost linearly dependent), due to the non self adjoint nature of
the Djrbashian-Caputo derivative operator $-\DDC 0 \beta$. Nonetheless, according to the discussions in Section \ref{ssec:eig}, the
eigenvalues $\{\lambda_j\}$ increase to infinity with the index $j$, and asymptotically lies on two rays. Hence, one naturally
expects that the backward problem is also exponentially ill-posed. However, the magnitudes (and the real parts) of the eigenvalues
grow at a rate slower than that of the standard Sturm-Liouville problem, and thus the space fractional backward problem
is less ill-posed than the classical one. To illustrate the point, we present the numerical results
in Fig. \ref{fig:backward_spacefrac}. For all fractional orders $\beta$, the singular values decay exponentially,
but the decay rate increases dramatically with the increase of the fractional order $\beta$ and the terminal
time $T$. Hence, anomalous superdiffusion does not change the exponentially ill-posed nature of the
backward problem, but numerically it does enable recovering more Fourier modes of the initial data $v$.

Last, we note that for other choices of the fractional derivative, e.g., the Riemann-Liouville fractional derivative and
the fractional Laplacian \cite{BlumenthalGetoor:1959,Kwasnicki:2012}, the magnitude of eigenvalues of the operator also
tends to infinity, and the growth rate increases with the fractional order $\beta$. Therefore, the preceding observations
on the space fractional backward problem are expected to be valid for these choices as well.

\subsection{Sideways problem}
Last we return to the classical sideways diffusion problem but now with a fractional derivative in space rather than in
time. Let $\Omega=(0,1)$ be the unit interval. Then the one-dimensional space fractional diffusion equation is given by
\begin{equation*}
  u_t - {\DDC 0 \beta} u = 0,\quad (x,t)\in \Omega\times (0,\infty),
\end{equation*}
where the fractional order $\beta\in(1,2)$. The equation is equipped with an initial condition
$u(x,0)=0$ and the following lateral Cauchy boundary conditions
\begin{equation*}
  u(0,t) = f(t)\ \quad \mbox{and}\quad  \ u_x(0,t) = g(t),\ \ t>0.
\end{equation*}
We wish to compute the solution at $x=1$, i.e., $h(t):=u(1,t)$. In the case $\beta=2$, the model recovers
the standard diffusion equation, and we have already discussed the severe ill-conditioning of the classical
case. Due to the nonlocal nature of the fractional derivative, one might expect that in the space fractional
case, the sideways problem is less ill-posed. To see this, we take Laplace transform in time to arrive at
(with $\ \widehat{} \ $ denoting the Laplace transform)
\begin{equation*}
  z \widehat{u}(x,z) - {\DDC 0 \beta } \widehat{u}(x,z) = 0,
\end{equation*}
with the initial conditions (at $x=0$)
\begin{equation*}
  \widehat{u}(0,z) = \widehat{f}(z)\quad \mbox{and}\quad \widehat{u}_x(0,z) = \widehat{g}(z).
\end{equation*}
The solution $\widehat{u}(x,z)$ to the initial value problem is given by
\begin{equation*}
  \widehat{u}(x,z) = \widehat{f}(z)E_{\beta,1}(zx^\beta) + \widehat{g}(z)xE_{\beta,2}(zx^\beta)
\end{equation*}
and thus
\begin{equation*}
  \widehat{h}(z) = \widehat{f}(z)E_{\beta,1}(z) + \widehat{g}(z)E_{\beta,2}(z).
\end{equation*}
Like before, the boundary condition $h(t)$ at $x=1$ can be found by an inverse Laplace transform
\begin{equation*}
  h(t) = \int_{Br} e^{zt}\widehat{h}(z) dz.
\end{equation*}
In case of $\beta=2$, this gives $\cosh \sqrt{z}$ and $\sinh \sqrt{z}/\sqrt{z}$ multipliers to the data
$\widehat{f}(z)$ and $\widehat{g}(z)$ resulting in the exponential ill-conditioning of the sideways
heat problem. In the case of a general $\beta\in(1,2)$, the exponential asymptotics in Lemma \ref{lem:mitlef}
indicates that the problem still suffers from exponentially growing multipliers to the data, and thus the
problem is still severely ill-conditioned. Simple computation shows that the multiplier is
asymptotically larger for the fractional order $\beta$ closer to unity. In other words, anomalous diffusion
in space does not mitigate the ill-conditioned nature of the sideways problem, but actually worsens the
conditioning severely.

\begin{figure}[hbt!]
  \centering
  \begin{tabular}{cc}
     \includegraphics[trim = .5cm .1cm 1cm 0cm, clip = true, width=7.5cm]{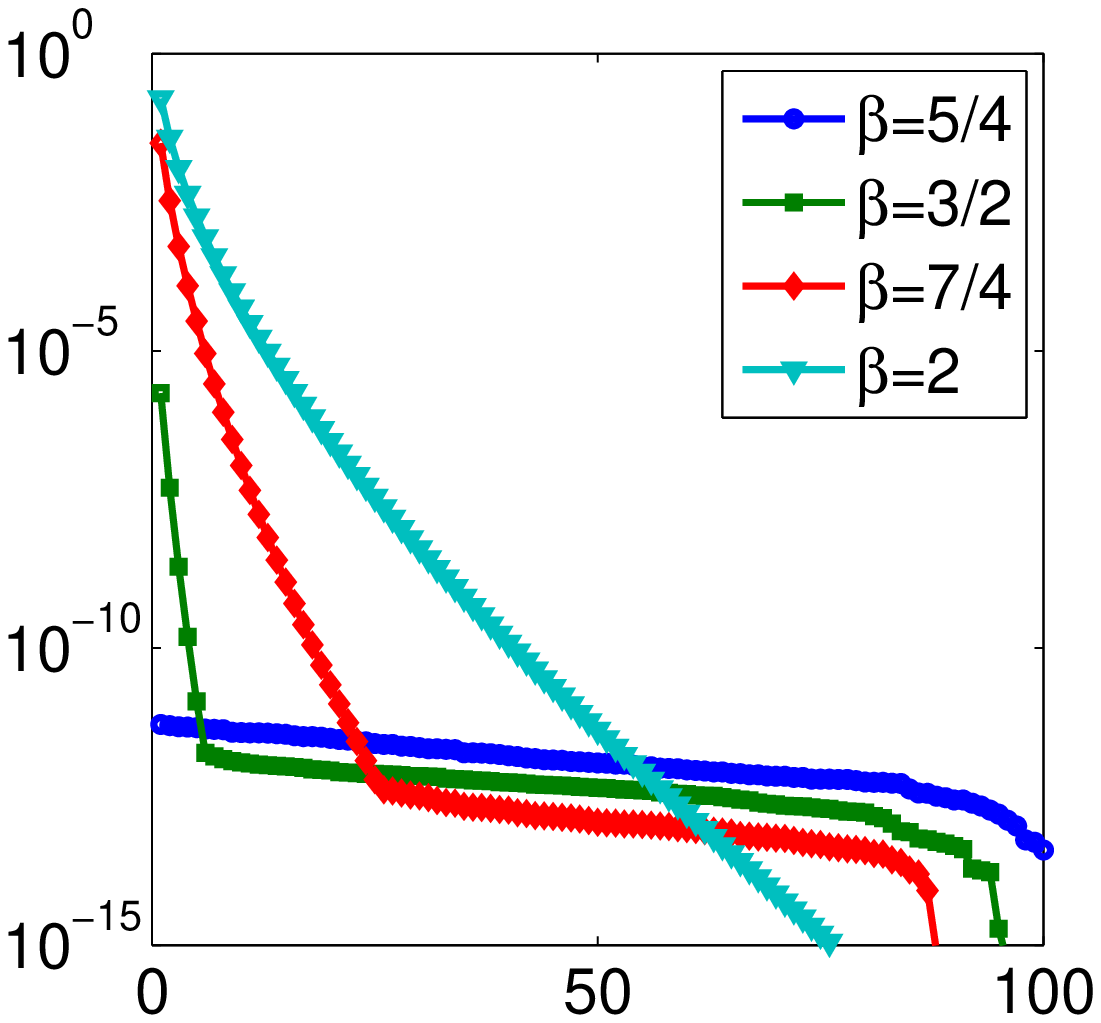} &
     \includegraphics[trim = .5cm .1cm 1cm 0cm, clip = true, width=7.5cm]{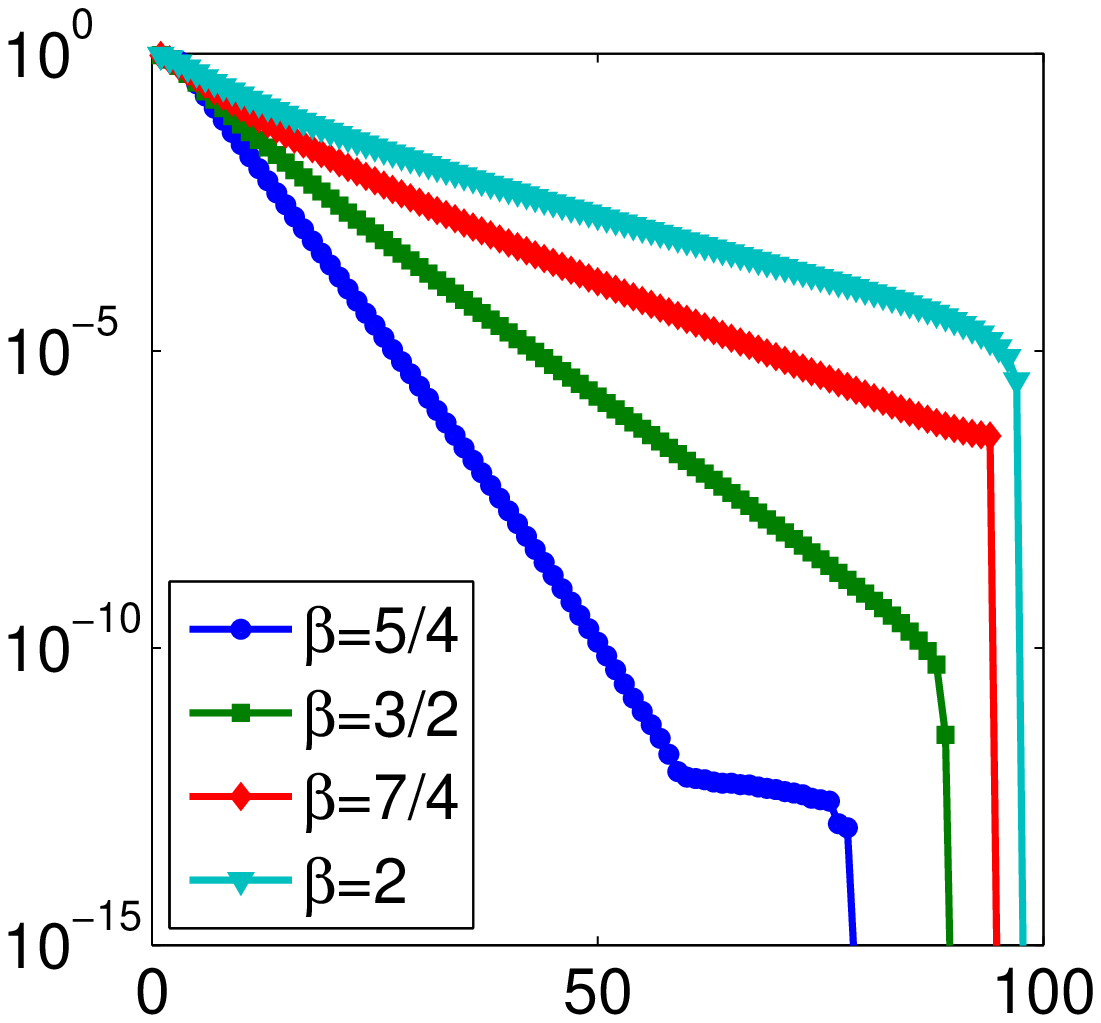}\\
     (a) $T=0.1$ & (b) $T=1$
  \end{tabular}
  \caption{Singular value spectrum of the forward map $F$ at times $T=0.1$ and $T=1$, for
  the sideways problem with Cauchy data at $x=0$.}\label{fig:sideways:space:left}
\end{figure}

To further illustrate the point, we compute the forward map $F$ from the Dirichlet boundary condition at
$x=1$ to the flux at $x=0$ numerically with a finite element in space/finite difference
in time scheme, cf. Appendix \ref{app:spacefrac} for the details. The numerical results are presented in Fig.
\ref{fig:sideways:space:left}. The singular value spectra clearly show the ill-posedness nature of the
space fractional sideways problem: as the fractional order $\beta$ increases from one to two, the
majority of the singular values move upward, the decay of the singular values slows down, and thus the
sideways problem becomes less and less ill-posed (but still severely so). Further, there are more tiny
singular values kicking in as the fractional order $\beta$ decreases to one, which indicates the inherent
rank deficiency of the forward map $F$ and might be relevant in the uniqueness of the inverse problem.
This confirms the preceding analysis: the degree of ill-posedness worsens with the decrease of the
fractional order $\beta$, and the fractional counterpart is more ill-posed than the classical one. In
other words, anomalous diffusion actually severely worsens the conditioning of the already very
ill-posed sideways problem. Further, the numerical results tend to indicate that the Djrbashian-Caputo
derivative with an order $\beta\in (1,2)$ acts as an interpolation between the diffusion and convection, which results in a history
mechanism in space: when the history piece runs from the left to the right, it is unlikely to
propagate the information in the reverse direction; and the closer is the fractional order $\beta$ to
unity, the stronger is the directional effect. The latter is not counterintuitive, since in the limit
of $\beta=1$, the Djrbashian-Caputo fractional derivative $\DDC 0 \beta u$ recovers the first order derivative $\frac{\partial u}{\partial x}$,
and the problem is of convection type, and surely no information can be convected backwards!

In the case $\beta=2$, one may equally measure the lateral Cauchy data at $x=1$, and aims at recovering the
Dirichlet boundary condition at $x=0$. Clearly, this does not change the nature of the inverse problem, and
it is equally ill-posed. Due to the directional nature of the Djrbashian-Caputo derivative $\DDC 0 \beta$,
one naturally wonders whether this ``directional'' feature does influence the ill-posed nature of the sideways
problem. To illustrate the point, we repeat the preceding arguments and deduce
\begin{equation*}
  z \widehat{u}(x,z) - {\DDC 0 \beta} \widehat{u}(x,z) = 0,
\end{equation*}
with the boundary conditions at $x=1$
\begin{equation*}
  \widehat{u}(1,z) = \widehat{f}(z) \quad \mbox{and}\quad \widehat{u}_x(1,z) = \widehat{g}(z).
\end{equation*}
To derive the solution, denote the initial conditions at $x=0$ by $\tilde f(z)=\widehat{u}(0,z)$ and
$\tilde g(z)=\widehat{u}_x(0,z)$. Then the solution $\widehat{u}(x,z)$ to the initial value problem is given by
\begin{equation*}
  \widehat{u}(x,z) = \tilde f(z)E_{\beta,1}(zx^\beta) + \tilde{g}xE_{\beta,2}(zx^\beta).
\end{equation*}
Use the differentiation formula $ \frac{d}{dx} x^{\gamma-1}E_{\beta,\gamma}(z x^\beta) =
z x^{\gamma-2}E_{\beta,\gamma-1}(z x^\beta)$ \cite[page 46]{KilbasSrivastavaTrujillo:2006}, we deduce that at $x=1$, there hold
\begin{equation*}
  \begin{aligned}
    \tilde f(z)E_{\beta,1}(z) + \tilde gE_{\beta,2}(z) &= \widehat f(z),\\
    \tilde f(z)E_{\beta,0}(z) + \tilde gE_{\beta,1}(z) &= z^{-1} \widehat g(z).
  \end{aligned}
\end{equation*}
Solving the linear system yields the solution to the sideways problem
\begin{equation*}
  \tilde f(z) = \frac{E_{\beta,1}(z)\widehat f(z) - z^{-1}E_{\beta,2}(z)\widehat g(z)}{
    E_{\beta,1}(z)^2 -E_{\beta,0}(z)E_{\beta,2}(z)},
\end{equation*}
and accordingly the solution $h(t)\equiv u(0,t)$ is given by an inverse Laplace transform.
The growth factors of the data $\widehat f$ and $\widehat g$ are $E_{\beta,1}(z)/(E_{\beta,1}(z)^2-
E_{\beta,0}(z)E_{\beta,2}(z))$ and $z^{-1}E_{\beta,2}(z)/(E_{\beta,1}(z)^2-E_{\beta,0}(z)E_{\beta,2}
(z))$, respectively. The growth of these factors at large $z$ argument determines the degree of
ill-conditioning of the sideways problem. To this end, we appeal to the exponential asymptotic of the
Mittag-Leffler function $E_{\alpha,\beta}(z)$, cf. Lemma \ref{lem:mitlef}, and note that Bromwhich
path lies in the sector $|\arg z|\leq \pi/2$ to deduce that for large $|z|$, there holds
\begin{equation*}
  \begin{aligned}
    E_{\beta,1}(z)^2 &\sim \frac{1}{\beta^2}e^{2z^{1/\beta}} - \frac{2}{\beta\Gamma(1-\beta)z} e^{z^{1/\beta}},\\
    E_{\beta,0}(z)E_{\beta,2}(z) & \sim \frac{1}{\beta^2}e^{2z^{1/\beta}} - \frac{1}{\beta\Gamma(2-\beta)}z^{1/\beta-1}e^{z^{1/\beta}}.
  \end{aligned}
\end{equation*}
Hence, the numerator $E_{\beta,1}(z)^2-E_{\beta,0}(z)E_{\beta,2}(z)$ behaves like
\begin{equation*}
  E_{\beta,1}(z)^2 - E_{\beta,0}(z)E_{\beta,2}(z) \sim \frac{1}{\beta\Gamma(2-\beta)}z^{1/\beta-1}e^{z^{1/\beta}}\qquad \mbox{as } |z|\to\infty.
\end{equation*}
This together with the exponential asymptotic of $E_{\beta,1}(z)$ and $E_{\beta,2}(z)$ from Lemma
\ref{lem:mitlef} indicates
that the multipliers for $\widehat f $ and $\widehat g$ are growing at most at a very low-order
polynomial rate, for large $z$. Hence, the high-frequency components of the data noise are not amplified
much (at most polynomially instead of exponentially). The analysis indicates that the sideways problem with
the lateral Cauchy data specified at the point $x=1$ is nearly well-posed, as long as the fractional
order $\beta$ is away from two, for which it recovers the classical ill-posed sideways problem for the
heat equation.

\begin{figure}[hbt!]
  \centering
  \begin{tabular}{cc}
     \includegraphics[trim = .5cm .1cm 1cm 0cm, clip = true, width=7.5cm]{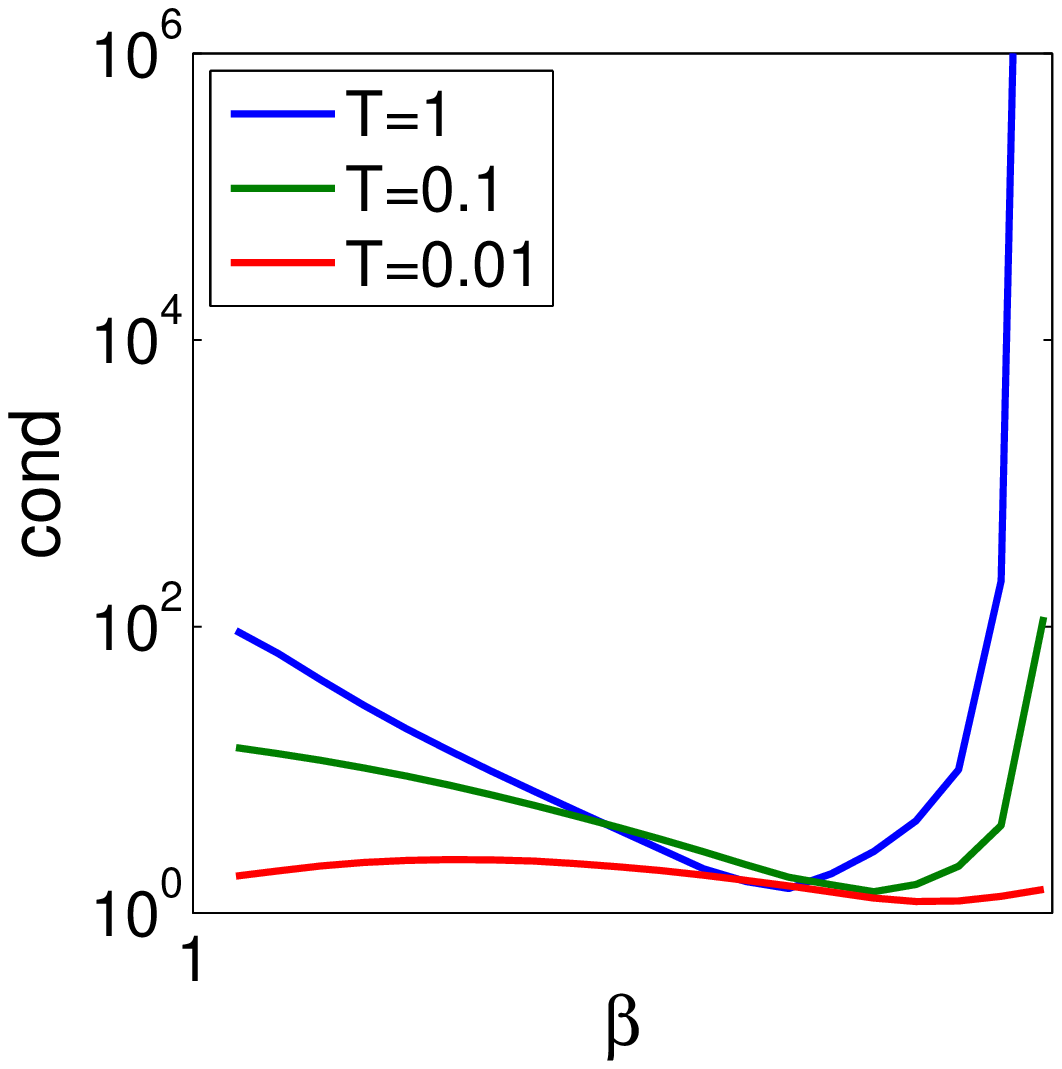} &
     \includegraphics[trim = .5cm .1cm 1cm 0cm, clip = true, width=7.5cm]{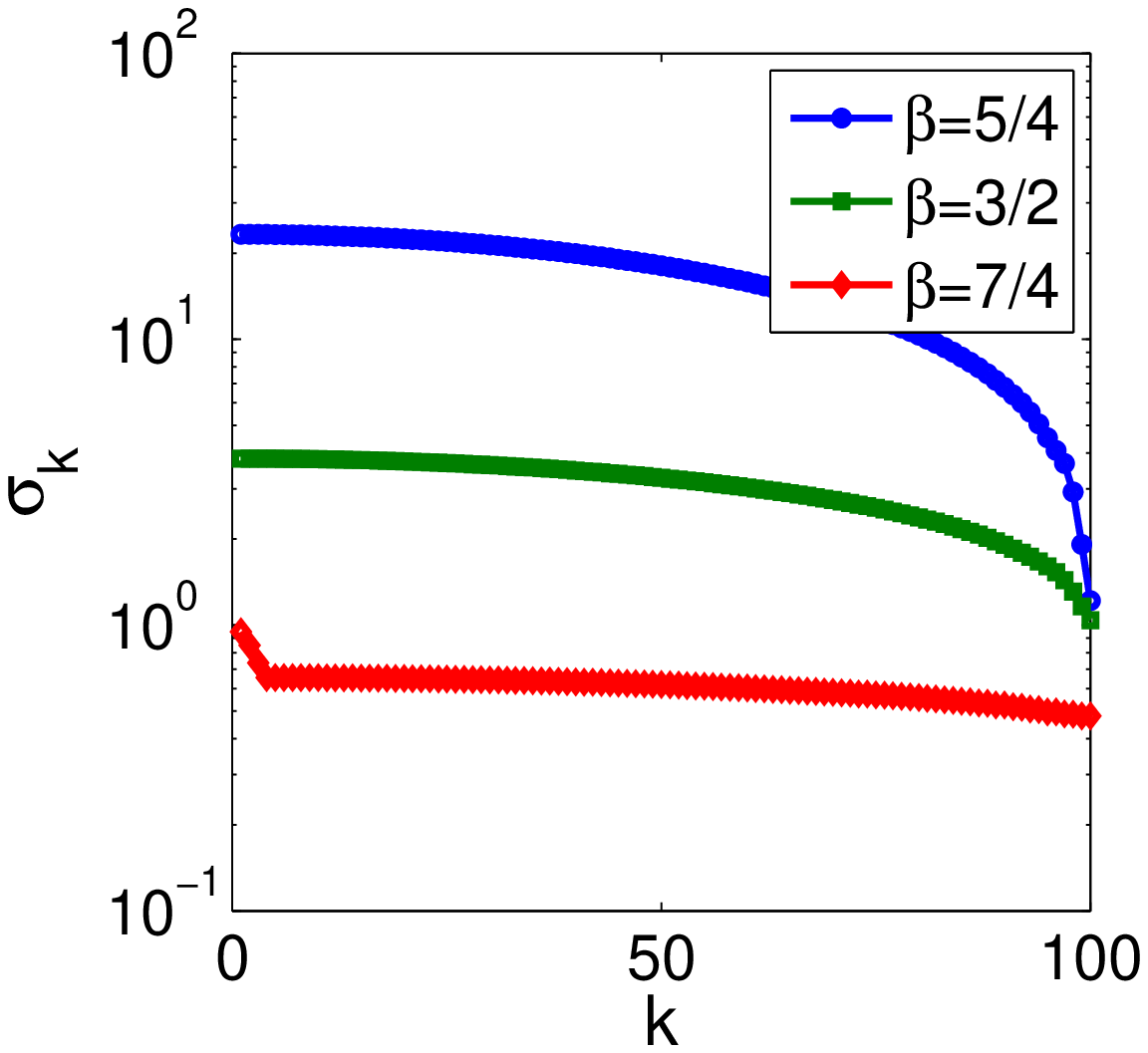}\\
     (a) condition number & (b) $T=1$
  \end{tabular}
  \caption{Numerical results for the space fractional sideways problem, with the lateral Cauchy data at the point $x=1$: (a) the condition number v.s. the fractional order $\beta$ and (b) the singular value spectrum at $T=1$.}\label{fig:sideways:space:right}
\end{figure}

Next we illustrate the preceding discussions numerically. The behavior of the forward map $F$ from the Dirichlet
boundary at $x=0$ to the flux data at $x=1$ is shown in Fig. \ref{fig:sideways:space:right}. For a wide range of
values of the fractional order $\beta$, the condition number of the forward map $F$ is of order $100$, which is fairly mild, in view of
the size of the linear system, i.e., $100\times 100$. When the fractional order $\beta$ increases towards two, the
inverse problem recovers the classical sideways problem, and as expected, the condition number increases dramatically.
However, the onset of the blowup depends on the terminal time $T$: the smaller is the time $T$, the smaller
seems the onset value. The precise mechanism for this phenomenon is still unknown. For $\beta\leq 7/4$, the singular
value spectrum only spans a narrow interval, resulting in a very small condition number. Physically, like before,
this can be explained as the ``convective'' nature of the Djrbashian-Caputo fractional derivative: as the fractional order $\beta$ tends
to unity, the information at $x=0$ is transported to $x=1$, free from distortion, and thus the inverse problem is almost
well-posed. In summary, depending on the location of the over-specified data, anomalous superdiffusion can either help or aggravate the
conditioning of the sideways problem.

Last, we would like to note that the study of space fractional inverse problems, either theoretical or numerical, is
fairly scarce. This is partly attributed to the relatively poor understanding of forward problems for FDEs with a space
fractional derivative: There are only a few mathematical studies on one-dimensional space fractional diffusion, and no
mathematical study on multi-dimensional problems involving space fractional derivatives (of either Riemann-Liouville or
Caputo type). Nonetheless, our preliminary numerical
experiments show distinct new features for related inverse problems, which motivate their analytical studies.

\section{Concluding remarks}
Anomalous diffusion processes arise in many disciplines, and the physics behind is very different
from normal diffusion. The unusual physics greatly influences the behavior of related forward problems.
Further, it is well known that backward fractional diffusion is much less ill-posed than the classical
backward diffusion, which has contributed to the belief that inverse problems for anomalous diffusion are
always better behaved than that for the normal diffusion. In this work we have examined several exemplary
inverse problems for anomalous diffusion processes in a numerical and semi-analytical manner. These
include the sideways problem, backward problem, inverse source problem, inverse Sturm-Liouville problem and Cauchy problems.
Our findings indicate that anomalous diffusion can give rise to very unusual new features, but they only
partially confirm the belief: depending on the data and unknown, it may influence either positively
or negatively the degree of ill-posedness of the inverse problem.

The mathematical study of inverse problems in anomalous diffusion is still in its infancy. There are
only a few rigorous theoretical results on the uniqueness, existence and stability, which mostly
focus on the one-dimensional case, and there are
many more open problems awaiting investigations. The development of stable and efficient reconstruction
procedures is an active ongoing research topic. However, due to the nonlocality of the forward model,
the construction of efficient schemes and their rigorous numerical analysis remain very challenging.
This is especially true for space fractional FDEs, and there are almost no theoretical or rigorous
numerical studies.

\section*{Acknowledgements}

The authors are grateful for the anonymous referees for their constructive comments.
The research of both authors is partly supported by NSF Grant DMS-1319052.

\appendix
\section{Numerical methods for special functions and FDEs}
\subsection{Computation of the Mittag-Leffler and Wright functions}\label{app:mitlef}

Like many special functions, the efficient and accurate numerical computation of the Mittag-Leffler function
$E_{\alpha,\beta}(z)$ is delicate \cite{GorenfloLoutchkoLuchko:2002,GorenfloLoutchkoLuchko:2003,SeyboldHilfer:2008}.
An efficient algorithm relies on partitioning the complex plane $\mathbb{C}$ into different regions, where
different approximations, i.e., power series, integral representation and exponential asymptotic for small
values of the argument, intermediate values and large values, respectively, are used for efficient numerical
computation; see \cite{SeyboldHilfer:2008} for the some partition and error estimates. The special case of
the Mittag-Leffler function $E_{\alpha,\beta}(z)$ with a real argument $z\in \mathbb{R}$, which plays a
predominant role in time-fractional diffusion, can also be efficiently computed with the Laplace transform
and suitable quadrature rules \cite{GarrappaPopolizio:2013}.

The computation of the Wright function $W_{\rho,\mu}(z)$ is even more delicate. In theory, like before, it can
be computed using power series for small values of the argument and a known asymptotic formula for large values,
and for the intermediate case, values are obtained by using an integral representation \cite{Luchko:2008}. The
integral representation for the Wright function $W_{\rho,\mu}(z)$ for intermediate values in the case of interest
for the fundamental solution in one-dimension (where $\rho=-\alpha/2<0$, $0<\mu=1+\rho<1$ and $z=-x$, $x>0$) is given by
\begin{equation*}
W_{\rho,\mu}(-x) = \int_0^\infty K(x,\rho,\mu,r)\,dr\qquad
\end{equation*}
where the kernel $K(x,\rho,\mu,r)$ is given by
\begin{equation*}
K(x,\rho,\mu,r) =
r^{-\mu} e^{-r + x\cos(\pi\rho)r^{-\rho}}
\sin(x\sin(\pi\rho) r^{-\rho} +\pi\mu).
\end{equation*}
This is a singular kernel with a leading order $r^{-\mu}=r^{-1-\rho}$, with successive singular kernels of the form
$r^{-1-2\rho}$, $r^{-1-3\rho}$ etc., upon expanding the terms. Hence, a direct treatment via numerical quadrature is inefficient. A more efficient
approach is to use the change of variable $s = r^{-\rho}$, i.e., $r=s^{-1/\rho}$, and the transformed kernel is
\begin{equation*}
  \widetilde{K}(x,\rho,\mu,s) = (-\rho)^{-1} s^{(-\rho)^{-1}(-\mu+1)-1}e^{-s^{-\frac{1}{\rho}}+x\cos(\pi\rho)s}\sin(x\sin(\pi\rho)s+\pi\mu).
\end{equation*}
The fundamental solution of the one-dimensional time-fractional diffusion equation is expressed in terms of a Wright function
$W_{\rho,\mu}(-x)$ with the choice $\rho=-\alpha/2$ and $\mu = 1 +\rho$, cf. \eqref{eqn:fundsol}. In this case
the resulting kernel $\tilde K$ simplifies to
\begin{equation*}
  \widetilde{K}(x,\rho,1+\rho,s) = (-\rho)^{-1}e^{-s^{-\frac{1}{\rho}}+x\cos(\pi\rho)s}\sin(x\sin(\pi\rho)s+(1+\rho)\pi).
\end{equation*}
This kernel is free from the grave singularity, and thus the quadrature method is quite effective.
In general, the integral can be computed efficiently via the Gauss-Jacobi quadrature, with the weight function
$s^{(-\rho)^{-1}(-\mu+1)-1}$. We note that an algorithm for the Wright function $W_{\rho,\mu}(z)$ over the whole complex
plane $\mathbb{C}$ with rigorous error analysis is still missing. The endeavor in this direction would almost certainly involve
dividing the complex domain $\mathbb{C}$ into different regions, and using different approximations on each region separately.

\subsection{Time fractional diffusion}\label{app:timefrac}
We describe a finite difference method for the initial boundary value problem for the one-dimensional
time-fractional diffusion equation
\begin{equation*}
    \partial_t^\alpha u - u_{xx} + qu = f(x,t) \quad (x,t)\in \Omega\times (0,T),
\end{equation*}
with the initial condition $u(x,0)=v$ and boundary conditions
\begin{equation*}
  u(0,t)=g(t)\quad \mbox{and}\quad u(1,t)=h(t),\quad t>0.
\end{equation*}
There are many efficient numerical schemes for discretizing the problem. The discretization in space can
be achieved by the standard central difference scheme, Galerkin finite element method \cite{JinLazarovZhou:2013}
or spectral method, and the discretization in time can be achieved with the L1 approximation \cite{SunWu:2006,
LinXu:2007} and convolution quadrature \cite{JinLazarovZhou:2014a}. We shall adopt the L1 approximation in time
and the central difference in space. Specifically, we divide the interval $[0,T]$ into uniform
subintervals, with nodes $t_k=k\tau$, $k=0,\ldots,K$, and a time step size $\tau=T/K$. Similarly, we partition
the spatial domain $\Omega$ into uniform subintervals, with grid points $x_i=ih$, $i=0,\ldots,N$, and mesh
size $h=1/N$. Then the L1-approximation of the Djrbashian-Caputo fractional derivative $\partial_t^\alpha u(x,t_k)$ developed in
\cite{SunWu:2006,LinXu:2007} is given by:
 \begin{equation}\label{eqn:L1approx}
   \begin{aligned}
     \partial_t^\alpha u(x,t_k)
     &\approx \tau^{-\alpha} \left[b_0u(x,t_k)-b_{k-1}u(x,t_0)+\sum_{j=1}^{k-1}(b_j-b_{j-1})u(x,t_{k-j})\right] 
   \end{aligned}
 \end{equation}
 where the weights $b_j$ are given by
\begin{equation*}
 b_j=((j+1)^{1-\alpha}-j^{1-\alpha})/\Gamma(2-\alpha),\ j=0,1,\ldots,K-1.
\end{equation*}
If the solution $u(x,t)$ is $C^2$ continuous in time, the local truncation error of the L1 approximation
is bounded by $c\tau^{2-\alpha}$ for some $c$ depending only on $u$ \cite[equation (3.3)]{LinXu:2007}. In general,
one can show that the scheme is only first-order accurate. Next with the central difference scheme in space and
the notation $u_i^k\approx u(x_i,t_k)$, we arrive at the following fully discrete scheme
\begin{equation*}
  \left[b_0u_i^k-b_{k-1}u_i^0+\sum_{j=1}^{k-1}(b_j-b_{j-1})u_i^{k-j}\right] + \frac{u_{i-1}^k-2u_i^k+u_{i+1}^k}{h^2} + q_iu_i^k = f_i^k, \ \ i =1,\ldots,N-1,
\end{equation*}
with $ q_i=q(x_i)$ and $f_i^k=f(x_i,t_k)$. We note that at each time step, one needs to solve a tridiagonal
linear system. However, the right hand side at the current step involves all previous steps, which can be
quite expensive for a small step size,  and this will most likely be required due to the first order in
time convergence. This history piece represents one of the main computational challenges for time fractional
differential equations. There are high-order schemes, e.g., convolution quadrature generated by
the second-order backward difference
formula \cite{JinLazarovZhou:2014a}. Further, we note that the finite difference scheme in space can be replaced
with the Galerkin finite element method, which is especially suitable for high dimensional problems on
a general domain and elliptic operator involving variable coefficients \cite{JinLazarovZhou:2013}.

\subsection{Space fractional diffusion} \label{app:spacefrac} Now we describe a fully discrete scheme based on the backward Euler
method in time and a Galerkin finite element method in space for the space fractional diffusion problem on the unit interval $\Omega=(0,1)$
\begin{equation*}
  u_t - {\DDC 0 \beta} u + qu = f\quad \mbox{ in }\Omega\times(0,T],
\end{equation*}
with the initial condition $u=v$ and the Dirichlet boundary condition
\begin{equation*}
  u(0,t) = g(t) \quad \mbox{and} \quad u(1,t)= h(t), \ \ t>0.
\end{equation*}
The Galerkin finite element method relies on the variational formulation for the fractional elliptic problem
\begin{equation*}
  - { \DDC 0 \beta } u + qu = f\quad \mbox{ in } \Omega,
\end{equation*}
with a homogeneous Dirichlet boundary condition $u(0)=u(1)=0$, recently developed in \cite{JinLazarovPasciak:2013a}.
The variational formulation of the problem is given by: find $u\in U\equiv H_0^{\beta/2}(\Omega)$ such that
\begin{equation*}
  -({\DDR 0 {\beta/2}} u,\ \DDR 1 {\beta/2} v) + (qu,v) = (f,v)\quad \forall v\in V,
\end{equation*}
where $\DDR 0 \gamma v$ and $\DDR 1 \gamma v$ are the left-sided and right-sided Riemann-Liouville derivative
of order $\gamma\in(0,1)$ defined by (with $c_\gamma=1/\Gamma(1-\gamma)$)
\begin{equation*}
  \DDR 0 \gamma v (x) = c_\gamma\frac{d}{dx}\int_0^x(x-s)^{-\gamma}v(s)ds\quad \mbox{and}\quad
  \DDR 1 \gamma v (x) = -c_\gamma\frac{d}{dx}\int_x^1(x-s)^{-\gamma}v(s)ds,
\end{equation*}
respectively, and the test space $V$ is given by
\begin{equation*}
  V = \left\{v\in H^{\beta/2}(\Omega):\ v(1) = 0,\ (x^{1-\beta},v)=0\right\}.
\end{equation*}
Now for the finite element discretization, we first divide the unit interval $\Omega$ into a uniform mesh, with the
grid points $x_i=ih$, $i=0,\ldots,N$ and mesh size $h=1/N$. Then for $U_h\subset U$ we take the continuous piecewise
linear finite element space, and for $V_h\subset V$ we construct it from $U_h$. Specifically, with the finite
element basis $\phi_i(x)$, $i=1,\ldots,N-1$, we take $\tilde\phi_i(x) = \phi_i(x) -\gamma_i(1-x) \in V_h$, where
the constant $\gamma_i$ is determined by the integral condition $(x^{1-\beta},\tilde\phi_i)=0$, i.e., $\gamma_i=
h^{2-\beta}((i-1)^{3-\beta}+(i+1)^{3-\beta}-2i^{3-\beta})$. The computation of the leading term in the stiffness
matrix and mass matrix can be carried analytically, and the part involving the potential $q$ can be computed efficiently
using quadrature rules; see \cite{JinLazarovPasciakRundell:2014} for details.

Now for the time-dependent problem, like before, we divide the time interval $[0,T]$ into uniform subintervals, with
$t_k=k\tau$, $k=0,\ldots,K$, and the time step size $\tau=T/K$. Then with the backward Euler method in time, and the
finite element method in space, the approximate solution $u_h^k$ at time $t_k$ can be split into $u_h^k=\tilde{u}_h^k
+ s^k$ with the particular solution $s^k=g(t_k) + (h(t_k)-g(t_k))x$ with the homogeneous solution $\tilde{u}_h^k\in U_h$ satisfying
\begin{equation*}
  \begin{aligned}
    \tau^{-1}(\tilde{u}_h^k,v_h) - ({\DDR 0 {\beta/2}}\tilde u_h^k,&\ {\DDR 1 {\beta/2}}v_h) + (q\tilde u_h^k,v_h) \\
    & =  (f^k,v_h) + \tau^{-1}(u_h^{k-1}-s^k,v_h) - (qs^k,v_h) \quad \forall v_h\in V_h.
  \end{aligned}
\end{equation*}
We note the resulting linear system is of lower Hessenberg form, due to the nonlocality of the fractional derivative
operator. However, the coefficient matrix does not change during the time stepping procedure, and
thus an LU factorization might be applied to speedup the computation.

\bibliographystyle{abbrv}
\bibliography{frac}

\end{document}